\documentclass[11pt]{article}

\usepackage{fullpage}

\usepackage{graphicx}
\usepackage{pgfplots}
\usetikzlibrary{decorations.pathreplacing}

\usepackage{amsmath,amssymb,amsthm,enumitem,mathtools}

\usepackage{multirow}
\usepackage{algorithm}\usepackage[noend]{algpseudocode} \makeatletter
\newlength{\AlgWrapKey}
\newcommand{\AlgWrap}[2]{\settowidth{\AlgWrapKey}{#1}#1\parbox[t]{\dimexpr\linewidth-\ALG@thistlm-\AlgWrapKey\relax}{#2\strut}}
\makeatother

\usepackage{booktabs}   \usepackage{adjustbox}  \usepackage{colortbl}   \definecolor{bestblue}{RGB}{226,237,249}
\definecolor{besttext}{RGB}{31,78,140}
\newcommand{\best}[1]{\hspace{-2.5pt}\tikz[baseline=(b.base)]\node[fill=bestblue, rounded corners=2.5pt, inner xsep=2.5pt, inner ysep=1.8pt, outer sep=0pt](b){\textcolor{besttext}{#1}};\hspace{-2.5pt}}
\usepackage{pifont}     \newcommand{\cmark}{\ding{51}}  \newcommand{\xmark}{\ding{55}}  

\usepackage[colorlinks=true, citecolor=cyan, linkcolor=cyan, urlcolor=cyan]{hyperref}

\usepackage{subcaption}

\usepackage[capitalize]{cleveref}
\Crefname{algorithm}{Algorithm}{Algorithms}

\usepackage[xindy, acronyms]{glossaries}   \makeglossaries

\allowdisplaybreaks

\let\oldparagraph\paragraph
\renewcommand{\paragraph}[1]{\oldparagraph{{{#1.}}}}

\newcommand{\ie}{{\it i.e.}}

\newcommand{\ones}{\mathbf 1}
\newcommand{\reals}{{\bf R}}

\newcommand{\symm}{{\bf S}}

\newcommand{\identity}{I}
\newcommand{\nullspace}{{\mathcal N}}
\newcommand{\range}{{\mathcal R}}

\newcommand{\Tr}{\mathop{\bf Tr}}
\newcommand{\diag}{\mathop{\bf diag}}
\newcommand{\lambdamax}{{\lambda_{\rm max}}}
\newcommand{\lambdamin}{\lambda_{\rm min}}

 \newcommand{\dist}{{\mathop{\bf dist}}}
\newcommand{\spann}{{\mathop{\bf span}}}
\newcommand{\argmin}{\mathop{\rm argmin}}

\newcommand{\nnz}{\mathop{\bf nnz}}

\newtheorem{theorem}{Theorem}[section]  \newtheorem{lemma}{Lemma}[section]  \newtheorem{corollary}{Corollary}[theorem]  \newtheorem{proposition}{Proposition}[section]
\newtheorem{assumption}{Assumption}[section]

\theoremstyle{definition}
\newtheorem{definition}{Definition}[section]

\theoremstyle{remark}

\newtheorem{remark}{Remark}[section]

\newacronym{LO}{LO}{linear optimization}
\newacronym{QO}{QO}{quadratic optimization}
\newacronym{MIQO}{MIQO}{mixed-integer quadratic optimization}
\newacronym{MIO}{MIO}{mixed-integer optimization}
\newacronym{MILO}{MILO}{mixed-integer linear optimization}
\newacronym{MINLO}{MINLO}{mixed-integer nonlinear optimization}
\newacronym{sBB}{sBB}{spatial branch and bound}
\newacronym{NLO}{NLO}{nonlinear optimization}
\newacronym{PWA}{PWA}{piecewise affine}
\newacronym{SVM}{SVM}{support vector machines}
\newacronym{ReLU}{ReLU}{rectified linear unit}
\newacronym{CPU}{CPU}{central processing unit}
\newacronym{GPU}{GPU}{graphics processing unit}
\newacronym{MPC}{MPC}{model predictive control}
\newacronym{ADMM}{ADMM}{alternating direction method of multipliers}
\newacronym{LP}{LP}{linear program}
\newacronym{QP}{QP}{quadratic program}
\newacronym{PDHG}{PDHG}{primal-dual hybrid gradient}
\newacronym{ADP}{ADP}{approximate dynamic programming}
\newacronym{FPGA}{FPGA}{field-programmable gate array}

\newcommand{\norm}[2][]{\left\| {#2} \right\|_{{#1}}}

\newcommand{\lagrangian}{\mathcal{L}}

\newcommand{\proj}{{\Pi}}
\newcommand{\opA}{{\mathbf{A}}}

\newcommand{\opH}{\mathbf{H}}

\newcommand{\opM}{\mathbf{M}}
\newcommand{\opN}{\mathbf{N}}

\newcommand{\opW}{\mathbf{W}}

\newcommand{\opT}{\mathbf{T}}
\newcommand{\opF}{\mathbf{F}}
\newcommand{\opTCV}{\opT_\mathrm{CV}}  

\newcommand{\xkkk}{x^{k+2}}
\newcommand{\xkk}{x^{k+1}}
\newcommand{\xk}{x^{k}}
\newcommand{\zkk}{z^{k+1}}
\newcommand{\zk}{z^{k}}

\newcommand{\ykk}{y^{k+1}}
\newcommand{\yk}{y^{k}}

\newcommand{\brho}{\boldsymbol{\rho}}
\newcommand{\balpha}{\boldsymbol{\alpha}}
\newcommand{\bc}{\boldsymbol{c}}

\newcommand{\xtkk}{\tilde{x}^{k+1}}

\newcommand{\cB}{\mathcal{B}}
\newcommand{\cC}{\mathcal{C}}

\newcommand{\cG}{\mathcal{G}}

\newcommand{\cN}{\mathcal{N}}

\newcommand{\cS}{\mathcal{S}}

\newcommand{\cX}{\mathcal{X}}

\definecolor{BrickRed}{cmyk}{0,0.89,0.94,0.28}  \colorlet{revcolor}{BrickRed}                   

\newcommand{\osqplf}{OSQP-LF}
\newcommand{\barc}{\operatorname{bar}}  \newcommand{\rec}{\operatorname{rec}}   

\newcommand{\hedgeSpeedup}{Against a single variant, the ensemble is up to five times faster on a benchmark set and never solves fewer instances (\Cref{sec:hedging-results}).}

\title{
    GPU-Accelerated Algorithm Hedging\\ 
    for Linear and Quadratic Programming
}
\author{
    \fontsize{13}{16}\selectfont
    Liubomir Maxim Baicev\textsuperscript{1}
    \enspace
    Jisun Park\textsuperscript{1}
    \enspace
    Ian McInerney\textsuperscript{2}
    \enspace
    Bartolomeo Stellato\textsuperscript{1}
    \\[1em]
    \small \textsuperscript{1}Department of Operations Research and Financial Engineering, Princeton University
    \\
    \small \textsuperscript{2}Department of Mechanical Engineering, Imperial College London
}

\begin{document}
\maketitle

\begin{abstract}
    We present a GPU-accelerated first-order solver for linear programs (LPs) and convex quadratic programs (QPs) that runs an ensemble of algorithm variants on the same problem as a single batch, and stops as soon as one of them converges.
    We call this \emph{algorithm hedging}, since no single variant is best on every problem and the ensemble protects the solver against a poor choice of variant.
    The variants differ in acceleration scheme, step-size rule, and relaxation parameter, but share the same problem matrices.
    On a GPU, all variants then share one sparse matrix product at a fraction of the cost of separate runs.
    The underlying method is the alternating direction method of multipliers, in which we replace the linear system solve with one preconditioned projected gradient step, so each iteration needs only matrix--vector products and simple projections.
    With Halpern iteration and restarts, the method converges linearly on problems with polyhedral constraints, and its iterates certify primal and dual infeasibility.
    The same batching also solves many instances of a parametric LP or QP at once, with several variants per instance.
    Our open-source implementation, \osqplf{}, is the new GPU mode of the OSQP solver.
    On five standard benchmark sets, MIPLIB and Mittelmann for LPs, and Maros--M\'esz\'aros, SuiteSparse, and QPLIB for QPs, it is the fastest GPU first-order solver in most cases, and hedging makes it up to five times faster than a single variant.
    In a portfolio backtest with thousands of QPs, warm starting and batching reduce the total solve time by one to two orders of magnitude compared with interior-point solvers.
\end{abstract}

\section{Introduction}\label{sec:intro}

We consider the convex quadratic program (QP)
\begin{equation}\label{prob:qp}
    \begin{array}{ll}
        \textnormal{minimize} & (1/2) x^T P x + q^T x \\
        \textnormal{subject to} & Ax \in \cC, \quad x \in \cX,
    \end{array}
    \tag{QP}
\end{equation}
where $x \in \reals^d$ is the optimization variable.
The problem data are a symmetric positive semidefinite matrix~$P \in \symm_+^d$, a vector~$q \in \reals^d$, a matrix~$A \in \reals^{m \times d}$, and nonempty \emph{boxes}~$\cC=[l,u] \subseteq \reals^m$ and~$\cX=[l_x,u_x] \subseteq \reals^d$.
The boxes collect the constraint bounds and the variable bounds, and every lower bound may equal~$-\infty$ and every upper bound may equal~$+\infty$.
If~$P = 0$, the problem is a linear program (LP).

Problems of this form arise in many fields.
In machine learning, they include support vector machines~\cite{cortesSupportvectorNetworks1995}, the Lasso~\cite{tibshiraniRegressionShrinkageSelection1996}, and Huber fitting~\cite{huberRobustEstimationLocation1964}.
In finance, portfolio optimization leads to QPs~\cite{markowitzPortfolioSelection1952,boydMultiPeriodTradingConvex2017}, and a backtest solves one of them for every trading day.
In control, model predictive control solves a QP at every time step~\cite{garciaModelPredictiveControl1989}.
LPs and QPs are also subproblems of other methods, such as the relaxations solved in branch and bound for mixed-integer programs~\cite{confortiIntegerProgramming2014}.
These applications call for solvers that work in two regimes.
In the first, a single instance has millions of variables and constraints.
In the second, many instances share the matrices~$P$ and~$A$ and differ only in the vector~$q$ and in the bounds.

On instances with millions of variables and constraints, the classical interior-point and simplex methods reach their limits.
Both rely on matrix factorizations, which become expensive in time and memory as the problem grows and are hard to parallelize~\cite{luCuPDLPjlGPUImplementation2024}.
First-order methods avoid most of this cost because their iterations are much cheaper, so recent solvers for large problems build on them.
Some of them, such as OSQP~\cite{stellatoOSQPOperatorSplitting2020} and SCS~\cite{odonoghueConicOptimizationOperator2016,odonoghueOperatorSplittingHomogeneous2021}, are based on the alternating direction method of multipliers (ADMM)~\cite{boydDistributedOptimizationStatistical2011}.
Each of their iterations solves a linear system with a fixed matrix and computes a projection.
Others, such as PDLP~\cite{applegatePDLPPracticalFirstOrder2025}, are based on the primal--dual hybrid gradient (PDHG) method~\cite{chambolleFirstOrderPrimalDualAlgorithm2011}, whose iterations need only matrix--vector products and projections.
These operations run efficiently on graphics processing units (GPUs).
As a result, GPU implementations of such methods are now competitive with commercial solvers on large LPs~\cite{luCuPDLPjlGPUImplementation2024}, later solvers improve on them further~\cite{luCuPDLPxFurtherEnhanced2025,chenHPRLPImplementationHPR2025}, and several have been extended to QPs~\cite{luPracticalOptimalFirstOrder2025,huangRestartedPrimalDualHybrid2025,chenHPRQPDualHalpern2025,luCuPDQPxGPUAcceleratedRestarted2026}.
Despite this progress, the current design of first-order solvers has three limitations.

\paragraph{Linear systems}
In ADMM-based solvers, the linear system is the main cost of an iteration.
A direct method factorizes the matrix of this system once and reuses the factors.
At large scale, however, the factorization is expensive in time and memory, and the triangular solves that follow it are hard to parallelize on a GPU~\cite{schubigerGPUAccelerationADMM2020}.
Indirect methods such as conjugate gradient replace the factorization with inner iterations that use only matrix--vector products~\cite{odonoghueConicOptimizationOperator2016,schubigerGPUAccelerationADMM2020,diamandisGeNIOSAlmostSecondorder2025}.
But each inner iteration costs as much as one step of a method without inner solves, and it only refines the solution of a subproblem that the next iteration replaces.
The inner solver also needs a stopping tolerance, which existing solvers set from the iteration count~\cite{odonoghueConicOptimizationOperator2016}, from the current residuals~\cite{schubigerGPUAccelerationADMM2020}, or from both~\cite{diamandisGeNIOSAlmostSecondorder2025}.
This tolerance is one more rule to tune.
In reported benchmarks, PDHG-based solvers outperform SCS with conjugate gradient on LPs~\cite{applegatePracticalLargescaleLinear2021} and on large QPs~\cite{huangRestartedPrimalDualHybrid2025}.
PDHG has no linear system to solve, because it is ADMM with a proximal term that cancels the matrix~$A^TA$ in the primal step~\cite{esserGeneralFrameworkClass2010,oconnorEquivalencePrimaldualHybrid2020}.
Its primal step is explicit for LPs.
For QPs, the step involves the matrix~$P$, and it is explicit only if~$P$ is diagonal.
First-order QP solvers therefore change the method, with inner iterations~\cite{huangRestartedPrimalDualHybrid2025,liPDHCGIIEnhancedVersion2026}, momentum~\cite{luPracticalOptimalFirstOrder2025}, an auxiliary variable with its own step size~\cite{luCuPDQPxGPUAcceleratedRestarted2026}, or a dual formulation with an additional block of variables~\cite{chenHPRQPDualHalpern2025}.

\paragraph{Batching on GPUs}
GPU solvers make little use of batching.
Most of them are built on the product of a sparse matrix with a single vector~\cite{luCuPDLPxFurtherEnhanced2025,chenHPRLPImplementationHPR2025,luCuPDQPxGPUAcceleratedRestarted2026}.
A GPU can multiply the same matrix by a block of vectors at little additional cost, since it reads the matrix from memory once for the whole block.
For example, Blin et al.~\cite{blinBatchedFirstOrderMethods2026} measure a nearly constant product time for blocks of~$32$ to~$512$ vectors on one LP.
In our experiments, an iteration with a batch of~$16$ vectors costs little more than an iteration with a single vector on instances with up to about~$10^5$ nonzeros (\Cref{sec:hedging-results}).
Deep learning relies on such batched operations~\cite{paszkePyTorchImperativeStyle2019}, and batching has recently entered mixed-integer programming~\cite{blinBatchedFirstOrderMethods2026,liuSequentialNodesGPU2026,lucasGPUacceleratedNonlinearBranchandBound2026} and model predictive control~\cite{hoseWarpMPCLargeBatchMPC2026}.
The LP and QP solvers that do batch use it to solve many instances at once with the same algorithm~\cite{gurungSimultaneousSolvingBatched2019,blinBatchedFirstOrderMethods2026,hoseWarpMPCLargeBatchMPC2026}.
To our knowledge, the batch has not been used to solve a single instance faster.

\paragraph{Parameter tuning}
The performance of a first-order method depends strongly on its parameters, such as step sizes, relaxation parameters, and restart rules~\cite{stellatoOSQPOperatorSplitting2020,ghadimiOptimalParameterSelection2015,giselssonLinearConvergenceMetric2017}.
Worst-case analysis gives optimal parameters for whole classes of problems~\cite{droriPerformanceFirstorderMethods2014,kimOptimizedFirstorderMethods2015,taylorSmoothStronglyConvex2017,dasguptaBranchandboundPerformanceEstimation2023}, but the best parameters for a given instance can be very different.
State-of-the-art solvers therefore adapt their parameters during the solve with restart tests and step-size updates whose constants are fixed in advance~\cite{applegatePDLPPracticalFirstOrder2025,luCuPDLPxFurtherEnhanced2025,chenHPRLPImplementationHPR2025,chenHPRQPDualHalpern2025}.
Learning-based methods show that tuning the parameters or the warm start to a family of similar instances can reduce the solve time considerably~\cite{ichnowskiAcceleratingQuadraticOptimization2021,sambharyaLearningWarmstartFixedpoint2024,sambharyaLearningAlgorithmHyperparameters2026,ranjanLearningDistributionallyRobust2026}, but they need training data from that family.
When no algorithm is best on every instance, a classical remedy is to run several algorithms at once and keep the first to finish~\cite{hubermanEconomicsApproachHard1997,gomesAlgorithmPortfolios2001}, as the concurrent modes of Gurobi~\cite{gurobioptimizationllcGurobiOptimizerReference2026} and cuOpt~\cite{nvidiaCuOptGPUacceleratedDecision2025} do with different LP algorithms.
To our knowledge, first-order solvers do not use this remedy across their own parameters.
For a general-purpose solver that faces a new instance, the best parameters thus remain unknown until the instance is solved.

In this paper, we address these limitations with a GPU-accelerated first-order method for LPs and QPs.
The method needs no linear system solve, and it handles the quadratic term with one gradient step, so the same iteration solves LPs and QPs.
It runs an ensemble of algorithm variants on each instance as a single batch, and stops as soon as one of them converges.

Our contributions are as follows.
\begin{itemize}
    \item \emph{Linear-system-free ADMM.}
    We propose \emph{linear-system-free ADMM}, a function-linearized proximal ADMM~\cite[Chapter~8]{ryuLargeScaleConvexOptimization2022} that replaces the linear system solve with one preconditioned projected gradient step on the augmented Lagrangian.
    Each iteration needs only matrix--vector products and projections, and the same iteration solves LPs and QPs: on LPs it coincides with PDHG, and on QPs it is a preconditioned Condat--V\~u iteration~\cite{condatPrimalDualSplitting2013,vuSplittingAlgorithmDual2013}.
    With restarted Halpern iteration, we prove linear convergence on problems with polyhedral constraints, and we show that the iterates certify primal and dual infeasibility.
    On equality-constrained LPs, the resulting complexity matches the lower bound for first-order methods up to constants and a logarithmic term.

    \item \emph{GPU batching and algorithm hedging.}
    Our GPU implementation runs many copies of the method as a single batch, in which one product of a sparse matrix with a block of vectors serves all copies.
    The copies can solve different instances that share the matrices~$P$ and~$A$.
    They can also run an \emph{ensemble} of \emph{algorithm variants} on a single instance, and stop as soon as one of them satisfies the termination criteria, which we call \emph{algorithm hedging}.
    The variants share a common structure that makes them batchable, and differ in acceleration scheme, step-size rule, and relaxation parameter.
    Each instance runs the first variants of a fixed order, which we compute offline by minimizing the total run time over perturbed benchmark instances.
    The solver chooses the number of variants for each instance at setup, from a model of the GPU and a few timed iterations.
    \hedgeSpeedup

    \item \emph{Numerical experiments.}
    We implement the method as \osqplf{}, a new linear-system-free mode of the OSQP solver.
    It is the fastest first-order solver on Maros--M\'esz\'aros, SuiteSparse, and the LP relaxations of MIPLIB, by factors between~$1.5$ and~$11$ over HPR-QP and cuPDLPx.
    On Mittelmann's LP benchmark it matches cuPDLPx and HPR-LP, and on QPLIB it is slower than HPR-QP at low accuracy and faster at high accuracy.
    Without presolve, it is also faster than the best interior-point solvers on every set except Maros--M\'esz\'aros.
    In a portfolio backtest with thousands of QPs, warm starting and batching reduce the total solve time by one to two orders of magnitude compared with interior-point solvers.
\end{itemize}

\osqplf{} is open source.
It is available in the development branch of OSQP and will be part of the OSQP~2.0 release.
The solver and its documentation are available at
\begin{center}
    \url{https://osqp.org},
\end{center}
and the code to reproduce all our experiments is available at
\begin{center}
    \url{https://github.com/osqp/osqp2-benchmarks}.
\end{center}

\subsection{Related work}\label{sec:related-work}

\paragraph{Proximal ADMM and PDHG}
ADMM~\cite{glowinskiLapproximationParElements1975,gabayDualAlgorithmSolution1976} is a proximal point method through Douglas--Rachford splitting~\cite{ecksteinDouglasRachfordSplitting1992}.
Semi-proximal ADMM~\cite{fazelHankelMatrixRank2013} adds a positive semidefinite proximal term to each subproblem.
With a suitable term, the quadratic terms cancel and the primal update becomes explicit, which is how our method avoids the linear system.
PDHG arises in the same way.
On LPs, linear-system-free ADMM is a linearized ADMM, which is PDHG once the splitting variable is eliminated~\cite{esserGeneralFrameworkClass2010}, \cite[\S8]{oconnorEquivalencePrimaldualHybrid2020}.
Conversely, Ma et al.~\cite{maUnderstandingConvergencePreconditioned2023} write preconditioned PDHG as an indefinite proximal ADMM.
\Cref{subsec:relation} gives the correspondence.
Chen et al.~\cite[Appendix~A.1]{chenHPRQPDualHalpern2025} also apply a proximal term of this kind to the primal QP, with a multiple of the identity in place of a preconditioner.
With Halpern iteration, the resulting method is linear-system-free ADMM with a scalar preconditioner and relaxation parameter~$2$.
On 15 Maros--M\'esz\'aros instances at tolerance~$10^{-8}$, they find it about two times slower than HPR-QP, which works on a dual formulation, and they attribute the difference to the large norm of the proximal term.
Our method differs in its diagonal preconditioner, which keeps the scale of each row of~$P$, and we give its convergence analysis on QPs.
In our experiments, the primal form with a diagonal preconditioner and hedging is faster than HPR-QP on the Maros--M\'esz\'aros and SuiteSparse sets (\Cref{sec:experiments}).

\paragraph{Acceleration via Halpern iteration and restarts}
Halpern iteration~\cite{halpernFixedPointsNonexpanding1967} pulls each step toward a fixed anchor point and reduces the fixed-point residual of a nonexpansive operator at the rate~$O(1/k)$~\cite{liederConvergenceRateHalperniteration2021}.
Sun et al.~\cite{sunAcceleratingPreconditionedADMM2025} write proximal ADMM as a degenerate preconditioned proximal point method~\cite{brediesDegeneratePreconditionedProximal2022} and prove~$O(1/k)$ rates for its Halpern-accelerated version.
HPR-LP~\cite{chenHPRLPImplementationHPR2025} and HPR-QP~\cite{chenHPRQPDualHalpern2025} build on this framework.
Restarts turn the~$O(1/k)$ rate into a linear one when the monitored quantity is \emph{sharp}, \ie, bounded below by a multiple of the distance to the solution set, as the normalized duality gap of an LP is~\cite{applegateFasterFirstorderPrimaldual2023}.
Lu and Yang~\cite{luRestartedHalpernPDHG2024} prove that the fixed-point residual of PDHG on an LP is sharp on bounded sets, so that restarted Halpern PDHG converges linearly even when the operator is reflected, as in cuPDLPx~\cite{luCuPDLPxFurtherEnhanced2025}.
Such sharpness follows from error bounds for polyhedral systems~\cite{robinsonContinuityPropertiesPolyhedral1981}.
Liu et al.~\cite{liuRestartedReflectedHalpern2026} derive it from local error bounds, so it holds only near the solution set.
Lu, Peng, and Yang~\cite{luCuPDQPxGPUAcceleratedRestarted2026} prove it for cuPDQPx on a bounded set that contains the iterates, and we do the same for linear-system-free ADMM.
For QPs, Lu and Yang~\cite{luPracticalOptimalFirstOrder2025} also prove lower bounds for span-respecting primal--dual methods, which PDQP matches up to a rescaling.

\paragraph{Infeasibility detection}
ADMM and PDHG detect infeasibility through the differences of successive iterates, which converge to certificates of primal or dual infeasibility~\cite{banjacInfeasibilityDetectionAlternating2019,applegateInfeasibilityDetectionPrimalDual2024}.
The limit of these differences is the minimal displacement vector of the operator.
Banjac~\cite{banjacMinimalDisplacementVector2021} splits this vector, for the Douglas--Rachford operator, orthogonally into primal and dual parts, and Jiang et al.~\cite{jiangRangeDisplacementOperator2026} characterize the range of the PDHG displacement.
Park and Ryu~\cite{parkAcceleratedInfeasibilityDetection2023} accelerate the convergence to the minimal displacement vector of any nonexpansive operator, and on LPs the fixed-point residual of restarted Halpern PDHG with fixed-frequency restarts converges to its limit linearly~\cite[Theorem~5]{luRestartedHalpernPDHG2024}.
A different route is a homogeneous self-dual embedding of the LP, to which H-PDLP~\cite{luHPDLPGPUAcceleratedLinear2026} applies a generalized PDHG, so that one run returns either a solution or a certificate.
We follow the first route, and test the primal and dual displacements of one step as certificates of dual and primal infeasibility, respectively.

\paragraph{GPU solvers and batching}
After cuPDLP~\cite{luCuPDLPjlGPUImplementation2024}, GPU first-order solvers for LPs include cuPDLPx~\cite{luCuPDLPxFurtherEnhanced2025}, HPR-LP~\cite{chenHPRLPImplementationHPR2025}, and the PDLP implementation of cuOpt~\cite{nvidiaCuOptGPUacceleratedDecision2025}, and PDNQP~\cite{chenPDNQPGPUbasedFactorizationfree2026} addresses nonconvex QPs.
Lu and Yang~\cite{luOverviewGPUbasedFirstOrder2025} survey GPU-based first-order methods for LPs and their extensions.
Batching has so far served to solve many instances at once.
Gurung and Ray~\cite{gurungSimultaneousSolvingBatched2019} run a dense simplex method on up to~$10^5$ small LPs.
Blin et al.~\cite{blinBatchedFirstOrderMethods2026} batch the LPs of strong branching and bound tightening, which share the constraint matrix, in a restarted Halpern PDHG solver.
WarpMPC~\cite{hoseWarpMPCLargeBatchMPC2026} batches model predictive control problems that share a sparsity pattern.

\paragraph{Concurrent solvers and hedging}
Several methods adapt the parameters of a first-order solver to the instance, but run one setting at a time.
ReLU-QP~\cite{bishopReLUQPGPUAcceleratedQuadratic2024}, a GPU ADMM solver for model predictive control, switches online among linear-system inverses precomputed for a fixed set of step sizes.
Lin et al.~\cite{linLearningOverRelaxationPolicies2026} train a policy offline that sets the relaxation parameter of OSQP during the solve.
Concurrent solvers instead run several algorithms at once.
The concurrent modes of Gurobi~\cite{gurobioptimizationllcGurobiOptimizerReference2026} and cuOpt~\cite{nvidiaCuOptGPUacceleratedDecision2025} run several LP algorithms, among simplex, interior-point, and PDHG methods, on one instance and keep the first to finish.
Our ensemble runs variants of one algorithm, which share every product with the problem matrices, as a single batch.
The order in which the variants enter the ensemble is fixed offline, separately for LPs and for QPs.

\paragraph{Outline}
\Cref{sec:lsf-osqp} derives linear-system-free ADMM and relates it to other first-order methods.
\Cref{sec:sharpness} analyzes its convergence, infeasibility detection, and termination criteria.
\Cref{sec:gpu-batching} describes the batched iteration and algorithm hedging, and~\Cref{sec:implementation} the implementation and its default settings.
\Cref{sec:experiments} reports the numerical experiments, and~\Cref{sec:conclusion} concludes.
The appendices contain the proofs, the algebra engine, the data of the experiments, and the detailed results.

\section{Linear-system-free ADMM}
\label{sec:lsf-osqp}

Throughout, we write the constraint~$Ax\in\cC$ of problem~\eqref{prob:qp} as~$Ax=z$ with~$z\in\cC$, so that the optimization variables are~$x\in\reals^d$ and~$z\in\reals^m$.

The dual problem of~\eqref{prob:qp} is
\begin{equation}
    \label{prob:qp-dual}
    \begin{array}{ll}
        \textnormal{maximize} & -(1/2)x^TPx - S_\cC(y) - S_\cX(w) \\
        \textnormal{subject to} & Px+q+A^Ty+w=0,
    \end{array}
    \tag{Dual-QP}
\end{equation}
where~$x\in\reals^d$ is the primal variable and~$y\in\reals^m$ and~$w\in\reals^d$ are the dual multipliers of~$Ax=z$ and of~$x\in\cX$, respectively.
The functions~$S_\cC(y)=\sup_{c\in\cC}y^Tc$ and~$S_\cX(w)=\sup_{x\in\cX}w^Tx$ are the support functions of~$\cC$ and~$\cX$.
For the box~$\cC$, the support function is~$S_\cC(y)=u^T(y)_++l^T(y)_-$, where~$(\cdot)_+=\max\{\cdot,0\}$ and~$(\cdot)_-=\min\{\cdot,0\}$ act entrywise and~$0\cdot(\pm\infty)=0$.
It is finite on the \emph{barrier cone}~$\barc\cC=\{ y\in\reals^m \mid S_\cC(y)<+\infty \}$~\cite{bauschkeConvexAnalysisMonotone2017}, \ie, when~$y_i\le0$ for every~$i$ with~$u_i=+\infty$ and~$y_i\ge0$ for every~$i$ with~$l_i=-\infty$.
The same holds for~$S_\cX$ and~$\barc\cX$.
Hence the dual problem implicitly restricts~$y$ and~$w$ to these cones.

The KKT optimality conditions of~\eqref{prob:qp} are
\begin{equation}
    \label{eq:qp-kkt}
    \begin{aligned}
        &Ax = z,
        \quad Px + q + A^T y + w = 0, \\
        &x \in \cX,
        \quad w \in \cN_\cX(x),
        \quad z \in \cC,
        \quad y \in \cN_{\cC}(z),
    \end{aligned}
\end{equation}
where~$\cN_\cX(x)$ and~$\cN_{\cC}(z)$ are the normal cones of~$\cX$ at~$x$ and of~$\cC$ at~$z$, respectively.
A point~$(x,z,y,w)$ satisfies these conditions if and only if~$(x,z)$ is feasible for~\eqref{prob:qp}, $(x,y,w)$ is feasible for~\eqref{prob:qp-dual}, and the two objective values are equal.
Hence a KKT point exists if and only if both problems have solutions and strong duality holds.
Since~$\cC$ and~$\cX$ are boxes, a KKT point exists whenever both problems are feasible~\cite{frankAlgorithmQuadraticProgramming1956}.

The rest of this section derives linear-system-free ADMM from ADMM, adds relaxation and Halpern iteration, and writes the step in primal--dual form.
It then relates the method to OSQP and to PDHG-based and HPR-based solvers.

\subsection{Derivation from ADMM}
\label{subsec:derivation}
The augmented Lagrangian of~\eqref{prob:qp} is
\begin{equation*}
    \lagrangian_\rho (x, z, y) =
        (1/2) x^T P x + q^T x
        + y^T (Ax - z)
        + (\rho/2) \norm{Ax - z}^2,
\end{equation*}
where~$y \in \reals^m$ is the multiplier of~$Ax = z$, $\rho > 0$ is the step size, and~$\norm{\cdot}$ is the Euclidean norm.
We consider the semi-proximal ADMM iteration~\cite{fazelHankelMatrixRank2013}
\begin{equation}\label{algo:prox-admm}
    \begin{array}{rcl}
        \xkk &=& \underset{x \in \cX}{\argmin} \left(
            \lagrangian_\rho (x, z^k, y^k)
            + (1/2) \norm[W]{x - x^k}^2
        \right) \\
        \zkk &=& \underset{z \in \cC}{\argmin} \, \lagrangian_\rho (x^{k+1}, z, y^k)
        = \proj_\cC\!\left(A\xkk + \rho^{-1}\yk\right) \\
        \ykk &=& \yk + \rho (A \xkk - \zkk),
    \end{array}
\end{equation}
where~$W \succeq 0$ is a proximal matrix, $\norm[W]{v}=(v^TWv)^{1/2}$, and~$\proj_\cC$ is the Euclidean projection onto~$\cC$.
With~$W = 0$, this is ADMM~\cite{boydDistributedOptimizationStatistical2011}, whose~$x$-update requires solving a linear system with the matrix~$P + \rho A^T A$ if~$\cX = \reals^d$, and a constrained QP otherwise.
We choose~$W$ to cancel this matrix.
Let~$R \succ 0$ be such that
\begin{equation*}
    W = R - (P + \rho A^T A) \succeq 0.
\end{equation*}
The objective of the~$x$-update is then a quadratic function of~$x$ with Hessian~$P + \rho A^T A + W = R$.
Its gradient at~$\xk$ is~$\nabla_x \lagrangian_\rho (\xk, \zk, \yk) = P\xk + q + A^T\big(\yk + \rho(A\xk - \zk)\big)$, because the gradient of the proximal term vanishes at~$\xk$.
Up to a constant, this objective therefore equals~$(1/2) \norm[R]{x - \hat x^k}^2$ with~$\hat x^k = \xk - R^{-1} \nabla_x \lagrangian_\rho (\xk, \zk, \yk)$, and the~$x$-update is the projected gradient step
\begin{equation*}
    \xkk = \proj_\cX^R \! \left( \xk - R^{-1} \nabla_x \lagrangian_\rho (\xk, \zk, \yk) \right),
    \qquad
    \proj_\cX^R (\hat x) = \argmin_{x \in \cX} \norm[R]{x - \hat x}.
\end{equation*}
This step needs only products with~$P$, $A$, and~$A^T$ and a projection.
If~$\cX$ is a box and~$R$ is diagonal, the projection~$\proj_\cX^R = \proj_\cX$ is a componentwise clipping.
The iteration~\eqref{algo:prox-admm} is then
\begin{equation}
    \label{algo:lsf-osqp-step}
    \begin{array}{rcl}
        \xkk &=& \proj_\cX \Big(
            \xk-R^{-1}\big[P\xk+q+A^T\big(\yk+\rho(A\xk-\zk)\big)\big]
        \Big)\\
        \zkk &=& \proj_\cC \big(
            A\xkk+\rho^{-1}\yk
        \big) \\
        \ykk &=& \yk+\rho\big(
            A\xkk-\zkk
        \big).
    \end{array}
\end{equation}
We choose the diagonal preconditioner
\begin{equation*}
    \lambda_A \ge \lambdamax (A^T A), \quad \lambda_A > 0,
    \qquad
    R = \diag \left( {\sum}_{j} |P_{ij}| \right) + \rho \lambda_A \identity,
\end{equation*}
which satisfies the condition~$W \succeq 0$.
Indeed, the matrix~$\diag \big( {\sum}_{j} |P_{ij}| \big) - P$ is diagonally dominant and hence positive semidefinite, so~$W = \diag \big( {\sum}_{j} |P_{ij}| \big) - P + \rho \left( \lambda_A \identity - A^T A \right) \succeq 0$.
If the inequality~$\lambda_A \ge \lambdamax (A^T A)$ is strict,
then~$\rho (\lambda_A \identity - A^T A) \succ 0$ and hence~$W \succ 0$.
Unlike a multiple of the identity such as~$\lambdamax (P) \identity$, the diagonal term keeps the scale of each row of~$P$, since~$\diag \big( {\sum}_{j} |P_{ij}| \big)$ is the smallest diagonal matrix whose difference with~$P$ is diagonally dominant with a nonnegative diagonal.
Forming~$R$ costs one pass over the nonzeros of~$P$ and one estimate of~$\lambda_A$.
The iteration~\eqref{algo:prox-admm} is also an instance of the function-linearized proximal ADMM (FLiP-ADMM) of Ryu and Yin~\cite[Chapter~8]{ryuLargeScaleConvexOptimization2022}, which allows a part of the objective to be linearized.
Linearizing the whole objective gives the same step with a different proximal matrix, and leads to a second step-size condition (\Cref{cor:lsf-averagedness}).
\Cref{app:flip} gives both derivations.

\subsection{Relaxation and Halpern iteration}
\label{subsec:fixed-point-formulation}
In terms of the iterate~$s^k = (x^k, v^k)$ with~$v^k = z^k + \rho^{-1} y^k$, the step~\eqref{algo:lsf-osqp-step} reads
\begin{equation}
    \label{eq:v-update}
    \opF s^k = \big(\xtkk,\, A\xtkk+\rho^{-1}\yk\big),
\end{equation}
where~$\xtkk$ is the output of the~$x$-update of~\eqref{algo:lsf-osqp-step}.
The second component equals~$\zkk+\rho^{-1}\ykk$ by the~$y$-update, and the~$z$-update projects it onto~$\cC$.
Hence the pair~$(z^k,y^k)$ is recovered from~$v^k$ by the Moreau decomposition
\begin{equation}\label{eq:osqp-recovery}
    z^k = \proj_\cC (v^k),
    \qquad y^k = \rho (v^k - z^k).
\end{equation}
With a \emph{relaxation parameter}~$\alpha\in(0,2]$, as in relaxed ADMM~\cite{ecksteinDouglasRachfordSplitting1992}, one iteration of \emph{linear-system-free ADMM} is
\begin{equation}
    \label{eq:relaxed-step}
    s^{k+1}=\opT s^k,
    \qquad
    \opT = (1-\alpha)\identity + \alpha \opF.
\end{equation}
The case~$\alpha=1$ is the step~$\opF$ itself, and the case~$\alpha=2$ is the \emph{reflected step}~$\opT=2\opF-\identity$.
We call~$s^{k+1}=\opT s^k$ the \emph{plain iteration}.
\Cref{algo:lsf-osqp} states it in the variables~$(x,z,y)$.
Each evaluation of~$\opT$ costs one product with each of~$P$, $A$, and~$A^T$ and two projections, since~$A\xkk$ and~$P\xkk$ follow from the products with~$\xtkk$ and~$\xk$.

\begin{algorithm}
    \begin{algorithmic}[1]
        \Statex \textbf{Input:} Initial iterate~$s^0=(x^0,v^0)$, parameters~$\rho>0$ and~$\alpha\in(0,2]$, and preconditioner~$R$
        \State $z^0\gets\proj_\cC(v^0)$; $y^0\gets\rho(v^0-z^0)$
        \Repeat
            \State $\xtkk \gets \proj_\cX \big( \xk-R^{-1}\big[P\xk+q+A^T\big(\yk+\rho(A\xk-\zk)\big)\big] \big)$
            \State $\xkk \gets \alpha\xtkk+(1-\alpha)\xk$
            \State $\zkk \gets \proj_\cC \big( A\xkk-(1-\alpha)(A\xk-\zk)+\rho^{-1}\yk \big)$
            \State $\ykk \gets \yk+\rho\big( A\xkk-(1-\alpha)(A\xk-\zk)-\zkk \big)$
        \Until{the termination criteria of~\Cref{subsec:termination-criteria} hold}
    \end{algorithmic}
    \caption{Linear-system-free ADMM}
    \label{algo:lsf-osqp}
\end{algorithm}

The plain iteration can be accelerated by Halpern iteration~\cite{halpernFixedPointsNonexpanding1967}, which averages each step with a fixed anchor point~$s^0$:
\begin{equation*}\label{eq:halpern}
    s^{k+1} = \frac{k+1}{k+2} \opT s^k + \frac{1}{k+2} s^0,
    \qquad
    k=0,1,\ldots.
    \tag{Halpern}
\end{equation*}
The weight of the anchor decreases as~$1/(k+2)$.
\Cref{algo:restarted-halpern} restarts this iteration from a new anchor whenever a restart condition holds.
Its iterates~$s^{n,k}$ carry two indices: $n$ counts the restarts, or \emph{outer iterations}, and~$k$ counts the \emph{inner iterations} since the last restart, so that~$s^{n,0}$ is the current anchor.
Below we drop the index~$n$ when it is not needed, and write~$s^k$ for the current iterate and~$s^0$ for its anchor.
A restart condition compares the fixed-point residual of the step,
\begin{equation}
    \label{eq:restart-residual}
    \omega^{k}
    = \sqrt{
        \norm[W]{\xtkk - x^{k}}^2
        + \rho \norm{A\xtkk - z^{k}}^2
    },
\end{equation}
with its value~$\omega^0$ at the anchor.
We restart if any of the conditions
\begin{equation}
    \label{eq:restart-conditions}
    \begin{aligned}
    &\omega^{k}\le\beta\,\omega^{0}
        &&\text{(sufficient decay)},\\
    &\omega^{k}\le\beta_\mathrm{nec}\,\omega^{0}
        \quad\text{and}\quad
        \omega^{k}>\omega^\mathrm{prev}
        &&\text{(necessary decay, no local progress)},\\
    &k+1\ge\beta_\mathrm{long}\,t
        &&\text{(long inner loop)}
    \end{aligned}
\end{equation}
holds, where~$\omega^\mathrm{prev}$ is the residual at the previous test, $k+1$ and~$t$ count the evaluations of~$\opT$ in the current outer iteration and in the whole run, and~$\beta<\beta_\mathrm{nec}$ and~$\beta_\mathrm{long}$ are constants in~$(0,1)$.
The sufficient-decay condition is the adaptive rule under which restarted Halpern iteration converges linearly (\Cref{thm:restart}).
The other two follow the restart strategies of PDLP~\cite{applegatePDLPPracticalFirstOrder2025}, HPR-LP~\cite{chenHPRLPImplementationHPR2025}, and HPR-QP~\cite{chenHPRQPDualHalpern2025}.
The necessary-decay condition restarts a run whose residual has decayed moderately but has grown since the previous test, and the long-inner-loop condition keeps a single outer iteration from taking up most of the run.
\Cref{thm:restart} covers the sufficient-decay condition with a fixed step size.
It does not cover the other two conditions or the step-size updates at restarts, which we use because they work well in practice, as these solvers also do.
\Cref{subsec:restart-conditions} gives the constants and the schedule of the tests in our implementation.

\begin{algorithm}
    \begin{algorithmic}[1]
        \Statex \textbf{Input:} Initial iterate~$s^{0,0}$
        \State $n \gets 0$; $k \gets 0$
        \Repeat
            \If{a restart condition~\eqref{eq:restart-conditions} holds}
                \State $s^{n+1,0} \gets \opF s^{n,k}$; $n \gets n+1$; $k \gets 0$
            \Else
                \State $s^{n,k+1} \gets \frac{k+1}{k+2} \opT s^{n,k} + \frac{1}{k+2} s^{n,0}$; $k \gets k+1$
            \EndIf
        \Until{the termination criteria of~\Cref{subsec:termination-criteria} hold}
    \end{algorithmic}
    \caption{Restarted Halpern iteration~\cite{luRestartedHalpernPDHG2024}}
    \label{algo:restarted-halpern}
\end{algorithm}

\subsection{Primal--dual form and relation to other methods}
\label{subsec:primal-dual}
\label{subsec:relation}
Eliminating~$z$ from the step~\eqref{algo:lsf-osqp-step} with its~$y$-update, which gives~$\rho(A\xkk-\zkk)=\ykk-\yk$, yields its primal--dual form:
\begin{equation}
    \label{eq:lsf-primal-dual}
    \begin{array}{rcl}
        \ykk &=& \rho\,(\identity-\proj_\cC)\big(
            A\xkk+\rho^{-1}\yk
        \big) \\
        \xkkk &=& \proj_\cX \Big(
            \xkk-R^{-1}\big[P\xkk+q+A^T(2\ykk-\yk)\big]
        \Big).
    \end{array}
\end{equation}
This iteration maps the pair~$(\xkk,\yk)$ to~$(\xkkk,\ykk)$, so it updates~$y$ before~$x$.
By the Moreau decomposition, the first line is the proximal operator of~$\rho S_\cC$ evaluated at~$\yk+\rho A\xkk$.
It is the preconditioned Condat--V\~u method~\cite[Algorithm~3.2]{condatPrimalDualSplitting2013} with dual step size~$\rho$ and the preconditioner~$R$ in place of~$\tau^{-1}\identity$.
If~$P=0$, then~$R=\rho\lambda_A\identity$, and~\eqref{eq:lsf-primal-dual} is a PDHG step with primal step size~$1/(\rho\lambda_A)$ and dual step size~$\rho$.

The two forms of the step connect linear-system-free ADMM to most first-order solvers for LPs and QPs.
The~$x$-update of the proximal ADMM~\eqref{algo:prox-admm} minimizes over~$\cX$ a quadratic function with Hessian~$P+\rho A^TA+W$, and the choice of the proximal matrix~$W$ distinguishes three methods.
OSQP~\cite{stellatoOSQPOperatorSplitting2020} takes~$W=\sigma\identity$ with a parameter~$\sigma>0$ and~$\cX=\reals^d$, and solves the resulting linear system by a factorization or conjugate gradient.
PDHG takes~$W=\tau^{-1}\identity-\rho A^TA$ with primal step size~$\tau>0$.
The~$x$-update then minimizes~$(1/2)x^TPx+q^Tx+(2\yk-y^{k-1})^TAx+(1/(2\tau))\norm{x-\xk}^2$ over~$\cX$.
With the dual step in~\eqref{eq:lsf-primal-dual}, this is PDHG~\cite{chambolleFirstOrderPrimalDualAlgorithm2011} in the form that extrapolates the dual variable~\cite{esserGeneralFrameworkClass2010}, \cite[\S8]{oconnorEquivalencePrimaldualHybrid2020}.
Linear-system-free ADMM takes~$W=R-P-\rho A^TA$, so~\Cref{algo:lsf-osqp} is the OSQP iteration~\cite[Algorithm~1]{stellatoOSQPOperatorSplitting2020} with the linear system solve replaced by a projected gradient step.

\Cref{tab:solver-comparison} lists these methods with the other first-order solvers for LPs and QPs, and the rest of this subsection compares them with ours.

\begin{table}
    \centering
    \small
    \caption{
        First-order solvers for LPs and QPs compared with \osqplf{}, which runs linear-system-free ADMM and solves both.
        OSQP, PDHG, and linear-system-free ADMM are the proximal ADMM~\eqref{algo:prox-admm} with the proximal matrices~$W=\sigma\identity$, $W=\tau^{-1}\identity-\rho A^TA$, and~$W=R-P-\rho A^TA$, respectively.
        All solvers except OSQP also restart.}
    \label{tab:solver-comparison}
    \begin{adjustbox}{max width=\linewidth}
    \begin{tabular}{@{}lllll@{}}
        \toprule
        Solver & Class & Base iteration & Treatment of~$P$ & Acceleration \\
        \midrule
        OSQP & QP & proximal ADMM~\eqref{algo:prox-admm}, $W=\sigma\identity$ & linear system with~$P+\rho A^TA+\sigma\identity$ & none \\
        \midrule
        PDLP & LP & PDHG & -- & averaging \\
        cuPDLPx & LP & PDHG & -- & reflected Halpern \\
        HPR-LP & LP & proximal ADMM on the dual LP & -- & reflected Halpern \\
        \midrule
        PDHCG & QP & PDHG & inner solve with~$P+\tau^{-1}\identity$ & averaging \\
        PDHCG-II & QP & PDHG & inner solve with~$P+\tau^{-1}\identity$ & reflected Halpern \\
        PDQP & QP & PDHG & gradient at an averaged point & momentum \\
        cuPDQPx & QP & PDHG with an auxiliary variable & auxiliary variable & reflected Halpern \\
        HPR-QP & QP & proximal ADMM on a Wolfe dual & proximal term~$P(\lambda_P\identity-P)$ & reflected Halpern \\
        \midrule
        \osqplf{} & LP, QP & proximal ADMM~\eqref{algo:prox-admm}, $W=R-P-\rho A^TA$ & projected gradient step & reflected Halpern or none \\
        \bottomrule
    \end{tabular}
    \end{adjustbox}
\end{table}

\paragraph{PDHG-based solvers on LPs}
The step-size condition~$\tau\rho\lambdamax(A^TA)\le1$ of PDHG is~$W\succeq0$ for its proximal matrix~$W=\tau^{-1}\identity-\rho A^TA$, so the PDHG step is firmly nonexpansive (\Cref{cor:lsf-averagedness}), and relaxed PDHG converges for every~$\alpha\in(0,2)$~\cite[\S4.3]{oconnorEquivalencePrimaldualHybrid2020}.
Since~$\lambda_A\ge\lambdamax(A^TA)$, the product~$1/\lambda_A$ of the two step sizes satisfies the step-size condition for every~$\rho$, and~$\rho$ sets only their ratio.
PDLP~\cite{applegatePDLPPracticalFirstOrder2025} sets this ratio by its primal weight.
It restarts from an average of the PDHG iterates and adapts the step size and the primal weight during the run.
cuPDLPx~\cite{luCuPDLPxFurtherEnhanced2025} runs the reflected restarted Halpern PDHG of Lu and Yang~\cite{luRestartedHalpernPDHG2024}.
Relaxation, Halpern iteration, and restarts carry over from one form of the step to the other (\Cref{thm:lsf-osqp-condat-vu}), so with~$\alpha=2$ restarted linear-system-free ADMM is this method with the roles of the primal and dual variables exchanged.
Since cuPDLPx extrapolates the primal variable, the iterates of the two solvers differ.

\paragraph{PDHG-based solvers on QPs}
If~$P\neq0$, the primal step of PDHG is no longer explicit.
PDHCG~\cite{huangRestartedPrimalDualHybrid2025} and PDHCG-II~\cite{liPDHCGIIEnhancedVersion2026} solve the primal step inexactly by inner iterations, whereas linear-system-free ADMM has no inner iterations and no inner tolerance.
In return, its step-size condition involves~$P$, through~$R\succeq P+\rho A^TA$ or~$R\succ P/2+\rho A^TA$ (\Cref{cor:lsf-averagedness}), whereas the step-size condition of PDHG involves only~$A$.
PDQP~\cite{luPracticalOptimalFirstOrder2025} restarts the accelerated PDHG of Chen et al.~\cite{chenOptimalPrimalDualMethods2014}, which adds momentum to the primal step by evaluating the gradient at a weighted combination of the current iterate and a running average.
It thus attains the optimal dependence on the condition number of~$P$ on strongly convex QPs, which restarted linear-system-free ADMM does not, as~\Cref{subsec:sharpness-free} shows.
cuPDQPx~\cite{luCuPDQPxGPUAcceleratedRestarted2026} writes the quadratic term as a maximum over an auxiliary variable and applies PDHG to the resulting saddle-point problem.
Its step is firmly nonexpansive, so it also uses reflected restarted Halpern iteration, at the cost of the auxiliary variable and its step size.

\paragraph{HPR-based solvers}
HPR-LP~\cite{chenHPRLPImplementationHPR2025} and HPR-QP~\cite{chenHPRQPDualHalpern2025} run the Halpern Peaceman--Rachford (HPR) method with semi-proximal terms, which is the accelerated proximal ADMM of Sun et al.~\cite{sunAcceleratingPreconditionedADMM2025} with relaxation parameter~$2$.
The step~$\opF$ is a proximal point step (\Cref{thm:lsf-osqp-condat-vu}), so restarted linear-system-free ADMM with~$\alpha=2$ is an HPR method in the same sense.
The three solvers differ in the problem to which they apply the HPR method and in the block that carries the semi-proximal term, which makes the update of that block explicit.
Linear-system-free ADMM works on the primal problem~\eqref{prob:qp} and puts the term~$W=R-P-\rho A^TA$ on the~$x$-block.
HPR-LP works on the dual LP, with the primal variable as the multiplier, and puts the term~$\lambda_A\identity-AA^T$ on the block of the constraint multipliers.
HPR-QP works on the restricted Wolfe dual of the QP, in which an auxiliary variable in the range of~$P$ carries the quadratic term; it puts the term~$P(\lambda_P\identity-P)$ with~$\lambda_P\ge\lambdamax(P)$ on this block and updates it twice per iteration, before and after the constraint multipliers, in a symmetric Gauss--Seidel sweep.
With~$P=0$, HPR-QP reduces to HPR-LP~\cite[Remark~2.1]{chenHPRQPDualHalpern2025}.
All three solvers test the restart conditions~\eqref{eq:restart-conditions} on the fixed-point residual and update the step size at restarts, as~\Cref{subsec:restart-conditions} describes for linear-system-free ADMM.

\section{Convergence analysis}
\label{sec:sharpness}
We first show that the step of linear-system-free ADMM is nonexpansive in a suitable metric.
We then show that its fixed-point residual is sharp, so that restarted Halpern iteration (\Cref{algo:restarted-halpern}) converges linearly~\cite{luRestartedHalpernPDHG2024}.
On equality-constrained LPs, the resulting complexity matches the lower bound of Applegate et al.~\cite{applegateFasterFirstorderPrimaldual2023} up to constants and a logarithmic term.
Finally, when~\eqref{prob:qp} is primal or dual infeasible, a linear transform of the fixed-point residual of the step converges to a limit that certifies this infeasibility.

\subsection{Nonexpansiveness}
\label{subsec:nonexpansiveness}
We write the primal--dual form~\eqref{eq:lsf-primal-dual} with monotone operators.
Over~$u=(x,y)\in\reals^d\times\reals^m$, define the operator~$\opA$ and the matrix~$\opM$ as
\begin{equation*}
    \opA(x,y)=\begin{bmatrix}
        Px+q+\cN_\cX(x)+A^Ty \\
        \partial S_\cC(y)-Ax
    \end{bmatrix},
    \qquad
    \opM=\begin{bmatrix}
        R-P & A^T \\
        A & \rho^{-1}\identity
    \end{bmatrix},
\end{equation*}
where~$\partial S_\cC$ is the subdifferential of the support function of~$\cC$, so that~$z\in\partial S_\cC(y)$ if and only if~$z\in\cC$ and~$y\in\cN_\cC(z)$.
The operator~$\opA$ is monotone, and its zeros are the pairs~$(x,y)$ of the KKT points of~\eqref{prob:qp}.
We also define
\begin{equation*}
    \opTCV=(\opM+\opA)^{-1}\opM,
    \qquad
    \Phi(x,y)=(x,\,Ax+\rho^{-1}y),
    \qquad
    \opW=\diag(W,\rho\identity),
\end{equation*}
where~$(\opM+\opA)^{-1}$ is the set-valued inverse of~$\opM+\opA$.
The operator~$\opTCV$ is a proximal point step on~$\opA$ in the metric~$\opM$, and its subscript refers to the Condat--V\~u method, of which~\eqref{eq:lsf-primal-dual} is an instance (\Cref{app:flip}).
The matrix~$W$ is the proximal matrix of~\Cref{subsec:derivation}.
The next theorem shows that~$\opTCV$ is the primal--dual form~\eqref{eq:lsf-primal-dual}, and that the step~$\opF$ is the same operator in other coordinates.
The linear map~$\Phi$ takes one to the other and is an isometry between their metrics.
Hence every statement about one operator, such as averagedness or sharpness, transfers verbatim to the other, and known results on proximal ADMM and on the Condat--V\~u method apply to linear-system-free ADMM.

\begin{theorem}[Primal--dual representation of linear-system-free ADMM]
    \label{thm:lsf-osqp-condat-vu}
    Suppose that~$W\succeq0$.
    Then the following hold.
    \begin{enumerate}
        \item[($a$)]
            The pairs~$u^k=(\xkk,\yk)$ of the primal--dual form~\eqref{eq:lsf-primal-dual} satisfy~$(\opM+\opA)u^{k+1} \ni \opM u^k$, \ie, $u^{k+1}=\opTCV u^k$.
        \item[($b$)]
            The matrix~$\opM$ is positive (resp.\ semi-)definite if and only if~$W\succ0$ (resp.\ $W\succeq0$).
        \item[($c$)]
            The map~$\Phi$ satisfies~$\Phi^T\opW\Phi=\opM$, so~$\norm[\opW]{\Phi u-\Phi\hat u}=\norm[\opM]{u-\hat u}$ for all~$u$ and~$\hat u$.
            Moreover, $\Phi u^k=(\xkk,v^{k+1})$ with~$v^{k+1}=\zkk+\rho^{-1}\ykk$, so~$\opF = \Phi \circ \opTCV \circ \Phi^{-1}$.
    \end{enumerate}
\end{theorem}

\begin{proof}
    For~($a$), we check the first and second components of the inclusion~$(\opM+\opA)u^{k+1}\ni\opM u^k$ in turn.
    The~$y$-update gives~$\rho \zkk = \rho A \xkk + \yk - \ykk$.
    Substituting~$\rho\zkk$ into the optimality condition of the~$x$-update of~\eqref{algo:prox-admm} at iteration~$k+1$ gives
    \begin{equation*}
        0\in P\xkkk+q+\cN_\cX(\xkkk)+A^T(2\ykk-\yk)+(W+\rho A^TA)(\xkkk-\xkk),
    \end{equation*}
    which is the first component, since~$W+\rho A^TA=R-P$.
    The optimality condition for the dual proximal step in~\eqref{eq:lsf-primal-dual} gives~$A\xkk+\rho^{-1}(\yk-\ykk)\in\partial S_\cC(\ykk)$, which is the second component.
    The second component uniquely determines~$\ykk$ by the dual proximal step, and the first then uniquely determines~$\xkkk$ by projection onto~$\cX$ in the~$R$-norm.
    Direct block multiplication gives~$\Phi^T\opW\Phi=\opM$, which yields the seminorm identity in~($c$).
    Since~$\Phi$ is invertible and~$\opW=\diag(W,\rho\identity)$, the same matrix identity gives~($b$).
    By~\eqref{eq:v-update} with~$\xtkk=\xkk$ for~$\alpha=1$, $v^{k+1}=A\xkk+\rho^{-1}\yk$ is the second component of~$\Phi u^k$.
    Hence~$\Phi\opTCV u^{k-1}=\Phi u^k=s^{k+1}=\opF s^k=\opF\Phi u^{k-1}$, and
    since~$\Phi$ is invertible, $\opF=\Phi\circ\opTCV\circ\Phi^{-1}$.
\end{proof}

\begin{corollary}[Averagedness of linear-system-free ADMM]
    \label{cor:lsf-averagedness}
    Let~$\opT=(1-\alpha)\identity+\alpha\opF$, with~$\opM$ and~$\opW$ as in~\Cref{thm:lsf-osqp-condat-vu}.
    \begin{enumerate}
        \item[($a$)]
            Suppose that~$W\succeq0$.
            Then~$\opTCV$ is firmly nonexpansive in the~$\opM$-seminorm, and~$\opF$ is firmly nonexpansive in the~$\opW$-seminorm.
            Hence, in the latter seminorm, the relaxed operator~$\opT$ is~$\alpha/2$-averaged for every~$\alpha\in(0,2)$, and the reflected operator~$\opT = 2 \opF - \identity$ is nonexpansive at~$\alpha=2$.
        \item[($b$)]
            Suppose that~$W_0=R-\rho A^TA\succ P/2$, and set~$L=\lambda_{\max}(W_0^{-1/2}PW_0^{-1/2})$ and~$\theta=2/(4-L)$.
            Then~$\opF$ is~$\theta$-averaged in the norm of~$\opW_0=\diag(W_0,\rho\identity)$.
            In this norm, $\opT$ is~$\alpha\theta$-averaged for~$0<\alpha<1/\theta$ and nonexpansive for~$\alpha=1/\theta$.
    \end{enumerate}
\end{corollary}

\begin{proof}
    For~($a$), \Cref{thm:lsf-osqp-condat-vu}($a$) gives~$\opM(u-u^+)\in\opA u^+$ for~$u^+=\opTCV u$.
    Applying monotonicity of~$\opA$ to this inclusion at~$u$ and~$\hat u$ gives~$\|u^+-\hat u^+\|_{\opM}^2\le(u^+-\hat u^+)^T\opM(u-\hat u)$, where~$\hat u^+=\opTCV\hat u$.
    \Cref{thm:lsf-osqp-condat-vu}($c$) transfers this bound to~$\opF$, and the relaxation bounds follow from the definition of averagedness.
    \Cref{app:flip} proves~($b$), by reading the same inclusion as a forward--backward step.
\end{proof}

The two parts give different step-size guarantees for the same update.
Part~($a$) is the proximal ADMM condition~$R\succeq P+\rho A^TA$~\cite{sunAcceleratingPreconditionedADMM2025}, and it guarantees nonexpansiveness for~$0<\alpha\le2$, including the reflected step at~$\alpha=2$.
Part~($b$) uses~$R\succ P/2+\rho A^TA$, the strict form of Condat's step-size condition~\cite[Theorem~3.1]{condatPrimalDualSplitting2013}.
For a scalar primal step~$R=\tau^{-1}\identity$, it can permit a larger~$\tau$, while the guaranteed relaxation range ends at~$\alpha=2-L/2$.

In the rest of this section, $\opT$ is the operator~\eqref{eq:relaxed-step} of linear-system-free ADMM with~$\alpha\in(0,2]$.
We use part~($a$) of~\Cref{cor:lsf-averagedness}, and suppose that~$R \succ 0$ satisfies the strict majorization condition
\begin{equation}
    \label{eq:free-strict-majorization}
    R \succ P + \rho A^T A,
\end{equation}
so that~$W\succ0$, the metric~$\opW$ is positive definite, and~$\norm[\opW]{\cdot}$ is a norm.
This holds for the preconditioner of~\Cref{subsec:derivation} when~$\lambda_A>\lambdamax(A^TA)$.
By~\Cref{cor:lsf-averagedness}($a$), the operator~$\opT$ is then nonexpansive in~$\norm[\opW]{\cdot}$.
More precisely, the operator~$\opT$ is~$\alpha/2$-averaged for every~$\alpha\in(0,2)$, so the plain iteration converges to a fixed point of~$\opT$ whenever one exists~\cite[\S5.2]{bauschkeConvexAnalysisMonotone2017}.
For the reflected step~$\alpha=2$, the operator~$\opT=2\opF-\identity$ is nonexpansive but in general not averaged, so the plain iteration need not converge.
Halpern iteration needs only nonexpansiveness.
For every~$\alpha\in(0,2]$, if~$\opT$ has a fixed point, the Halpern iterates converge to a fixed point of~$\opT$~\cite{halpernFixedPointsNonexpanding1967}, \cite[Theorem~30.1]{bauschkeConvexAnalysisMonotone2017}, and the fixed-point residual~$\opT s^k-s^k$ vanishes at the rate~$O(1/k)$~\cite{liederConvergenceRateHalperniteration2021}.
\Cref{subsec:sharpness-free} shows that restarts turn this rate into a linear one.

\subsection{Fixed points and KKT residuals}
\label{subsec:kkt-residuals}
We call~$(\opT-\identity)s$ the \emph{fixed-point residual} at~$s$ and~$(\opF-\identity)s$ the fixed-point residual of~$\opF$.
Since~$\opT-\identity=\alpha(\opF-\identity)$ with~$\alpha>0$, the operators~$\opT$ and~$\opF$ have the same fixed-point set~$\cS^\star=\{ s \mid s = \opT s \}$.
We show that the fixed-point residual of~$\opF$ is an affine expression in terms of the \emph{primal and dual residuals}
\begin{equation*}
    r_p(x,z) = A x - z,
    \qquad
    r_d(x,y,w) = P x + q + A^T y + w.
\end{equation*}
At the iterate, we write~$r_p^k=r_p(x^k,z^k)$ and~$r_d^k=r_d(x^k,y^k,w^{k+1})$,
where the vector
\begin{equation}
    \label{eq:free-normal-vector}
    w^{k+1}
    = R(x^k-\xtkk)-\nabla_x\lagrangian_\rho(x^k,z^k,y^k)
    = R(x^k-\xtkk)-\big(Px^k+q+A^Ty^k+\rho A^T(Ax^k-z^k)\big)
\end{equation}
takes the place of the multiplier~$w$ of the constraint~$x\in\cX$ in~\eqref{eq:qp-kkt}.
By the optimality condition of the projection onto~$\cX$ in the~$x$-update,
this vector lies in~$\cN_\cX(\xtkk)$.
It is zero if~$\cX=\reals^d$.

\begin{lemma}[Fixed-point and KKT residuals]
    \label{lem:free-residual-map}
    The fixed-point residual of~$\opF$ and the stacked residual~$r^k=(r_d^k,\,r_p^k)$ satisfy
    \begin{equation}
        \label{eq:free-residual-matrix}
        (\opF-\identity)s^k
        = -U r^k,
        \quad
        U
        = \begin{bmatrix}
            R^{-1} & \rho R^{-1}A^T \\
            AR^{-1} & \rho AR^{-1}A^T-\identity
        \end{bmatrix}
        = \begin{bmatrix}\identity\\A\end{bmatrix}R^{-1}\begin{bmatrix}\identity&\rho A^T\end{bmatrix}
            -\begin{bmatrix}0&0\\0&\identity\end{bmatrix}.
    \end{equation}
    The matrix~$U$ is invertible, with
    \begin{equation}
        \label{eq:inverse-kkt-residual-map}
        r^k
        = -U^{-1}(\opF-\identity)s^k,
        \qquad
        U^{-1}
        = \begin{bmatrix}
            R-\rho A^TA & \rho A^T \\
            A & -\identity
        \end{bmatrix}.
    \end{equation}
    The iterate~$s^k$ is a fixed point of~$\opF$ if and only if~$r^k=0$, which in turn holds if and only if~$(x^k,z^k,y^k,w^{k+1})$ satisfies the KKT conditions~\eqref{eq:qp-kkt}.
\end{lemma}

\begin{proof}
    Set~$\delta s^k=(\opF-\identity)s^k=(\delta x^k,\delta v^k)$.
    The definition~\eqref{eq:free-normal-vector} of~$w^{k+1}$ and~\eqref{eq:v-update} give
    \begin{align*}
        R\delta x^k
        &= R(\xtkk-x^k)
        = -\big(r_d^k+\rho A^Tr_p^k\big) \\
        \delta v^k
        &= A\xtkk+\rho^{-1}y^k-v^k
        = A\xtkk-z^k
        = A\delta x^k+r_p^k.
    \end{align*}
    Substituting~$\delta x^k=-R^{-1}(r_d^k+\rho A^Tr_p^k)$ into the second identity gives~\eqref{eq:free-residual-matrix}.
    Conversely, the second identity gives~$r_p^k=\delta v^k-A\delta x^k$, and substituting it into the first gives~$r_d^k=-(R-\rho A^TA)\delta x^k-\rho A^T\delta v^k$.
    This is~\eqref{eq:inverse-kkt-residual-map}, and it shows that~$U$ is invertible.

    If~$\delta s^k=0$, $x^k=\xtkk\in\cX$, $w^{k+1}\in\cN_\cX(x^k)$, and~$r_p^k=r_d^k=0$.
    Together with~$z^k\in\cC$ and~$y^k\in\cN_\cC(z^k)$ from~\eqref{eq:osqp-recovery}, these show that~$(x^k,z^k,y^k,w^{k+1})$ satisfies the KKT conditions~\eqref{eq:qp-kkt}.
    Conversely, these conditions give~$r_p^k=r_d^k=0$, so~\eqref{eq:free-residual-matrix} gives~$\delta s^k=0$.
\end{proof}

By~\Cref{lem:free-residual-map}, the Moreau decomposition~\eqref{eq:osqp-recovery}, and the optimality condition for the projection onto~$\cX$, the change of variables~$\Phi$ of~\Cref{thm:lsf-osqp-condat-vu}($c$) maps the pairs~$(x,y)$ of the KKT points of~\eqref{prob:qp} onto~$\cS^\star$.
Hence~$\cS^\star \neq \emptyset$ if and only if~\eqref{prob:qp} has a KKT point.
The point~$(\xtkk,z^k,y^k,w^{k+1})$ satisfies every inclusion of the KKT conditions~\eqref{eq:qp-kkt}.

\subsection{Sharpness and linear convergence}
\label{subsec:sharpness-free}
Restarted Halpern iteration converges linearly when the fixed-point residual is sharp~\cite{luRestartedHalpernPDHG2024}.
We show that the fixed-point residual of linear-system-free ADMM is sharp, and we compare the resulting complexity with lower bounds for first-order methods.
Let~$\dist_{\opW}(s,\cS^\star)=\min_{\hat s\in\cS^\star}\norm[\opW]{s-\hat s}$.
For an operator~$\opT$ that is nonexpansive in~$\norm[\opW]{\cdot}$ and has a fixed point, the Halpern iterates~\eqref{eq:halpern} satisfy~\cite{liederConvergenceRateHalperniteration2021}
\begin{equation*}
    \norm[\opW]{\opT s^k - s^k} \le \frac{2}{k+1} \dist_{\opW} (s^0, \cS^\star).
\end{equation*}

\begin{assumption}[Sharpness]
    \label{assumption:sharpness}
    For some~$s^\star\in\cS^\star$, let~$\cB_0=\{ s \mid \norm[\opW]{s-s^\star}\le\norm[\opW]{s^{0,0}-s^\star} \}$.
    There exists~$\eta>0$ such that~$\norm[\opW]{\opT s-s}\ge\eta\,\dist_{\opW}(s,\cS^\star)$ for every~$s\in\cB_0$.
\end{assumption}
Since~$\opT$ is nonexpansive, every iterate of~\Cref{algo:restarted-halpern} lies in the compact ball~$\cB_0$~\cite[Lemma~2]{luRestartedHalpernPDHG2024}, so sharpness is needed only there.

\begin{theorem}[Restarted Halpern under sharpness, {\cite[Theorems~1--2; Remark~1]{luRestartedHalpernPDHG2024}}]
    \label{thm:restart}
    Suppose that~\Cref{assumption:sharpness} holds.
    Fix~$\beta\in(0,1)$ and set~$k^\star=\lceil2/(\eta\beta)\rceil$.
    \begin{enumerate}
        \item[$(i)$]
            The \emph{fixed rule}, which restarts after~$k^\star$ inner iterations, gives~$\dist_{\opW}(s^{n,0},\cS^\star)\le\beta^n\dist_{\opW}(s^{0,0},\cS^\star)$ for~$n\ge0$.
        \item[$(ii)$]
            The \emph{adaptive rule} restarts at the first~$k$ with~$\norm[\opW]{\opT s^{n,k}-s^{n,k}}\le\beta\norm[\opW]{\opT s^{n,0}-s^{n,0}}$.
            It gives~$\dist_{\opW}(s^{n,0},\cS^\star)\le\beta^{n+1}k^\star\dist_{\opW}(s^{0,0},\cS^\star)$ for~$n\ge0$, and each of its outer iterations has at most~$k^\star$ inner iterations.
    \end{enumerate}
    Consequently, either rule reaches~$\dist_{\opW}(s^{n,0},\cS^\star)\le\varepsilon\dist_{\opW}(s^{0,0},\cS^\star)$ for~$\varepsilon\in(0,1/e]$ within~$O\!\left(\eta^{-1}\log\big(1/(\eta\varepsilon)\big)\right)$ evaluations of~$\opT$, with a constant that depends only on~$\beta$.
\end{theorem}

\Cref{thm:restart} is adapted from Lu and Yang~\cite{luRestartedHalpernPDHG2024}, who state it for Halpern PDHG on LPs.
Their proof uses only the nonexpansiveness of~$\opT$, the Halpern rate, \Cref{assumption:sharpness}, and two properties of the restart point.
\Cref{algo:restarted-halpern} restarts from~$\opF s^{n,k}$, which has both: since~$\opF$ is nonexpansive in~$\norm[\opW]{\cdot}$ and has the fixed points of~$\opT$, this point is no farther from~$\cS^\star$ than~$s^{n,k}$, and its fixed-point residual is no larger, \ie, $\norm[\opW]{\opT\opF s-\opF s}\le\norm[\opW]{\opT s-s}$.
Hence the theorem also holds for restarted linear-system-free ADMM.
What remains is to verify~\Cref{assumption:sharpness} for this operator, on QPs and not only on LPs.
The next theorem does so.

\begin{theorem}[Sharpness of linear-system-free ADMM]
    \label{thm:sharpness}
    Suppose that~$\cS^\star\neq\emptyset$.
    For every bounded set~$\cB\subset\reals^{d+m}$, there exists~$\eta_{\opW,\cB}>0$ such that
    \begin{equation}
        \label{eq:free-sharpness-Q-exact}
        \norm[\opW]{\opT s-s}
        \ge\eta_{\opW,\cB}\dist_{\opW}(s,\cS^\star),
        \qquad \text{for all } s\in\cB.
    \end{equation}
\end{theorem}

\begin{proof}
    Let~$G=\opT-\identity$, so that~$\cS^\star=G^{-1}(0)$.
    Since the projections~$\proj_\cC$ and~$\proj_\cX^R$ onto boxes are piecewise affine, so is~$G$, and its inverse~$G^{-1}$ is a polyhedral multifunction.
    Robinson's error bound~\cite{robinsonContinuityPropertiesPolyhedral1981}, which holds in any norm, therefore gives~$\varepsilon>0$ and~$\kappa>0$ such that every~$s$ with~$\norm[\opW]{Gs}\le\varepsilon$ satisfies~$\dist_{\opW}(s,\cS^\star)\le\kappa\norm[\opW]{Gs}$.
    This bound requires only a small residual, not that~$s$ lie near~$\cS^\star$.
    Since~$\cB$ is bounded and~$\cS^\star$ is nonempty, the number~$D=\sup_{s\in\cB}\dist_{\opW}(s,\cS^\star)$ is finite.
    If~$D=0$, then~$\cB\subseteq\cS^\star$, and~\eqref{eq:free-sharpness-Q-exact} holds with~$\eta_{\opW,\cB}=1/\kappa$.
    Otherwise, every~$s\in\cB$ with~$\norm[\opW]{Gs}>\varepsilon$ satisfies~$\dist_{\opW}(s,\cS^\star)\le D<(D/\varepsilon)\norm[\opW]{Gs}$, so~\eqref{eq:free-sharpness-Q-exact} holds with~$\eta_{\opW,\cB}=\min\{1/\kappa,\varepsilon/D\}$.
\end{proof}

Applying~\Cref{thm:sharpness} with~$\cB$ equal to the compact ball~$\cB_0$ of~\Cref{assumption:sharpness} verifies~\Cref{assumption:sharpness} with~$\eta=\eta_{\opW,\cB_0}>0$.
Hence~\Cref{thm:restart} applies to restarted linear-system-free ADMM, and its linear rate is governed by that constant.
The proof of~\Cref{thm:sharpness} uses only that~$\cC$ and~$\cX$ are polyhedral, so the same holds for general polyhedral constraints.
By~\Cref{lem:projected-residual}, the KKT residual norm at the projected point computed from the restart point~$s^{n,0}$ is~$O(\beta^n)$.

\paragraph{Explicit rates and lower bounds}
\Cref{thm:sharpness} shows that a sharpness constant exists, but it does not give its value, so the rate of~\Cref{thm:restart} is not explicit.
On two problem classes the operator~$\opT$ is affine, and with the scalar preconditioner~$R=\tau^{-1}\identity$ we compute the best sharpness constant in closed form (\Cref{app:exact-sharpness}).
This gives explicit bounds on the number of evaluations of~$\opT$ that reduce the distance to~$\cS^\star$ by a factor~$\varepsilon\in(0,1/e]$, which we compare with the lower bounds for first-order methods.
\begin{itemize}
    \item \emph{Equality-constrained LPs.}
        With~$\tau^{-1}=c\rho\norm{A}^2$ for a fixed~$c>1$, restarted linear-system-free ADMM needs at most~$O(\kappa_A\log(\kappa_A/\varepsilon))$ evaluations of~$\opT$, where~$\kappa_A=\norm{A}/\sigma_{\min}^+(A)$, and~$\norm{A}$ and~$\sigma_{\min}^+(A)$ are the largest and the smallest positive singular values of~$A$.
        The lower bound for first-order primal--dual methods is~$\Omega(\kappa_A\log(1/\varepsilon))$~\cite{applegateFasterFirstorderPrimaldual2023}.
        Hence the method is optimal on this class, up to constants and a logarithmic term (\Cref{cor:lp-optimality}).
    \item \emph{Unconstrained strongly convex QPs.}
        With~$\tau^{-1}=\Theta(\lambdamax(P))$ and~$\rho\norm{A}^2=O(\lambdamax(P))$, the method needs at most~$O(\kappa_P\log(\kappa_P/\varepsilon))$ evaluations, where~$\kappa_P=\lambdamax(P)/\lambdamin(P)$.
        The lower bound for first-order methods is~$\Omega(\sqrt{\kappa_P}\log(1/\varepsilon))$~\cite[Theorem~2.1.13]{nesterovLecturesConvexOptimization2018}, and momentum methods such as PDQP attain it~\cite{luPracticalOptimalFirstOrder2025}.
        Hence our bound is larger than the lower bound by a factor~$\sqrt{\kappa_P}$.
        This gap is not an artifact of the analysis: no restart rule for Halpern iteration improves the factor~$\kappa_P$ (\Cref{thm:restart-tightness}).
\end{itemize}

\subsection{Infeasibility detection}
\label{subsec:infeasible}
We say that~\eqref{prob:qp} is \emph{primal infeasible} if no~$x\in\cX$ satisfies~$Ax\in\cC$, and \emph{dual infeasible} if no~$(x,y,w)$ with~$y\in\barc\cC$ and~$w\in\barc\cX$ satisfies the constraint of~\eqref{prob:qp-dual}.
In either case~\eqref{prob:qp} has no KKT point, so~$\cS^\star=\emptyset$ by~\Cref{subsec:kkt-residuals}.
The norm of the iterates of~\Cref{algo:restarted-halpern} then tends to infinity.
Without restarts and for a fixed step size, the fixed-point residuals still converge (\Cref{prop:candidate-convergence}).
A linear transform of their limit certifies primal and dual infeasibility, each whenever it occurs.

The certificates use the barrier cones of~\Cref{sec:lsf-osqp} and the \emph{recession cone}~$\rec\cC$ of a box~\cite{bauschkeConvexAnalysisMonotone2017}, which collects the directions along which the box is unbounded,
\begin{equation*}
    \rec\cC=\{ e\in\reals^m \mid e_i\le0 \text{ if } u_i<+\infty, \quad e_i\ge0 \text{ if } l_i>-\infty \}.
\end{equation*}
The cone~$\rec\cX$ is defined in the same way.

\begin{lemma}[Infeasibility certificates]
    \label{lem:infeasibility-certificates}
    The following hold.
    \begin{enumerate}
        \item[($a$)]
            If there exists~$\delta y\in\reals^m$ such that
            \begin{equation}
                \label{eq:primal-certificate}
                S_\cC(\delta y)+S_\cX(-A^T\delta y)<0,
            \end{equation}
            then~\eqref{prob:qp} is primal infeasible.
        \item[($b$)]
            If there exists~$\delta x\in\reals^d$ such that
            \begin{equation}
                \label{eq:dual-certificate}
                \delta x\in\rec\cX, \qquad A\delta x\in\rec\cC, \qquad P\delta x=0, \qquad q^T\delta x<0,
            \end{equation}
            then~\eqref{prob:qp} is dual infeasible.
            If, in addition, \eqref{prob:qp} is primal feasible, then its optimal value is~$-\infty$.
    \end{enumerate}
\end{lemma}

\begin{proof}
    For~($a$), suppose that~\eqref{prob:qp} is feasible, \ie, there exists~$x\in\cX$ with~$Ax\in\cC$.
    Then~$\delta y^T(Ax)\le S_\cC(\delta y)$ and~$(-A^T\delta y)^Tx\le S_\cX(-A^T\delta y)$, and adding the two inequalities gives~$0\le S_\cC(\delta y)+S_\cX(-A^T\delta y)<0$, a contradiction.
    For~($b$), suppose that~$(x,y,w)$ satisfies the constraint of~\eqref{prob:qp-dual} with~$y\in\barc\cC$ and~$w\in\barc\cX$.
    By the sign patterns of the cones, $y^T(A\delta x)\le0$ and~$w^T\delta x\le0$.
    The inner product of the constraint with~$\delta x$ and~$P\delta x=0$ then give~$0=q^T\delta x+y^T(A\delta x)+w^T\delta x\le q^T\delta x<0$, a contradiction.
    If~$x_0$ is primal feasible, then so is~$x_0+t\delta x$ for every~$t\ge0$, and its objective value, which equals that of~$x_0$ plus~$t\,q^T\delta x$, tends to~$-\infty$.
\end{proof}

The inequality~\eqref{eq:primal-certificate} requires both support functions to be finite, so it implies~$\delta y\in\barc\cC$ and~$-A^T\delta y\in\barc\cX$.
For~$\cX=\reals^d$, the cone~$\barc\cX$ is~$\{0\}$, and~($a$) and~($b$) are the certificates that OSQP tests~\cite[Section~5.2]{banjacInfeasibilityDetectionAlternating2019}.
For a general~$\cX$, they are the certificates of Banjac~\cite[Proposition~4.1]{banjacMinimalDisplacementVector2021}.

\begin{definition}
    \label{def:candidates}
    The \emph{certificate candidates} at the iterate~$s^k=(x^k,v^k)$ are
    \begin{equation*}
        (\delta x^k,\delta y^k)
        =\Phi^{-1}(\opF-\identity)s^k
        =(\opTCV-\identity)\Phi^{-1}s^k,
    \end{equation*}
    where the second equality follows from~$\opF=\Phi\circ\opTCV\circ\Phi^{-1}$ (\Cref{thm:lsf-osqp-condat-vu}($c$)).
\end{definition}

By~\Cref{lem:free-residual-map} and the definition of~$\Phi^{-1}$, the certificate candidates are
\begin{equation}
    \label{eq:candidates}
    \delta x^k=\xtkk-x^k=-R^{-1}(r_d^k+\rho A^Tr_p^k),
    \qquad
    \delta y^k=\rho\,(Ax^k-z^k)=\rho\,r_p^k,
\end{equation}
where~$r_p^k$ and~$r_d^k$ are the primal and dual residuals of~\Cref{subsec:kkt-residuals}.
For the plain iteration~$s^{k+1}=\opT s^k$, the certificate candidates equal~$\alpha^{-1}\Phi^{-1}(s^{k+1}-s^k)$, so they are successive differences of the iterates, as in Banjac et al.~\cite{banjacInfeasibilityDetectionAlternating2019}.
The certificate candidates, however, depend only on the current iterate~$s^k$ and not on the update that produced it, so the same definition serves the Halpern iteration~\eqref{eq:halpern}.
Computing the certificate candidates costs nothing extra: the solver forms~$\xtkk$ and~$r_p^k$ when it evaluates~$\opT s^k$.

Consider the sets of values that the primal and dual residuals of~\Cref{subsec:kkt-residuals} can take,
\begin{equation}
    \label{eq:gap-sets}
    \cG_\mathrm{p}=\{ Ax-z \mid x\in\cX,\ z\in\cC \},
    \qquad
    \cG_\mathrm{d}=\{ Px+q+A^Ty+w \mid x\in\reals^d,\ y\in\barc\cC,\ w\in\barc\cX \}.
\end{equation}
Problem~\eqref{prob:qp} is primal feasible if and only if~$0\in\cG_\mathrm{p}$, and dual feasible if and only if~$0\in\cG_\mathrm{d}$.
Both sets are nonempty, convex, and closed (\Cref{lem:range-saddle}), so each has a unique element of minimal norm.
We define the \emph{minimal residuals} and the vector~$(\delta x^\star,\delta y^\star)$ as
\begin{equation}
    \label{eq:closed-form}
    r_p^\star=\argmin_{r\in\cG_\mathrm{p}}\norm{r},
    \qquad
    r_d^\star=\argmin_{r\in\cG_\mathrm{d}}\norm[R^{-1}]{r},
    \qquad
    \delta x^\star=-R^{-1}r_d^\star,
    \qquad
    \delta y^\star=\rho\,r_p^\star.
\end{equation}
By~\eqref{eq:closed-form}, $\rho^{-1}\norm{\delta y^\star}=\norm{r_p^\star}=\dist(A\cX,\cC)$ and~$\norm[R]{\delta x^\star}=\norm[R^{-1}]{r_d^\star}$ measure how far~\eqref{prob:qp} is from primal and dual feasibility.

\begin{theorem}[Minimal displacement of linear-system-free ADMM]
    \label{thm:minimal-displacement}
    The vector~$(\delta x^\star,\delta y^\star)$ of~\eqref{eq:closed-form} has the following properties.
    \begin{enumerate}
        \item[($a$)]
            The vector~$\delta x^\star$ satisfies~$\delta x^\star\in\rec\cX$, $A\delta x^\star\in\rec\cC$, $P\delta x^\star=0$, and
            \begin{equation}
                \label{eq:dual-identity}
                q^T\delta x^\star=-\norm[R]{\delta x^\star}^2.
            \end{equation}
        \item[($b$)]
            The vector~$\delta y^\star$ satisfies~$\delta y^\star\in\barc\cC$, $-A^T\delta y^\star\in\barc\cX$, and
            \begin{equation}
                \label{eq:primal-identity}
                S_\cC(\delta y^\star)+S_\cX(-A^T\delta y^\star)=-\rho^{-1}\norm{\delta y^\star}^2.
            \end{equation}
        \item[($c$)]
            The two blocks satisfy~$(A\delta x^\star)^T\delta y^\star=0$ and~$\norm[\opM]{(\delta x^\star,\delta y^\star)}^2=\norm[R]{\delta x^\star}^2+\rho^{-1}\norm{\delta y^\star}^2$.
        \item[($d$)]
            The range of~$\opTCV-\identity$ is closed and convex, and~$(\delta x^\star,\delta y^\star)$ is its unique element of minimal~$\opM$-norm, the \emph{minimal displacement} of~$\opTCV$.
        \end{enumerate}
        The vector~$\delta y^\star$ is nonzero if and only if~\eqref{prob:qp} is primal infeasible.
        The vector~$\delta x^\star$ is nonzero if and only if~\eqref{prob:qp} is dual infeasible.
        The vector~$(\delta x^\star,\delta y^\star)$ vanishes if and only if~$\cS^\star \neq \emptyset$.
\end{theorem}

The minimal residuals~$r_p^\star$ and~$r_d^\star$ bound the residuals of the two KKT equations at the iterates from below.
By~\Cref{subsec:kkt-residuals}, the point~$(\xtkk,z^k,y^k,w^{k+1})$ satisfies every inclusion of~\eqref{eq:qp-kkt}, and the normal cones of a box lie in its barrier cone.
Hence~$A\xtkk-z^k\in\cG_\mathrm{p}$ and~$P\xtkk+q+A^Ty^k+w^{k+1}\in\cG_\mathrm{d}$, so every iterate~$s^k$ satisfies
\begin{equation*}
    \norm{A\xtkk-z^k}\ge\norm{r_p^\star},
    \quad
    \norm[R^{-1}]{P\xtkk+q+A^Ty^k+w^{k+1}}\ge\norm[R^{-1}]{r_d^\star}.
\end{equation*}

By~\Cref{thm:minimal-displacement}($a$) and~($b$), a nonzero~$\delta x^\star$ satisfies~\eqref{eq:dual-certificate} and a nonzero~$\delta y^\star$ satisfies~\eqref{eq:primal-certificate}.
Hence the converse of~\Cref{lem:infeasibility-certificates} also holds: a certificate~\eqref{eq:primal-certificate} exists whenever~\eqref{prob:qp} is primal infeasible, and a certificate~\eqref{eq:dual-certificate} exists whenever it is dual infeasible.
Banjac et al.~\cite[Section~3.2]{banjacInfeasibilityDetectionAlternating2019} note this converse for a polyhedral or compact~$\cC$ and~$\cX=\reals^d$.
\Cref{thm:minimal-displacement} is the counterpart, for linear-system-free ADMM, of Banjac's decomposition of the minimal displacement of the Douglas--Rachford operator~\cite[Propositions~3.1, 3.2, 4.2]{banjacMinimalDisplacementVector2021}.

\begin{proposition}[Convergence of the certificate candidates]
    \label{prop:candidate-convergence}
    Fix~$\rho$ and let~$s^0$ be arbitrary.
    Then the certificate candidates~\eqref{eq:candidates} converge to the minimal displacement~$(\delta x^\star,\delta y^\star)$ in each of the following cases:
    \begin{enumerate}
        \item[($a$)]
            the plain iteration~$s^{k+1}=\opT s^k$ with~$\alpha\in(0,2)$;
        \item[($b$)]
            the Halpern iteration~\eqref{eq:halpern} with~$\alpha\in(0,2]$.
    \end{enumerate}
\end{proposition}

Hence the certificate candidates converge to a certificate of each infeasibility that occurs.
The proof in~\Cref{app:infeasibility} derives both cases from the bounds of Park and Ryu~\cite{parkAcceleratedInfeasibilityDetection2023}.
Case~($a$) is the counterpart of~\cite[Lemma~5.1]{banjacInfeasibilityDetectionAlternating2019}, which covers relaxed ADMM with~$\alpha\in(0,2)$; case~($b$) also covers the reflected step~$\alpha=2$.
The proposition assumes a fixed step size~$\rho$ and no restarts, while our implementation restarts and updates~$\rho$ at restarts (\Cref{subsec:step-sizes}), so it applies between two consecutive restarts and does not cover the whole run.

\subsection{Termination criteria}
\label{subsec:termination-criteria}
At a point~$(x,y,w)$, the termination criteria use the primal and dual residuals
\begin{align*}
    r_{\mathrm{prim}}
    &=
    \begin{bmatrix}
        A x - \proj_\cC (A x) \\
        x - \proj_\cX (x)
    \end{bmatrix},&
    r_{\mathrm{dual}}
    &=
        P x + q + A^T y + w.
\end{align*}
The duality gap at this point is
\begin{equation*}
    r_{\mathrm{gap}} =
        x^T P x
        + q^T x
        + S_\cC (y) + S_\cX (w).
\end{equation*}
This gap is the difference between the primal objective~$(1/2)x^TPx+q^Tx$ and the objective~$-(1/2)x^TPx-S_\cC(y)-S_\cX(w)$ of the dual problem~\eqref{prob:qp-dual}.
The solver terminates when the residuals and duality gap satisfy
\begin{equation}
    \label{eq:termination}
    \norm[\infty]{r_{\mathrm{prim}}}\le\varepsilon_\mathrm{prim},
    \qquad
    \norm[\infty]{r_{\mathrm{dual}}}\le\varepsilon_\mathrm{dual},
    \qquad
    |r_{\mathrm{gap}}|\le\varepsilon_\mathrm{gap},
\end{equation}
where,
given the absolute and relative tolerances~$\varepsilon_\mathrm{abs}>0$ and~$\varepsilon_\mathrm{rel}>0$,
\begin{align*}
    \varepsilon_\mathrm{prim}
    &=
        \varepsilon_\mathrm{abs}
        + \varepsilon_\mathrm{rel} \max \{
            \norm[\infty]{Ax},
            \norm[\infty]{x}
        \},\\
    \varepsilon_\mathrm{dual}
    &=
        \varepsilon_\mathrm{abs}
        + \varepsilon_\mathrm{rel} \max \{
            \norm[\infty]{P x},
            \norm[\infty]{A^T y + w},
            \norm[\infty]{q}
        \},\\
    \varepsilon_\mathrm{gap}
    &=
        \varepsilon_\mathrm{abs}
        + \varepsilon_\mathrm{rel} \max \{
            |x^T P x|,
            |q^Tx|,
            |S_\cC(y)+S_\cX(w)|
        \}.
\end{align*}
These criteria follow those of OSQP~\cite{stellatoOSQPOperatorSplitting2020} and extend them with the variable bounds, their dual variable~$w$, and a duality gap test of the same form.
\Cref{app:relation} compares the two sets of criteria.
At iteration~$k$, the solver evaluates the optimality criterion~\eqref{eq:termination} at~$x=\xtkk$, $y=y^k$, and~$w=w^{k+1}$, where~\Cref{lem:projected-residual} bounds the KKT residuals by the fixed-point residual.
Since~$\xtkk\in\cX$ and~$z^k\in\cC$, these KKT residuals bound~$r_\mathrm{prim}$ and~$r_\mathrm{dual}$, so a small fixed-point residual guarantees small primal and dual residuals.
Moreover, $y^k\in\cN_\cC(z^k)$ and~$w^{k+1}\in\cN_\cX(\xtkk)$, so the support functions in the duality gap are finite, with~$S_\cC(y^k)=(y^k)^Tz^k$ and~$S_\cX(w^{k+1})=(w^{k+1})^T\xtkk$.

Criteria of the same kind as~\eqref{eq:termination} detect infeasibility.
For vectors~$\delta x\in\reals^d$ and~$\delta y\in\reals^m$, we define the infeasibility residuals as
\begin{align*}
    r_{\mathrm{pinf}}
    &=
    \begin{bmatrix}
        \delta y - \proj_{\barc\cC}(\delta y) \\
        -A^T\delta y - \proj_{\barc\cX}(-A^T\delta y)
    \end{bmatrix},&
    r_{\mathrm{dinf}}
    &=
    \begin{bmatrix}
        P\delta x \\
        A\delta x - \proj_{\rec\cC}(A\delta x) \\
        \delta x - \proj_{\rec\cX}(\delta x)
    \end{bmatrix}.
\end{align*}
The solver reports \emph{primal infeasibility} and terminates if
\begin{equation}
    \label{eq:primal-infeasibility}
    \norm[\infty]{r_{\mathrm{pinf}}}\le\varepsilon_\mathrm{pinf}\norm[\infty]{\delta y},
    \qquad
    S_\cC \big(
        \proj_{\barc\cC}(\delta y)
    \big) + S_\cX \big(
        \proj_{\barc\cX}(-A^T\delta y)
    \big) < 0,
\end{equation}
and it reports \emph{dual infeasibility} and terminates if
\begin{equation}
    \label{eq:dual-infeasibility}
    \norm[\infty]{r_{\mathrm{dinf}}}\le\varepsilon_\mathrm{dinf}\norm[\infty]{\delta x},
    \qquad
    q^T\delta x<0,
\end{equation}
where~$\varepsilon_\mathrm{pinf}>0$ and~$\varepsilon_\mathrm{dinf}>0$ are tolerances.
If~$r_{\mathrm{pinf}}=0$, the projections leave~$\delta y$ and~$-A^T\delta y$ unchanged, and~\eqref{eq:primal-infeasibility} reduces to~\eqref{eq:primal-certificate}; likewise, \eqref{eq:dual-infeasibility} with~$r_{\mathrm{dinf}}=0$ is~\eqref{eq:dual-certificate}.
At every termination check, the solver evaluates these criteria at the certificate candidates~$(\delta x^k,\delta y^k)$ of~\eqref{eq:candidates}, which converge to~$(\delta x^\star,\delta y^\star)$ by~\Cref{prop:candidate-convergence}.
As Banjac et al.~\cite[Section~5.2 and Remark~5.2]{banjacInfeasibilityDetectionAlternating2019} recommend, the criteria test whether a certificate candidate nearly satisfies the conditions of~\Cref{lem:infeasibility-certificates}, not whether successive candidates are close.

\begin{remark}[Boxes]
    \label{rem:boxes}
    By the sign patterns of the recession and barrier cones, the projections onto~$\rec\cC$, $\barc\cC$, $\rec\cX$, and~$\barc\cX$ act entrywise: they set to zero the entries with a forbidden sign and keep the others, so each block of~$r_{\mathrm{pinf}}$ and~$r_{\mathrm{dinf}}$ collects the sign violations over the entries.
    The projections onto the barrier cones zero exactly the entries that would multiply an infinite bound, so, by the formula for~$S_\cC$ and~$S_\cX$ in~\Cref{sec:lsf-osqp}, the second inequality in~\eqref{eq:primal-infeasibility} equals
    \begin{equation*}
        u^T (\delta y)_+
        + l^T (\delta y)_- 
        + u_x^T (- A^T \delta y)_+
        + l_x^T (- A^T \delta y)_-
        < 0,
    \end{equation*}
    where the inner products are taken over the finite entries of~$l$, $u$, $l_x$, and~$u_x$ only.
\end{remark}

\section{GPU batching and algorithm hedging}
\label{sec:gpu-batching}
Linear-system-free ADMM consists only of matrix--vector products, projections onto boxes, and vector operations, so many copies of it can run on a GPU as a single batch.
We first describe the batched iteration, which solves several instances that share the matrices~$P$ and~$A$.
We then use the same batch to run several variants of the algorithm on one instance, which we call \emph{algorithm hedging}.

\subsection{Batched iteration}
\label{subsec:batch-mechanism}
Consider~$B$ instances of~\eqref{prob:qp} that share the matrices~$P$ and~$A$,
\begin{equation*}
    \begin{array}{ll}
        \textnormal{minimize} & (1/2) x_i^T P x_i + q_i^T x_i \\
        \textnormal{subject to} & Ax_i \in \cC_i, \quad x_i \in \cX_i,
    \end{array}
    \qquad i=1,\ldots,B,
\end{equation*}
where instance~$i$ has its own optimization variable~$x_i\in\reals^d$, its own vector~$q_i$, and its own boxes~$\cC_i$ and~$\cX_i$.
Instance~$i$ also runs with its own step size~$\rho_i$, relaxation parameter~$\alpha_i$, and diagonal preconditioner~$R_i\succ P+\rho_iA^TA$.
The instances need not be distinct.
In general, the batch holds~$N_\mathrm{inst}$ distinct instances, and each of them occupies~$E$ columns that run different variants of the algorithm (\Cref{sec:hedging}), so that~$B=N_\mathrm{inst}E$.
This section treats every column as a separate instance.
We store the batch in matrices with one column per instance,
\begin{equation*}
    Q=\begin{bmatrix}q_1&\cdots&q_B\end{bmatrix}\in\reals^{d\times B},
    \qquad
    X=\begin{bmatrix}x_1&\cdots&x_B\end{bmatrix}\in\reals^{d\times B},
    \qquad
    V,Z,Y\in\reals^{m\times B},
\end{equation*}
where~$V$, $Z$, and~$Y$ collect the iterates~$v$, $z$, and~$y$ of the instances in the same way.
We write~$M_i$ for column~$i$ of a matrix~$M$, and~$S=(X,V)$ for the state of the batch.
We collect the parameters of the instances in the vectors
\begin{equation*}
    \brho=(\rho_1,\ldots,\rho_B),
    \qquad
    \brho^{-1}=(1/\rho_1,\ldots,1/\rho_B),
    \qquad
    \balpha=(\alpha_1,\ldots,\alpha_B).
\end{equation*}

Because all instances share~$P$ and~$A$, the products~$PX$, $AX$, and~$A^TY$ apply these matrices to every column at once, and each of them reads~$P$ or~$A$ from memory once for the whole batch.
All other operations act on each column separately.
For a vector~$\bc\in\reals^B$ and a matrix~$M$ with~$B$ columns, we define them through their columns as
\begin{equation*}
    (\bc\otimes M)_i=c_iM_i,
    \quad
    (R^{-1}\otimes M)_i=R_i^{-1}M_i,
    \quad
    \proj_\cC(M)_i=\proj_{\cC_i}(M_i),
    \quad
    \proj_\cX(M)_i=\proj_{\cX_i}(M_i),
\end{equation*}
and we write~$\bc\otimes S=(\bc\otimes X,\bc\otimes V)$.
\Cref{algo:batched-h-osqp} states one iteration.
Lines~1 to~4 are the relaxed step of~\Cref{algo:lsf-osqp}, applied to all columns at once.
Line~5 applies the \emph{acceleration scheme}~$\Gamma_i$ of each column, which is either the Halpern update or the identity:
\begin{equation*}
    \Gamma_i(\hat s,s^0,k)=\frac{k+1}{k+2}\,\hat s+\frac{1}{k+2}\,s^0
    \qquad\text{or}\qquad
    \Gamma_i(\hat s,s^0,k)=\hat s.
\end{equation*}

\begin{algorithm}
    \begin{algorithmic}[1]
        \Statex \textbf{Input:} State~$S=(X,V)$, anchors~$S^{0}$, and counters~$k_1,\ldots,k_B$
        \State $Z\gets\proj_{\cC}(V)$; $Y\gets \brho \otimes (V-Z)$
        \State $\bar X\gets \proj_{\cX} \big( X-R^{-1} \otimes \!\left[PX+Q+A^T\!\left(Y+\brho \otimes (AX-Z)\right)\right] \big)$
        \State $\bar S\gets(\bar X,\;A\bar X+\brho^{-1}\otimes Y)$
        \State $\hat S\gets(\ones-\balpha) \otimes S+\balpha \otimes \bar S$
        \State $S_i\gets \Gamma_i\big(\hat S_i, S_i^{0}, k_i\big)$; $k_i\gets k_i+1$, for~$i=1,\ldots,B$
    \end{algorithmic}
    \caption{One iteration of batched linear-system-free ADMM}
    \label{algo:batched-h-osqp}
\end{algorithm}

Within an iteration the columns exchange no information, so each column evolves exactly as a single unbatched run would.
Every~$K$ iterations, the solver evaluates the termination criteria of~\Cref{subsec:termination-criteria} and the restart conditions~\eqref{eq:restart-conditions} for all columns in parallel.
We call this a \emph{check}.
Restarts and step-size updates act on the columns that meet a restart condition and leave the other columns unchanged.

When the columns hold different instances, such as the instances of a parametric problem, the solver runs until every column satisfies the termination criteria.
A converged instance keeps its stored solution, but its column stays in the batch until the solver \emph{repacks} it, which happens once the number of unconverged instances has fallen to half of its value at the previous repack.
A repack removes the columns of the converged instances, which makes the remaining iterations cheaper.
It can also give the freed columns to the remaining instances as additional variants of the ensemble described next, each started from the best column of its instance.

\subsection{Algorithm hedging}
\label{sec:hedging}
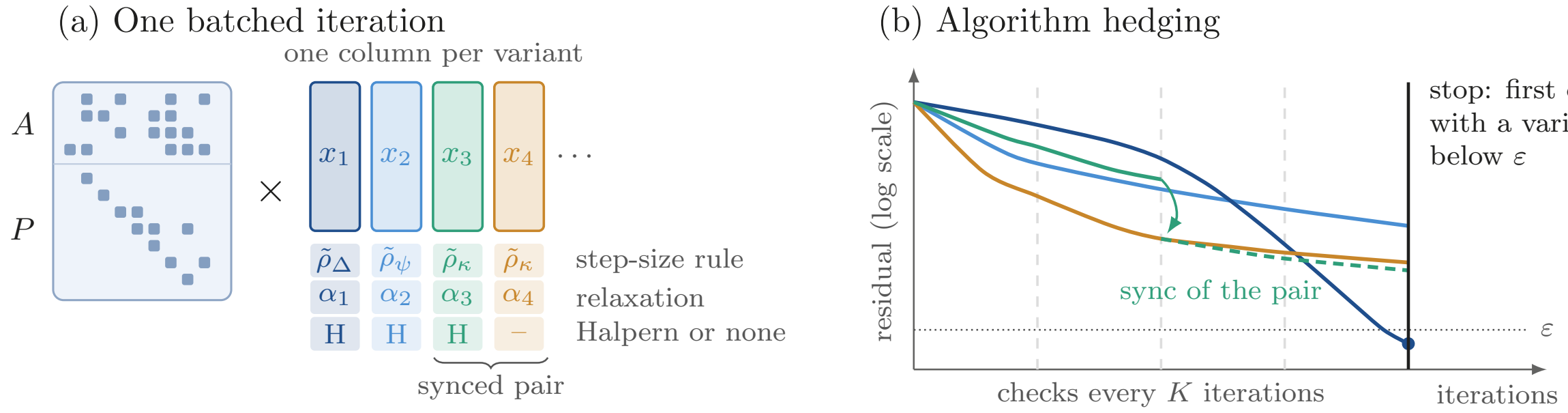
\begin{figure}[!ht]
    \centering
    \resizebox{\linewidth}{!}{\definecolor{hdgInk}{HTML}{1F1F1E}\definecolor{hdgMuted}{HTML}{555555}\definecolor{hdgGrid}{HTML}{DDDDDD}
\definecolor{hdgBand}{HTML}{EEF3FA}
\definecolor{hdgA}{HTML}{1F4E8C}\definecolor{hdgB}{HTML}{4A8FD3}\definecolor{hdgC}{HTML}{2F9E7A}\definecolor{hdgD}{HTML}{C8862A}
\begin{tikzpicture}[x=1cm,y=1cm,font=\footnotesize, >=latex,
  lab/.style={font=\scriptsize, text=hdgMuted},
  chip/.style={rounded corners=1.2pt, minimum width=.5cm, minimum height=.34cm, inner sep=0pt, font=\scriptsize},
]
\node[anchor=west, font=\small, text=hdgInk] at (-0.1,4.1) {(a) One batched iteration};
\fill[hdgBand, rounded corners=3pt] (0,1.3) rectangle (1.8,3.5);
\draw[hdgA!50, rounded corners=3pt, line width=.6pt] (0,1.3) rectangle (1.8,3.5);
\draw[hdgA!30] (0,2.675) -- (1.8,2.675);
\fill[hdgA!55, rounded corners=.7pt] (0.283,3.272) rectangle ++(0.115,0.115);
\fill[hdgA!55, rounded corners=.7pt] (0.623,3.272) rectangle ++(0.115,0.115);
\fill[hdgA!55, rounded corners=.7pt] (1.132,3.272) rectangle ++(0.115,0.115);
\fill[hdgA!55, rounded corners=.7pt] (1.472,3.272) rectangle ++(0.115,0.115);
\fill[hdgA!55, rounded corners=.7pt] (0.283,3.103) rectangle ++(0.115,0.115);
\fill[hdgA!55, rounded corners=.7pt] (0.453,3.103) rectangle ++(0.115,0.115);
\fill[hdgA!55, rounded corners=.7pt] (0.963,3.103) rectangle ++(0.115,0.115);
\fill[hdgA!55, rounded corners=.7pt] (1.132,3.103) rectangle ++(0.115,0.115);
\fill[hdgA!55, rounded corners=.7pt] (0.623,2.933) rectangle ++(0.115,0.115);
\fill[hdgA!55, rounded corners=.7pt] (0.963,2.933) rectangle ++(0.115,0.115);
\fill[hdgA!55, rounded corners=.7pt] (1.132,2.933) rectangle ++(0.115,0.115);
\fill[hdgA!55, rounded corners=.7pt] (1.302,2.933) rectangle ++(0.115,0.115);
\fill[hdgA!55, rounded corners=.7pt] (0.113,2.763) rectangle ++(0.115,0.115);
\fill[hdgA!55, rounded corners=.7pt] (0.283,2.763) rectangle ++(0.115,0.115);
\fill[hdgA!55, rounded corners=.7pt] (1.132,2.763) rectangle ++(0.115,0.115);
\fill[hdgA!55, rounded corners=.7pt] (1.302,2.763) rectangle ++(0.115,0.115);
\fill[hdgA!55, rounded corners=.7pt] (1.472,2.763) rectangle ++(0.115,0.115);
\fill[hdgA!55, rounded corners=.7pt] (0.283,2.472) rectangle ++(0.115,0.115);
\fill[hdgA!55, rounded corners=.7pt] (0.793,2.132) rectangle ++(0.115,0.115);
\fill[hdgA!55, rounded corners=.7pt] (0.453,2.302) rectangle ++(0.115,0.115);
\fill[hdgA!55, rounded corners=.7pt] (0.793,1.962) rectangle ++(0.115,0.115);
\fill[hdgA!55, rounded corners=.7pt] (1.472,1.622) rectangle ++(0.115,0.115);
\fill[hdgA!55, rounded corners=.7pt] (1.302,1.962) rectangle ++(0.115,0.115);
\fill[hdgA!55, rounded corners=.7pt] (0.623,2.132) rectangle ++(0.115,0.115);
\fill[hdgA!55, rounded corners=.7pt] (1.302,1.452) rectangle ++(0.115,0.115);
\fill[hdgA!55, rounded corners=.7pt] (0.963,1.792) rectangle ++(0.115,0.115);
\fill[hdgA!55, rounded corners=.7pt] (1.132,1.622) rectangle ++(0.115,0.115);
\fill[hdgA!55, rounded corners=.7pt] (0.963,1.962) rectangle ++(0.115,0.115);
\node[text=hdgInk, anchor=east] at (-0.05,3.08) {$A$};\node[text=hdgInk, anchor=east] at (-0.05,2.02) {$P$};
\node[text=hdgInk] at (2.2,2.4) {\large$\times$};
\foreach \i/\c/\n in {0/hdgA/1,1/hdgB/2,2/hdgC/3,3/hdgD/4}{
  \pgfmathsetmacro{\xl}{2.6+.62*\i}
  \fill[\c!20, rounded corners=1.5pt] (\xl,2.0) rectangle ({\xl+.5},3.5);
  \draw[\c, line width=.7pt, rounded corners=1.5pt] (\xl,2.0) rectangle ({\xl+.5},3.5);
  \node[text=\c] at ({\xl+.25},2.75) {$x_\n$};
}
\node[text=hdgMuted] at (5.3,2.75) {$\cdots$};
\node[lab] at (3.85,3.78) {one column per variant};
\foreach \i/\c/\n/\h/\r in {0/hdgA/1/H/\Delta,1/hdgB/2/H/\psi,2/hdgC/3/H/\kappa,3/hdgD/4/{--}/\kappa}{
  \pgfmathsetmacro{\xc}{2.85+.62*\i}
  \node[chip, fill=\c!14, text=\c] at (\xc,1.7) {$\tilde\rho_{\r}$};
  \node[chip, fill=\c!14, text=\c] at (\xc,1.33) {$\alpha_\n$};
  \node[chip, fill=\c!14, text=\c] at (\xc,0.96) {\h};
}
\node[lab, anchor=west] at (5.15,1.7) {step-size rule};
\node[lab, anchor=west] at (5.15,1.33) {relaxation};
\node[lab, anchor=west] at (5.15,0.96) {Halpern or none};
\draw[decorate, decoration={brace, mirror, amplitude=3pt}, hdgMuted, line width=.5pt] (3.85,0.72) -- (5.0,0.72);
\node[lab] at (4.43,0.42) {synced pair};

\begin{scope}[xshift=8.7cm]
\node[anchor=west, font=\small, text=hdgInk] at (-0.5,4.1) {(b) Algorithm hedging};
\draw[->, hdgInk!70, line width=.5pt] (0,0.6) -- (6.4,0.6);
\node[lab] at (5.9,0.35) {iterations};
\draw[->, hdgInk!70, line width=.5pt] (0,0.6) -- (0,3.65);
\node[lab, rotate=90] at (-0.3,2.1) {residual (log scale)};
\foreach \x in {1.25,2.5,3.75,5.0}{ \draw[hdgGrid, line width=.7pt, dashed] (\x,0.6) -- (\x,3.5); }
\node[lab] at (2.5,0.35) {checks every $K$ iterations};
\draw[hdgMuted, line width=.5pt, densely dotted] (0,1.0) -- (6.2,1.0) node[right, lab] {$\varepsilon$};
\draw[hdgB, line width=1.1pt] plot[smooth] coordinates {(0,3.3) (0.8,2.85) (1.25,2.68) (2.5,2.42) (3.75,2.22) (5.0,2.05)};
\draw[hdgA, line width=1.1pt] plot[smooth] coordinates {(0,3.3) (0.8,3.16) (1.25,3.07) (2.0,2.9) (2.5,2.73) (3.0,2.45) (3.5,2.07) (4.0,1.65) (4.4,1.3) (4.75,1.0) (5.0,0.86)};
\draw[hdgD, line width=1.1pt] plot[smooth] coordinates {(0,3.3) (0.7,2.6) (1.25,2.35) (2.0,2.05) (2.5,1.92) (3.2,1.84) (3.75,1.78) (5.0,1.68)};
\draw[hdgC, line width=1.1pt] plot[smooth] coordinates {(0,3.3) (0.9,2.95) (1.25,2.85) (2.0,2.62) (2.5,2.52)};
\draw[->, hdgC, line width=.7pt] (2.5,2.52) to[bend left=40] (2.56,1.97);
\draw[hdgC, line width=1.1pt, densely dashed] plot[smooth] coordinates {(2.5,1.92) (3.2,1.8) (3.75,1.72) (5.0,1.6)};
\node[lab, anchor=west, text=hdgC] at (1.95,1.38) {sync of the pair};
\fill[hdgA] (5.0,0.86) circle (.065);
\draw[hdgInk, line width=.9pt] (5.0,0.6) -- (5.0,3.5);
\node[lab, anchor=west, align=left, text=hdgInk] at (5.08,3.1) {stop: first check\\with a variant\\below $\varepsilon$};
\end{scope}
\end{tikzpicture}
 }
    \caption{Batched iteration and algorithm hedging.
    (a) One iteration of the batch reads~$A$ and~$P$ once for all columns, and each column runs a variant with its own step-size rule, relaxation parameter, and acceleration scheme.
    Here~$\tilde\rho_\Delta$, $\tilde\rho_\psi$, and~$\tilde\rho_\kappa$ are the candidate step sizes of the displacement, residual, and curvature rules of~\Cref{tab:step-rules}.
    (b) Schematic residuals of the variants: the variants run side by side, the members of a synced pair synchronize at the checks, and the solver stops at the first check at which some variant meets the tolerance.}
    \label{fig:batching}
\end{figure}

Algorithm hedging uses the~$E$ columns of an instance in a different way: they hold the same instance, and each of them runs a different algorithm.
To simplify the exposition, this section considers a single instance,~$N_\mathrm{inst}=1$, so that~$B=E$.
We call each algorithm that the solver runs on the instance an \emph{algorithm variant}, and the set of variants that run together an \emph{ensemble}.
No single variant performs best on every instance, and the best one is rarely known in advance.
The solver stops at the first check at which some variant satisfies the termination criteria.
It thus runs only as many iterations as the best variant of the ensemble needs, without knowing in advance which variant that is.
\Cref{fig:batching} illustrates one batched iteration and algorithm hedging.

The variants are not arbitrary algorithms.
To run as one batch, they must execute the same sequence of operations, with the same sparse products, and differ only in data that each column carries.
Every variant therefore applies the operator~$\opF$ with the shared matrices~$P$ and~$A$, and tests the same restart conditions~\eqref{eq:restart-conditions} on its own residual.
The variants differ in three components, each of which has the same form in every column.
\begin{itemize}
    \item \emph{Acceleration scheme.}
        Every scheme~$\Gamma_i$ is an affine combination of the relaxed iterate and the anchor, with coefficients that depend only on the counter~$k_i$.
        A variant chooses the scheme and the relaxation parameter~$\alpha_i$.
    \item \emph{Step-size rule.}
        Every variant updates~$\rho_i$ by the same formula, from a candidate step size that it computes from its own iterates.
        A variant chooses the rule for the candidate and the constants of the formula.
    \item \emph{Synced pairs.}
        Variants exchange information only at the checks, by copying columns.
        Two variants may form a pair that does so.
\end{itemize}
Since each component acts on each column separately, a single GPU kernel executes it for the whole ensemble without branching, and the sparse products of~\Cref{algo:batched-h-osqp} are unchanged.
This common structure is what makes the variants batchable.
The data preconditioning, in contrast, changes the matrices~$P$ and~$A$ that all columns share, so it is the same for all variants.
\Cref{subsec:variants} describes the three components, and~\Cref{subsec:hedging-algorithm} states the resulting algorithm.

\subsection{Algorithm variants}
\label{subsec:variants}
\paragraph{Acceleration scheme}
\label{subsec:acceleration-scheme}
Each variant fixes its acceleration scheme~$\Gamma$, the Halpern update or the identity, and its relaxation parameter~$\alpha$.
\begin{itemize}
    \item \emph{Halpern variants} use~$\alpha=2$, so they run the Halpern iteration on the reflected operator~$2\opF-\identity$, as the reflected restarted Halpern PDHG of~\cite{luRestartedHalpernPDHG2024} does.
        They have the worst-case guarantees of~\Cref{subsec:sharpness-free}, and their certificate candidates converge by~\Cref{prop:candidate-convergence}($b$).
    \item \emph{Halpern-free variants} run the plain iteration~$s^{k+1}=\opT s^k$.
        Its worst-case guarantee is weaker, but on unconstrained strongly convex QPs it can match the worst-case order of restarted Halpern iteration for fixed~$\alpha\in(0,1]$ (\Cref{app:tightness}).
        In practice, the plain iteration is typically much faster on problems with local strong convexity.
        On LPs, these variants use~$\alpha=1.8$, for which~$\opT$ is averaged (\Cref{cor:lsf-averagedness}($a$)), so the iteration converges and so do its certificate candidates (\Cref{prop:candidate-convergence}($a$)).
        On QPs, they use~$\alpha=2$, for which~$\opT$ is nonexpansive but in general not averaged, so the iteration need not converge.
        We include them because~$\opT$ is often contractive on strongly convex QPs.
\end{itemize}

\paragraph{Step-size rule}
\label{subsec:step-sizes}
The step size~$\rho$ sets the relative length of the primal and the dual steps: a larger~$\rho$ moves the dual iterates faster and the primal iterates slower.
Each variant starts from its own initial step size and updates it at its restarts with its own rule.
Every rule computes a \emph{candidate step size}~$\tilde\rho$ from the iterates of the variant, and the variant moves its step size toward the candidate by
\begin{equation}
    \label{eq:step-update}
    \rho\gets\big(\xi\tilde\rho\big)^{\theta}\rho^{1-\theta},
    \qquad
    \theta=\exp\!\big(-\zeta\,\omega^k/\omega^\mathrm{first}\big).
\end{equation}
The multiplier~$\xi>0$ scales the candidate, and the weight~$\theta\in(0,1]$ smooths its fluctuations between consecutive updates.
Here~$\omega^k$ is the residual~\eqref{eq:restart-residual}, $\omega^\mathrm{first}$ is its value at the first check, and~$\xi$ and~$\zeta>0$ are constants of the variant.
The weight increases toward~$1$ as the residual decreases, so the step size follows the candidate more closely near a solution.

The candidates combine two quantities.
Let~$\Delta\tilde x$, $\Delta y$, $\Delta z$, and~$\Delta w$ be the changes of~$\tilde x$, $y$, $z$, and~$w$ since the last restart, on the scaled data.
The \emph{displacement ratio} and the \emph{residual ratio} are
\begin{equation*}
    \rho_\Delta=\frac{\|\Delta y\|}{\sqrt{\lambda_A}\,\|\Delta\tilde x\|},
    \qquad
    \psi=\sqrt{r_\mathrm{prim}^\mathrm{rel}/r_\mathrm{dual}^\mathrm{rel}},
\end{equation*}
where~$\lambda_A$ is the constant of~\Cref{subsec:derivation}, and
\begin{equation*}
    r_\mathrm{prim}^\mathrm{rel}=\frac{\|r_\mathrm{prim}\|_\infty}{1+\max\{\|Ax\|_\infty,\|x\|_\infty\}},
    \qquad
    r_\mathrm{dual}^\mathrm{rel}=\frac{\|r_\mathrm{dual}\|_\infty}{1+\max\{\|Px\|_\infty,\|A^Ty+w\|_\infty,\|q\|_\infty\}}
\end{equation*}
are the relative primal and dual residuals of~\Cref{subsec:termination-criteria} in the original variables.
The displacement ratio balances how far the primal and the dual iterates have moved, since it minimizes~$\rho\lambda_A\|\Delta\tilde x\|^2+\rho^{-1}\|\Delta y\|^2$ over~$\rho>0$.
The residual ratio exceeds~$1$ when the primal residual dominates, which calls for a larger step size, and it is below~$1$ when the dual residual dominates.
\Cref{tab:step-rules} lists the rules.
The cross-term rule uses
\begin{align*}
    c_\mathrm{d}&=\big|\Delta y^T(A\Delta\tilde x+\lambda_A\Delta z-\lambda_A\Delta y)\big|, \\
    c_\mathrm{p}&=\big|\Delta\tilde x^TP\Delta\tilde x+(A\Delta\tilde x)^T(\Delta z-\Delta y)-\Delta z^T\Delta y+\lambda_A\|\Delta z\|^2\big|,
\end{align*}
and the curvature rule uses
\begin{equation*}
    \kappa=\frac{\sum_i\big(\sum_j|P_{ij}|\big)(\Delta\tilde x_i)^2}{\|\Delta\tilde x\|^2},
    \qquad
    \sigma=\frac{\|A^T\Delta y\|}{\|\Delta y\|},
\end{equation*}
with~$\kappa=1$ on LPs.
Two further rules, each used by one QP variant, add a multiple of~$\psi$ to a multiple of~$\rho_\Delta$.
The displacement rule is the most common one in the ensemble.
Its variant with hold rejects a decrease of~$\rho$ while the dual residual is the smallest of the residuals, and an increase while it is the largest.
The curvature rule is the only one that updates~$\rho$ at every check, whether or not the variant restarts, because the Halpern-free variants restart rarely.
Any other candidate fits the update~\eqref{eq:step-update}, so a new rule extends the ensemble without changes to the rest of the solver.

\begin{table}
    \centering
    \small
    \caption{Step-size rules and their candidate step sizes~$\tilde\rho$.}
    \label{tab:step-rules}
    \begin{tabular}{lll}
        \toprule
        Rule & $\tilde\rho$ & Description \\
        \midrule
        Displacement & $\rho_\Delta$ & balances the primal and the dual displacements \\
        \quad with bound multipliers & $\rho_\Delta$ & $\|\Delta y\|^2+\|\Delta w\|^2$ replaces~$\|\Delta y\|^2$ in~$\rho_\Delta$ \\
        \quad with hold & $\rho_\Delta$ & rejects an update that moves~$\rho$ against the residuals \\
        Residual & $\rho\,\psi$ & increases~$\rho$ when the primal residual dominates \\
        Cross-term & $\sqrt{c_\mathrm{d}/c_\mathrm{p}}$ & keeps the products through~$A$ and~$P$ \\
        Curvature & $\kappa/(2\sigma\sqrt{\lambda_A})$ & curvature~$\kappa$ of the objective against the gain~$\sigma$ of~$A^T$ \\
        \bottomrule
    \end{tabular}
\end{table}

\paragraph{Synced pairs}
\label{subsec:synchronization}
A \emph{synced pair} consists of a Halpern variant and a Halpern-free variant that are otherwise equal, except possibly for~$\alpha$.
Which of the two makes more progress can change during a run.
The Halpern-free member is typically faster where the problem is locally strongly convex, while the Halpern member secures the worst-case rate.
Synchronization lets the pair follow the member that currently makes more progress.
At every check, before the restart test, the solver compares the residuals~$\omega^k$ of~\eqref{eq:restart-residual} of the two members.
The member with the larger residual copies the iterate, the anchor, the counter, the step size, and the restart residuals of the other, with ties resolved in favor of the second member of the pair.
Both members therefore make the same restart decision.
Since the benefit comes from combining the two members, a pair is split only at the one width where no candidate can take its place (\Cref{subsec:variant-order}).

\subsection{Hedging algorithm}
\label{subsec:hedging-algorithm}
At each check, the solver measures every variant by its \emph{normalized KKT error}
\begin{equation}
    \label{eq:kkt-error}
    e=\max\left\{
        \frac{\norm[\infty]{r_{\mathrm{prim}}}}{\varepsilon_{\mathrm{prim}}},\,
        \frac{\norm[\infty]{r_{\mathrm{dual}}}}{\varepsilon_{\mathrm{dual}}},\,
        \frac{|r_{\mathrm{gap}}|}{\varepsilon_{\mathrm{gap}}}
    \right\},
\end{equation}
the largest ratio of a residual of~\Cref{subsec:termination-criteria} to its tolerance.
It is at most~$1$ if and only if the variant satisfies the optimality criteria~\eqref{eq:termination}, so it compares variants by how close they are to termination.
\Cref{algo:hedging} states algorithm hedging, where~$e_i$ is the normalized KKT error and~$\omega_i$ the residual~\eqref{eq:restart-residual} of variant~$i$.
Its iterations in line~2 include the acceleration scheme~$\Gamma_i$ of every variant, which is line~5 of~\Cref{algo:batched-h-osqp}.

\begin{algorithm}[!htb]
    \begin{algorithmic}[1]
        \Statex \textbf{Input:} Ensemble of~$B$ variants, check interval~$K$
        \Repeat
            \State run~$K$ iterations of~\Cref{algo:batched-h-osqp}
            \If{some variant satisfies~\eqref{eq:primal-infeasibility} or~\eqref{eq:dual-infeasibility}}
                \Comment{infeasibility}
                \State \Return its certificate
            \EndIf
            \State $j\gets\argmin_i e_i$
                \Comment{termination}
            \If{$e_j\le1$}
                \State \Return $(x_j,y_j,w_j)$
            \EndIf
            \For{every synced pair~$(i,i')$ with~$\omega_{i'}\le\omega_i$}
                \Comment{synchronization}
                \State $(S_i,S_i^0,k_i,\rho_i,\omega_i^0,\omega_i^\mathrm{prev})\gets(S_{i'},S_{i'}^0,k_{i'},\rho_{i'},\omega_{i'}^0,\omega_{i'}^\mathrm{prev})$
            \EndFor
            \For{every variant~$i$ that satisfies~\eqref{eq:restart-conditions}}
                \Comment{restart}
                \State update~$\rho_i$ by~\eqref{eq:step-update}; $S_i\gets\bar S_i$; $S_i^0\gets\bar S_i$; $k_i\gets0$
            \EndFor
        \Until{the iteration or time limit is reached}
        \State \Return $(x_j,y_j,w_j)$
    \end{algorithmic}
    \caption{Algorithm hedging}
    \label{algo:hedging}
\end{algorithm}

The infeasibility criteria do not depend on the variant that produced the certificate candidates, so the ensemble hedges infeasibility detection as well.
If a norm or an inner product in~\eqref{eq:kkt-error} is not finite for some variant, the solver sets~$e_i=+\infty$ for this variant, without affecting the others.
Ties in the choice of~$j$ are broken by the lowest index.
A variant that uses the curvature rule also updates~$\rho_i$ at the checks at which it does not restart, which~\Cref{algo:hedging} does not show.

\section{Implementation}
\label{sec:implementation}
We implement the method as \emph{\osqplf{}}, a new linear-system-free mode of the OSQP solver, which runs on GPUs.
The code and its documentation are available at \url{https://osqp.org}.
Each iteration of~\Cref{algo:batched-h-osqp} needs two sparse products: one with the stacked matrix of~$A$ and~$P$, and one with~$A^T$ (\Cref{app:algebra-engine}).
An algebra engine evaluates them and selects at setup, for each product, between the cuSPARSE library~\cite{nvidiaCuSPARSELibrary2019} and custom kernels (\Cref{app:algebra-engine}).
The rest of this section gives the default settings of the solver: the scaling of the problem data, the restart schedule, and the choice of the ensemble.

\subsection{Data preconditioning}
\label{subsec:preconditioning}
As is standard for first-order solvers, we scale the problem data with positive diagonal matrices~$D_\mathrm{c}\in\reals^{d\times d}$ and~$D_\mathrm{r}\in\reals^{m\times m}$, which act on the columns and on the rows of~$A$.
The solver runs on the scaled data
\begin{equation*}
    \hat P=D_\mathrm{c}PD_\mathrm{c},
    \quad
    \hat q=D_\mathrm{c}q,
    \quad
    \hat A=D_\mathrm{r}AD_\mathrm{c},
    \quad
    [\hat l,\hat u]=[D_\mathrm{r}l,D_\mathrm{r}u],
    \quad
    [\hat l_x,\hat u_x]=[D_\mathrm{c}^{-1}l_x,D_\mathrm{c}^{-1}u_x],
\end{equation*}
which correspond to the change of variables~$x=D_\mathrm{c}\hat x$.
The solver evaluates the termination criteria of~\Cref{subsec:termination-criteria} in the original variables.

We compute~$D_\mathrm{c}$ and~$D_\mathrm{r}$ in a sequence of passes, starting from the identity.
Each pass computes a factor~$\delta_{\mathrm{r},i}>0$ for every row~$i$ and a factor~$\delta_{\mathrm{c},j}>0$ for every column~$j$ from the current scaled matrices~$\hat A$ and~$\hat P$, and divides the corresponding diagonal entries of~$D_\mathrm{r}$ and~$D_\mathrm{c}$ by them.
A zero row or column has factor~$1$.
We use three kinds of passes, in this order.
\begin{itemize}
    \item \emph{Geometric mean scaling~\cite{tomlinScalingLinearProgramming1975}, 12 passes, for LPs only.}
        The row factor is~$\delta_{\mathrm{r},i}=(\min_j|\hat A_{ij}|\,\max_j|\hat A_{ij}|)^{1/2}$, where the minimum and the maximum are over the nonzero entries of row~$i$.
        The column factors are the same quantities for the columns of the matrix scaled by the new row factors.
    \item \emph{Ruiz equilibration~\cite{ruizScalingAlgorithmEquilibrate2001}, 10 passes.}
        The factors are~$\delta_{\mathrm{r},i}=\|\hat A_{i,:}\|_\infty^{1/2}$ and~$\delta_{\mathrm{c},j}=\max\{\|\hat A_{:,j}\|_\infty,\|\hat P_{:,j}\|_\infty\}^{1/2}$, where~$\hat A_{i,:}$ and~$\hat A_{:,j}$ are row~$i$ and column~$j$ of~$\hat A$.
    \item \emph{Pock--Chambolle scaling~\cite{pockDiagonalPreconditioningFirst2011} with exponent~$1$, 1 pass.}
        The factors are~$\delta_{\mathrm{r},i}=\|\hat A_{i,:}\|_1^{1/2}$ and~$\delta_{\mathrm{c},j}=(\|\hat A_{:,j}\|_1+\|\hat P_{:,j}\|_1)^{1/2}$.
\end{itemize}

\paragraph{Limit on the scaling factors}
A variable with finite lower and upper bounds already has a scale, given by its box.
If the passes above scale such a variable by several orders of magnitude, the rows in which it appears receive very large factors.
Their dual variables then move much faster than the others, which can stall the solver.
We therefore keep some entries of~$D_\mathrm{c}$ and~$D_\mathrm{r}$ in the interval~$[10^{-2},10^{2}]$ during the first two kinds of passes.
We call a variable \emph{bounded} if both of its bounds are finite, and a row of~$A$ \emph{wide} if its largest nonzero entry is at least~$10^4$ times its smallest in magnitude.
In the geometric mean passes, the limit applies to the bounded variables and to the rows that contain one.
In the Ruiz passes, it applies to every column except those of the variables that are not bounded and appear in a wide row, which may need a large factor.
The Pock--Chambolle pass has no limit.

\subsection{Restarts}
\label{subsec:restart-conditions}
Each variant tests the restart conditions~\eqref{eq:restart-conditions} on its own iterates, with its own residual~$\omega^k$ of~\eqref{eq:restart-residual}.
The tests take place every~$K=100$ iterations, together with the termination tests of~\Cref{algo:hedging}.
Hence~$\omega^\mathrm{prev}=\omega^{k-K}$, or~$+\infty$ if the current outer iteration has no previous test.
At the first test every variant restarts, and its residual~$\omega^\mathrm{first}$ serves as a reference scale for the rest of the run.
We use~$\beta=0.2$ and~$\beta_\mathrm{nec}=0.8$.
The fraction~$\beta_\mathrm{long}$ equals~$0.5$ if~$\omega^{k}>0.1\,\omega^\mathrm{first}$ and~$0.2$ otherwise, so restarts become more frequent once the residual is small.
Since~$\opT-\identity=\alpha(\opF-\identity)$, the residual~$\omega^k$ equals~$\alpha^{-1}\norm[\opW]{\opT s^k-s^k}$ in the metric~$\opW$ of~\Cref{subsec:nonexpansiveness}, so the sufficient-decay condition is the adaptive rule of~\Cref{thm:restart}($ii$), tested every~$K$ iterations only.

At a restart, the variant updates its step size by~\eqref{eq:step-update}, resets~$k$ to~$0$, and sets both the iterate and the anchor~$s^0$ to the restart point~$\opF s^k$ of~\Cref{algo:restarted-halpern}.
The solver evaluates this point as~$(\xtkk,\,A\xtkk+\rho^{-1}y^{k})$, where~$\xtkk$ and~$y^k$ are the outputs of the last step and~$\rho$ is the updated step size.

\subsection{Ensemble of algorithm variants}
\label{subsec:hedging-implementation}
For each problem class, \osqplf{} has a fixed list of~$16$ candidate variants.
The candidates combine the components of~\Cref{subsec:variants}, and each list contains one synced pair.\footnote{The OSQP documentation lists the candidate variants and their components.}
The ensemble of an instance is determined by two choices: the order in which the candidates enter it, and its size~$E$, which we call the \emph{ensemble size}.

\paragraph{Variant order}
\label{subsec:variant-order}
An instance with ensemble size~$E$ runs the first~$E$ candidates of a fixed \emph{variant order}.
The only exception occurs when the~$E$th column would separate the synced pair of~\Cref{subsec:synchronization}: the solver then skips the pair and runs the next candidate, or, when no candidate follows the pair, its first member alone.
We compute the order offline, separately for QPs and LPs, from~$N$ randomly perturbed copies of the MIPLIB, Mittelmann, Maros--M\'esz\'aros, and SuiteSparse instances, which~\Cref{app:variant-order} describes.
We choose the order that minimizes the total run time over these instances, where we approximate the run time of an ensemble by the smallest run time among its members.

\emph{Data.}
The synced pair always enters an ensemble as a whole, so we treat it separately from the other~$M=14$ candidates, which we call \emph{single}.
We run every candidate on its own on every one of these instances.
The matrix~$\Theta\in\reals^{N\times M}$ collects the run times of the single candidates, with~$\Theta_{ni}$ the run time of single candidate~$i$ on instance~$n$, and the vector~$\theta\in\reals^N$ collects the run times of the pair.
We also record the ensemble size~$E_n$ that the solver selects for instance~$n$ at setup.

\emph{Variables.}
The order is described by two variables.
The permutation matrix~$\Pi\in\{0,1\}^{M\times M}$ orders the single candidates: $\Pi_{ij}=1$ if single candidate~$i$ is the~$j$th of them.
The vector~$p\in\{0,1\}^{M+1}$ places the pair among them: it has one nonzero entry, and~$p_j=1$ if the pair comes right after the first~$j-1$ single candidates.
Two more variables select, for every instance, the member of its ensemble whose run time counts: $\Xi_{ni}=1$ if it is single candidate~$i$, and~$\xi_n=1$ if it is the pair, with~$\Xi\in\reals^{N\times M}$ and~$\xi\in\reals^N$.

\emph{Ensembles.}
The order determines the ensemble of instance~$n$ as follows.
The pair takes two columns, so it fits in the first~$E_n$ columns if at most~$E_n-2$ single candidates come before it, \ie, if~$p_j=1$ for some~$j\le E_n-1$.
In that case, the ensemble consists of the pair and the first~$E_n-2$ single candidates.
Otherwise, it consists of the first~$E_n$ single candidates.
Three 0--1 matrices encode these three sets of positions,
\begin{equation*}
    F_{nj}=\begin{cases}1 & j\le E_n\\ 0 & \text{otherwise,}\end{cases}
    \qquad
    H_{nj}=\begin{cases}1 & j\le E_n-2\\ 0 & \text{otherwise,}\end{cases}
    \qquad
    G_{nj}=\begin{cases}1 & j\le E_n-1\\ 0 & \text{otherwise,}\end{cases}
\end{equation*}
where~$F$ and~$H$ have~$M$ columns and~$G$ has~$M+1$.
Then entry~$n$ of~$Gp$ equals~$1$ if the ensemble of instance~$n$ contains the pair.
Entry~$(n,i)$ of~$F\Pi^T$ equals~$1$ if single candidate~$i$ is among the first~$E_n$ single candidates, and that of~$H\Pi^T$ if it is among the first~$E_n-2$.

\emph{Problem.}
The variant order solves the mixed-integer linear program
\begin{equation}
    \label{prob:variant-order}
    \begin{array}{ll}
        \textnormal{minimize} & \Tr(\Theta^T\Xi)+\theta^T\xi \\
        \textnormal{subject to} & \Pi\ones=\ones, \quad \Pi^T\ones=\ones, \quad \ones^Tp=1, \\
        & \Xi\ones+\xi=\ones, \\
        & \xi\le Gp, \\
        & \Xi\le F\Pi^T, \\
        & \Xi\le H\Pi^T+(\ones-Gp)\ones^T, \\
        & \Pi\in\{0,1\}^{M\times M}, \quad p\in\{0,1\}^{M+1}, \quad \Xi\ge0, \quad \xi\ge0,
    \end{array}
\end{equation}
where the inequalities hold entrywise.
The objective is the total of the selected run times.
The first line of constraints states that~$\Pi$ and~$p$ describe an order, and the second that every instance selects one candidate.
The third line allows an instance to select the pair only if its ensemble contains it.
The fourth line allows it to select single candidate~$i$ only if~$i$ is among its first~$E_n$ single candidates, and the fifth tightens this to the first~$E_n-2$ when the ensemble contains the pair.
Since we minimize, every instance selects the fastest member of its ensemble.
The resulting orders are fixed default settings of~\osqplf{} and are not recomputed for new instances or benchmarks.

\paragraph{Batch size}
\label{subsec:batching}
The run time is the iteration count of the best variant in the ensemble times the time of one batched iteration.
A larger ensemble size~$E$ makes it more likely that the ensemble contains a variant that suits the instance, which lowers the first factor.
Once the GPU saturates, however, every additional column raises the second factor.
The saturation point depends on the instance and on the hardware, so we choose~$E$ for each instance, or for each batch of instances, at setup, in two stages.
All instances of a batch run the same variants, so the batch size is~$B=N_\mathrm{inst}E$ (\Cref{subsec:batch-mechanism}).
A hardware model first gives an upper bound~$B_\mathrm{hw}$, and timed iterations then select~$B\le B_\mathrm{hw}$, unless the batch holds more than~$B_\mathrm{hw}$ instances.

\emph{Upper bound.}
The bound is
\begin{equation*}
    B_\mathrm{hw}=\max\left\{1,\;\min\left\{B_{\max},\,B_\mathrm{thr},\,B_\mathrm{L2},\,B_\mathrm{roof}\right\}\right\},
\end{equation*}
where~$B_{\max}=N_\mathrm{inst}E_{\max}$, and the maximum ensemble size~$E_{\max}$ is a setting of the solver.
A single instance has~$N_\mathrm{inst}=1$.
The other three numbers estimate the batch sizes at which three resources of the GPU saturate.
\begin{itemize}
    \item \emph{Threads.}
        A vector operation of the batched iteration updates~$(d+m)B$ entries, one per thread.
        The GPU keeps~$N_\mathrm{thr}$ threads resident, so all entries are updated in one wave of threads if~$B\le B_\mathrm{thr}=\lfloor N_\mathrm{thr}/(d+m)\rfloor$.
    \item \emph{Cache.}
        The two sparse products read and write~$2(d+m)$ double-precision entries per column, \ie, $16(d+m)$ bytes.
        The vectors of all columns fit in the L2 cache, of~$C_\mathrm{L2}$ bytes, if~$B\le B_\mathrm{L2}=\lfloor C_\mathrm{L2}/(16(d+m))\rfloor$.
    \item \emph{Memory bandwidth.}
        In the roofline model~\cite{williamsRooflineInsightfulVisual2009}, a kernel takes the larger of two times: the bytes it moves divided by the memory bandwidth~$\pi_\mathrm{bw}$, and the floating-point operations it performs divided by the peak double-precision throughput~$\pi_\mathrm{flop}$ of the GPU.
        Let~$\nu$ be the number of nonzeros that the two sparse products of an iteration read, and~$\mu$ the number of bytes that store them.
        At batch size~$B$, the products read the matrices once and the vectors of every column, \ie, $\mu+16(d+m)B$ bytes, and they perform~$2\nu B$ operations.
        For small~$B$ the first of the two times dominates and grows slowly with~$B$, so additional columns are almost free.
        The two times are equal at
        \begin{equation*}
            B_\mathrm{roof}=\left\lfloor\frac{\mu/\pi_\mathrm{bw}}{2\nu/\pi_\mathrm{flop}-16(d+m)/\pi_\mathrm{bw}}\right\rfloor,
        \end{equation*}
        and beyond this batch size the time of an iteration grows in proportion to~$B$.
        If the denominator is not positive, the products are bound by the memory bandwidth at every batch size, and we drop this term.
\end{itemize}
The solver obtains~$N_\mathrm{thr}$, $C_\mathrm{L2}$, $\pi_\mathrm{bw}$, and~$\pi_\mathrm{flop}$ from the device at setup.

\emph{Timed selection.}
Since~$B=N_\mathrm{inst}E$, the batch size is a multiple of~$N_\mathrm{inst}$.
We then time one iteration at the batch sizes~$B_0<B_1<\cdots<B_J$, which are~$N_\mathrm{inst}$ times the powers of two~$1,2,4,\ldots$ below~$B_J$, followed by
\begin{equation*}
    B_J=N_\mathrm{inst}\max\left\{1,\;\lfloor B_\mathrm{hw}/N_\mathrm{inst}\rfloor\right\}.
\end{equation*}
For a single instance,~$B_J=B_\mathrm{hw}$.
For a batch,~$B_J$ is the largest multiple of~$N_\mathrm{inst}$ not above~$B_\mathrm{hw}$.
In a batch of more than~$B_\mathrm{hw}$ instances, every instance still runs one variant, so~$E=1$ and~$B=N_\mathrm{inst}$.
Let~$t(B)$ be the measured time of one iteration at batch size~$B$, so that~$B/t(B)$ is the throughput in column iterations per unit of time.
Starting from~$B_0$, we move from~$B_{j-1}$ to~$B_j$ as long as
\begin{equation*}
    \frac{B_j}{t(B_j)}\ge1.2^{\log_2(B_j/B_{j-1})}\,\frac{B_{j-1}}{t(B_{j-1})},
\end{equation*}
and we set~$B$ to the last batch size reached, so that~$E=B/N_\mathrm{inst}$.
Hence every doubling of the batch size must raise the throughput by at least~$20\%$, and the exponent adjusts this requirement for the last step, which need not be a doubling.
The threshold of~$20\%$ is a tuning parameter of the solver.

A repack frees columns, so at setup the solver also chooses the ensemble size~$E$ after each repack.
It repeats the timed selection for a batch of~$\lfloor N_\mathrm{inst}/2\rfloor$ instances, then of~$\lfloor N_\mathrm{inst}/4\rfloor$, and so on, with the same bound~$B_\mathrm{hw}$ and the same~$E_{\max}$.
Each of these selections starts from the ensemble size that the previous one chose, so a repack never removes a variant.
A repack that leaves fewer instances than planned uses the choice for the smallest planned batch that holds them.

\section{Experimental results}
\label{sec:experiments}

This section evaluates \osqplf{} on LP and QP benchmarks.
We first describe the experimental setup and summarize the results.
We then discuss the QP and the LP benchmarks, the batching of parametric instances, and the effect of batching and hedging.

\paragraph{Tolerance}
We run every solver with tolerances~$\varepsilon_{\textnormal{abs}} = \varepsilon_{\textnormal{rel}} \in \{10^{-4}, 10^{-8} \}$.
The solvers use different termination criteria, so the statuses they report are not comparable.
We therefore modify each open-source solver to use the termination criteria of this paper, and leave everything else unchanged.
Gurobi cannot be modified in this way.
When the solution it returns does not satisfy our criteria, we multiply its tolerances by~$0.9$ and solve again.
We stop when the solution satisfies our criteria, when Gurobi reports an error or reaches the time limit, or when its tolerances reach the smallest value it accepts, $10^{-9}$.
The time we report for Gurobi is that of its last solve.

\paragraph{Settings of \osqplf{}}
We run \osqplf{} with the settings of~\Cref{sec:implementation} and a maximum ensemble size~$E_{\max}=16$.
The solver then chooses the ensemble size of each instance at setup (\Cref{subsec:batching}).

\paragraph{Solvers}
We use Gurobi with crossover and presolve disabled, the C implementations of HPR-LP and HPR-QP, and the GPU implementation of Clarabel.\footnote{We use Gurobi~13.0.2, HPR-LP~0.1.3, and HPR-QP~0.1.0. For cuPDLPx and cuOpt we use development versions after their releases~0.3.0 and~26.10.00a, and for Clarabel the CuClarabel branch of Clarabel.jl~0.11.0. The repository \url{https://github.com/osqp/osqp2-benchmarks} lists the commit of every open-source solver.}
Every solver uses its default infeasibility tolerances and no iteration limit, so a run ends at optimality, with an infeasibility certificate, or at the time limit.
The time of \osqplf{} includes its setup, and with it the selection of the ensemble size of~\Cref{subsec:batching}.

\paragraph{Time limit}
As in~\cite{luCuPDLPxFurtherEnhanced2025,chenHPRLPImplementationHPR2025}, the time limit depends on the size of the instance.
The instances with~$\nnz(P) + \nnz(A) > 10^7$ have a time limit of 18,000 seconds, and all other instances 3,600 seconds.
We also use 18,000 seconds on the Mittelmann set, where these two papers use 15,000.

\paragraph{Shifted geometric mean}
We compare solve times with the shifted geometric mean
\begin{equation*}
    g_s = \exp\bigg(\frac{1}{|\mathcal{P}|}\sum_{p\in\mathcal{P}} \ln(t_{p,s} + k)\bigg) - k,
\end{equation*}
where~$t_{p,s}$ is the time solver~$s$ takes on instance~$p$, $\mathcal{P}$ is the set of instances, and~$k=10$ is the shift.
If solver~$s$ fails on instance~$p$, we set~$t_{p,s}$ to the time limit of that instance.
A solver fails if it reaches the time limit, reports a numerical error or a status other than optimal, or exceeds~$450$ GB of memory.

\paragraph{Presolving}
We run every solver on the original instances and on the instances obtained after Gurobi's presolve routine.
The presolved sets are smaller for SuiteSparse and MIPLIB, by 3 and 19 instances, because we drop the instances that presolve eliminates or cannot process.

\paragraph{Computing environment}
We run all experiments on a single node with two 32-core AMD EPYC 9334 processors and an NVIDIA H100 GPU with 94 GB of memory, using CUDA 13.3.
Each solver may use one GPU and at most 32 CPU cores.
The code to reproduce the experiments, including our modifications to the other solvers, is available at
\begin{center}
    \url{https://github.com/osqp/osqp2-benchmarks}.
\end{center}

\paragraph{Size groups}
We split each benchmark set into small, medium, and large instances by the number of nonzeros~$\nnz(P) + \nnz(A)$ of the original instance.
The thresholds are~$10^4$ and~$10^5$ nonzeros for QPs, and~$10^6$ and~$10^7$ for LPs.
\Cref{app:detailed-results} reports the number of solved instances and the shifted geometric mean of every solver, for each size group, on the original and on the presolved instances.

\begin{table}[!ht]
    \centering
    \caption{Summary of the benchmark results.
    Ratio is the time of the other solver divided by the time of \osqplf{}, so a ratio above~$1$ favors \osqplf{}.
    Gurobi runs on the CPU and every other solver on the GPU.}
    \label{tab:summary}

    \begin{subtable}{\linewidth}
    \centering
    \caption{LP and QP benchmark sets without presolve, with the number of instances in parentheses.
    For each tolerance~$\varepsilon$, the number of solved instances (\#) and the shifted geometric mean of the solve times in seconds (SGM) of \osqplf{}, of the fastest other first-order solver, and of the fastest interior-point solver.}
    \label{tab:summary-benchmarks}
    \footnotesize
    \setlength{\tabcolsep}{3pt}
    \begin{tabular*}{\linewidth}{@{\extracolsep{\fill}}llrrlrrrlrrr@{}}
        \toprule
        & & \multicolumn{2}{c}{\osqplf{}} & \multicolumn{4}{c}{Fastest other first-order solver} & \multicolumn{4}{c}{Fastest interior-point solver}\\
        \cmidrule{3-4}\cmidrule{5-8}\cmidrule{9-12}
        Benchmark set & $\varepsilon$ & \# & SGM & Solver & \# & SGM & Ratio & Solver & \# & SGM & Ratio\\
        \midrule
        \multirow{2}{*}{Maros--M\'esz\'aros (138)} & $10^{-4}$ & 138 & 0.76 & HPR-QP & 130 & 8.69 & 11.4 & cuOpt-IP & 137 & 0.58 & 0.76\\
         & $10^{-8}$ & 138 & 5.48 & HPR-QP & 126 & 12.69 & 2.31 & cuOpt-IP & 137 & 0.61 & 0.11\\
        \midrule
        \multirow{2}{*}{SuiteSparse (52)} & $10^{-4}$ & 52 & 1.57 & HPR-QP & 51 & 3.94 & 2.51 & Clarabel & 52 & 3.56 & 2.27\\
         & $10^{-8}$ & 50 & 5.21 & HPR-QP & 50 & 7.90 & 1.52 & Clarabel & 50 & 6.25 & 1.20\\
        \midrule
        \multirow{2}{*}{QPLIB (36)} & $10^{-4}$ & 35 & 4.86 & HPR-QP & 36 & 3.51 & 0.72 & Gurobi-IP & 35 & 5.33 & 1.10\\
         & $10^{-8}$ & 34 & 13.43 & HPR-QP & 33 & 17.45 & 1.30 & Clarabel & 30 & 18.22 & 1.36\\
        \midrule
        \multirow{2}{*}{Mittelmann (49)} & $10^{-4}$ & 48 & 17.55 & cuPDLPx & 48 & 16.80 & 0.96 & Gurobi-IP & 42 & 72.55 & 4.13\\
         & $10^{-8}$ & 45 & 57.72 & HPR-LP & 44 & 57.29 & 0.99 & Gurobi-IP & 38 & 141.01 & 2.44\\
        \midrule
        \multirow{2}{*}{MIPLIB (992)} & $10^{-4}$ & 986 & 2.07 & cuPDLPx & 979 & 3.15 & 1.52 & Gurobi-IP & 967 & 5.44 & 2.62\\
         & $10^{-8}$ & 971 & 4.66 & cuPDLPx & 949 & 7.28 & 1.56 & Gurobi-IP & 926 & 9.58 & 2.06\\
        \bottomrule
    \end{tabular*}
    \end{subtable}

    \vspace{1.2em}
    \begin{subtable}{\linewidth}
    \centering
    \caption{Portfolio backtest of~\Cref{sec:parametric-results} with~$n=3{,}000$ stocks.
    Total solve time in seconds of the 19,200 QPs for \osqplf{} with cold starts, with warm starts, and batched with warm starts, and for the fastest interior-point solver, with the ratio to the batched run.}
    \label{tab:summary-portfolio}
    \footnotesize
    \setlength{\tabcolsep}{3pt}
    \begin{tabular*}{\linewidth}{@{\extracolsep{\fill}}llrrrlrr@{}}
        \toprule
        & & \multicolumn{3}{c}{\osqplf{}, total time} & \multicolumn{3}{c}{Fastest interior-point solver}\\
        \cmidrule{3-5}\cmidrule{6-8}
        Portfolio backtest & $\varepsilon$ & Cold start & Warm start & Batched, warm start & Solver & Total time & Ratio\\
        \midrule
        \multirow{2}{*}{$k=50$ factors} & $10^{-4}$ & 1512.57 & 73.84 & 19.13 & cuOpt-IP & 1401.61 & 73.3\\
         & $10^{-8}$ & 1914.21 & 409.19 & 73.11 & Gurobi-IP & 2258.76 & 30.9\\
        \midrule
        \multirow{2}{*}{Full covariance} & $10^{-4}$ & 3317.90 & 207.88 & 148.20 & cuOpt-IP & 10571.94 & 71.3\\
         & $10^{-8}$ & 5857.66 & 1933.31 & 1516.41 & cuOpt-IP & 14407.90 & 9.50\\
        \bottomrule
    \end{tabular*}
    \end{subtable}
\end{table}
\paragraph{Summary of the results}
\Cref{tab:summary} compares \osqplf{} on each benchmark set with the fastest other first-order solver and with the fastest interior-point solver, \ie, the ones with the lowest shifted geometric mean.
Four observations stand out.

\begin{itemize}
    \item \emph{Robustness.}
        On the LP and QP sets, \osqplf{} solves the most instances in 9 of the 10 rows, in three of them tied with another solver.
        It is the only solver that solves all 138 Maros--M\'esz\'aros instances, which it does at both tolerances.
    \item \emph{First-order methods.}
        \osqplf{} is the fastest first-order solver in 7 of these 10 rows.
        It is~$2.3$ to~$11$ times faster than HPR-QP on Maros--M\'esz\'aros, $1.5$ to~$2.5$ times faster on SuiteSparse, and~$1.5$ to~$1.6$ times faster than cuPDLPx on MIPLIB.
        On the Mittelmann set it is within~$5\%$ of the fastest solver at both tolerances.
        It is clearly slower only on QPLIB at~$\varepsilon=10^{-4}$, by a factor of~$1.4$.
    \item \emph{Interior-point methods.}
        \osqplf{} is faster than the fastest interior-point solver in 8 of these 10 rows, by factors between~$1.1$ and~$4.1$.
        The exception is Maros--M\'esz\'aros, where the interior-point method of cuOpt is~$1.3$ times faster at~$\varepsilon=10^{-4}$ and~$9$ times faster at~$\varepsilon=10^{-8}$.
    \item \emph{Parametric problems.}
        On the portfolio backtest, where \osqplf{} solves the 20 QPs of each day as one batch and warm starts them, it is~$71$ to~$73$ times faster than the fastest interior-point solver at~$\varepsilon=10^{-4}$, and~$10$ to~$31$ times faster at~$\varepsilon=10^{-8}$.
\end{itemize}

\noindent Two things contribute to these results.
With a single variant, \osqplf{} is already the fastest first-order solver on Maros--M\'esz\'aros, SuiteSparse, and MIPLIB, by factors between~$1.2$ and~$2.4$.
Algorithm hedging adds a further factor of up to~$5.5$ on the instances where an ensemble is cheap, which are the smaller ones (\Cref{sec:hedging-results}).
On the largest instances the ensemble shrinks to a few variants, and \osqplf{} is up to~$2.2$ times slower than the fastest first-order solver, or~$2.5$ times after presolve.

\subsection{QP benchmarks}
\label{sec:qp-benchmarks}
We use three standard sets of QPs: the Maros--M\'esz\'aros set~\cite{marosRepositoryConvexQuadratic1999}, the Lasso and Huber problems built from least-squares matrices of the SuiteSparse collection~\cite{davisUniversityFloridaSparse2011,kolodziejSuiteSparseMatrixCollection2019}, and the convex instances of QPLIB~\cite{furiniQPLIBLibraryQuadratic2019}.
\Cref{tab:qp-size} reports their sizes.
We compare \osqplf{} with the GPU first-order solver HPR-QP~\cite{chenHPRQPDualHalpern2025}, the GPU interior-point method of cuOpt~\cite{nvidiaCuOptGPUacceleratedDecision2025}, and the CPU interior-point method of Gurobi~\cite{gurobioptimizationllcGurobiOptimizerReference2026}.
We also compare with the GPU versions of the interior-point solver Clarabel~\cite{goulartClarabelInteriorpointSolver2026,chenCuClarabelGPUAcceleration2026} and of the ADMM-based conic solver SCS~\cite{odonoghueConicOptimizationOperator2016,odonoghueOperatorSplittingHomogeneous2021}.
Among the other first-order QP solvers of~\Cref{subsec:relation}, we compare only with HPR-QP, the fastest in published comparisons: it is about three times faster than PDQP on Maros--M\'esz\'aros~\cite{chenHPRQPDualHalpern2025}, and faster than PDQP and PDHCG on Maros--M\'esz\'aros, QPLIB, and synthetic QPs~\cite{luCuPDQPxGPUAcceleratedRestarted2026}.
We do not include cuPDQPx~\cite{luCuPDQPxGPUAcceleratedRestarted2026}, whose code was not publicly available at the time of writing.

\begin{table}
    \centering
    \caption{Size groups of the QP benchmark sets by the number of nonzeros $\nnz(P) + \nnz(A)$ of the original instances, and the ranges of the number of variables $d$, the number of constraints $m$, and the number of nonzeros.}
    \label{tab:qp-size}
    \footnotesize
    \setlength{\tabcolsep}{3pt}
    \begin{adjustbox}{max width=\linewidth}
    \begin{tabular}{lccccrrr}
        \toprule
        \multirow{2}{*}{Set} & Small & Medium & Large & \multirow{2}{*}{Total} & \multicolumn{3}{c}{Range}\\
        \cmidrule(lr){6-8}
         & $\nnz<10^4$ & $10^4 \leq \nnz \leq 10^5$ & $\nnz>10^5$ & & \multicolumn{1}{c}{$d$} & \multicolumn{1}{c}{$m$} & \multicolumn{1}{c}{$\nnz(P) + \nnz(A)$}\\
        \midrule
        Maros--M\'esz\'aros & 83 & 47 & 8 & 138 & 2--93{,}263 & 1--186{,}531 & 4--2{,}257{,}500\\
        SuiteSparse & 22 & 8 & 22 & 52 & 60--6{,}043{,}555 & 32--3{,}943{,}590 & 354--15{,}702{,}708\\
        QPLIB & 2 & 21 & 13 & 36 & 150--1{,}009{,}306 & 0--1{,}001{,}000 & 4{,}000--9{,}653{,}690\\
        \bottomrule
    \end{tabular}
    \end{adjustbox}
\end{table}

\paragraph{Maros--M\'esz\'aros}
The Maros--M\'esz\'aros set~\cite{marosRepositoryConvexQuadratic1999} has 138 QPs, most of them small (\Cref{tab:maros-sgm}).
\osqplf{} solves all of them at both tolerances, while HPR-QP solves 130 and 126, and SCS 127 and 119.
At~$\varepsilon=10^{-4}$ its time of~$0.76$ seconds is~$11$ times lower than that of HPR-QP and more than~$16$ times lower than that of SCS.
It matches Gurobi ($0.77$ seconds) and is within a factor of~$1.3$ of cuOpt, the fastest solver.
This set is where hedging matters most: with a single variant, \osqplf{} is~$5.5$ times slower (\Cref{sec:hedging-results}), and still two times faster than HPR-QP.
At~$\varepsilon=10^{-8}$, \osqplf{} remains the fastest first-order solver, $2.3$ times faster than HPR-QP, but the interior-point methods are ahead.
The interior-point method of cuOpt takes~$0.61$ seconds against~$5.48$ for \osqplf{}, and the difference comes from the medium instances, which take~$1.6$ and~$17$ seconds.
Presolve changes these times by less than~$15\%$ for \osqplf{} and does not change the ranking among the first-order solvers.

\paragraph{SuiteSparse matrix collection}
We take 26 least-squares problems from the SuiteSparse collection~\cite{davisUniversityFloridaSparse2011,kolodziejSuiteSparseMatrixCollection2019} and build a Lasso and a Huber fitting problem from each, as in~\cite{stellatoOSQPOperatorSplitting2020}, which gives 52 QPs (\Cref{tab:suite-sgm}).
On this set \osqplf{} is the fastest of all solvers at both tolerances.
At~$\varepsilon=10^{-4}$ it solves all 52 instances in~$1.57$ seconds, which is~$2.5$ times faster than HPR-QP, $2.3$ times faster than Clarabel, the fastest interior-point solver, and~$8.5$ times faster than Gurobi.
The 22 large instances, with up to~$1.6\cdot10^7$ nonzeros, account for this: \osqplf{} takes~$3.96$ seconds on them, against~$11.44$ for HPR-QP, $10.37$ for Clarabel, and~$63.41$ for Gurobi, which solves 19 of them.
On the small instances the interior-point methods are faster, but every solver except SCS and cuOpt takes less than~$0.2$ seconds there.
At~$\varepsilon=10^{-8}$, \osqplf{} is~$1.5$ times faster than HPR-QP, $1.2$ times faster than Clarabel, and~$6.5$ times faster than Gurobi.
It is also the fastest solver on the presolved instances, by factors of~$1.5$ and~$1.4$ over Clarabel.
Over the whole set, hedging changes the time by at most~$3\%$ (\Cref{sec:hedging-results}), so the advantage comes from the method with a single variant.

\paragraph{QPLIB}
We take the 36 instances of QPLIB~\cite{furiniQPLIBLibraryQuadratic2019} with a convex quadratic objective and linear or bound constraints (\Cref{tab:qplib-sgm}).
Of these, 17 have binary variables, which we relax to the interval~$[0, 1]$.
This is the only set on which another first-order solver is clearly faster.
At~$\varepsilon=10^{-4}$, HPR-QP solves all 36 instances in~$3.51$ seconds, and \osqplf{} solves 35 in~$4.86$ seconds, ahead of Gurobi ($5.33$ seconds).
The two size groups behave in opposite ways.
On the 21 medium instances, \osqplf{} is~$3.1$ times faster than HPR-QP and~$7.8$ times faster than Gurobi.
On the 13 large instances it is~$2.2$ times slower than HPR-QP and solves one instance fewer.
At~$\varepsilon=10^{-8}$ the order reverses: \osqplf{} is the fastest of all solvers and solves the most instances, 34.
It is~$1.3$ times faster than HPR-QP and~$1.4$ times faster than Clarabel and Gurobi.
On the presolved instances, \osqplf{} is~$1.4$ times faster than HPR-QP at both tolerances, while an interior-point solver is the fastest: Gurobi at~$\varepsilon=10^{-4}$ and Clarabel at~$\varepsilon=10^{-8}$, each by a factor of~$1.3$.

\subsection{LP benchmarks}
\label{sec:lp-benchmarks}
We use two standard sets of LPs: Mittelmann's LP benchmark~\cite{mittelmannBenchmarksOptimizationSoftware2026} and the LP relaxations of the MIPLIB 2017 collection~\cite{gleixnerMIPLIB2017Datadriven2021}.
\Cref{tab:lp-size} reports their sizes.
We compare \osqplf{} with the GPU first-order solvers cuPDLPx~\cite{luCuPDLPxFurtherEnhanced2025} and HPR-LP~\cite{chenHPRLPImplementationHPR2025}, the PDLP implementation of cuOpt~\cite{nvidiaCuOptGPUacceleratedDecision2025}, and the PDHG method of Gurobi~\cite{gurobioptimizationllcGurobiOptimizerReference2026}.
We also report the CPU interior-point method of Gurobi and the GPU interior-point method of Clarabel.

\begin{table}
    \centering
    \caption{Size groups of the LP benchmark sets by the number of nonzeros $\nnz(A)$ of the original instances, and the ranges of the number of variables $d$, the number of constraints $m$, and the number of nonzeros.}
    \label{tab:lp-size}
    \footnotesize
    \setlength{\tabcolsep}{3pt}
    \begin{adjustbox}{max width=\linewidth}
    \begin{tabular}{lccccrrr}
        \toprule
        \multirow{2}{*}{Set} & Small & Medium & Large & \multirow{2}{*}{Total} & \multicolumn{3}{c}{Range}\\
        \cmidrule(lr){6-8}
         & $\nnz<10^6$ & $10^6 \leq \nnz \leq 10^7$ & $\nnz>10^7$ & & \multicolumn{1}{c}{$d$} & \multicolumn{1}{c}{$m$} & \multicolumn{1}{c}{$\nnz(A)$}\\
        \midrule
        Mittelmann's LP & 11 & 27 & 11 & 49 & 6{,}624--38{,}868{,}107 & 4{,}284--30{,}000{,}600 & 94{,}950--93{,}001{,}800\\
        MIPLIB 2017 & 845 & 123 & 24 & 992 & 3--38{,}868{,}107 & 1--19{,}912{,}111 & 3--256{,}963{,}402\\
        \bottomrule
    \end{tabular}
    \end{adjustbox}
\end{table}

\paragraph{Mittelmann's LP}
On the 49 instances of this set (\Cref{tab:Mittelmann-LP-sgm}), the three fastest first-order solvers are within~$14\%$ of one another.
At~$\varepsilon=10^{-4}$, cuPDLPx takes~$16.80$ seconds, \osqplf{}~$17.55$, and HPR-LP~$18.59$, and each solves 48 instances.
At~$\varepsilon=10^{-8}$, HPR-LP takes~$57.29$ seconds, \osqplf{}~$57.72$, and cuPDLPx~$65.12$, and \osqplf{} solves the most instances, 45.
Hedging matters little here, since 43 of the 49 instances are too large for a batch and run a single variant (\Cref{sec:hedging-results}): it gains a factor of~$1.13$ at~$\varepsilon=10^{-4}$ and none at~$10^{-8}$.
Hence \osqplf{} matches the state of the art on this set with one variant.
By size, it is the fastest solver on the 27 medium instances, $1.5$ and~$1.4$ times faster than cuPDLPx, and about~$2$ times slower than HPR-LP on the 11 large ones.
On the small instances at~$\varepsilon=10^{-4}$ it fails on one of 11, which dominates its mean for that group.
All these first-order solvers are well ahead of the interior-point methods: Gurobi is~$4.1$ and~$2.4$ times slower than \osqplf{} and solves 42 and 38 instances.
Presolve makes \osqplf{}~$2.3$ and~$2.1$ times faster.
On the presolved instances, HPR-LP is the fastest at~$\varepsilon=10^{-4}$ ($6.97$ seconds against~$7.57$ for \osqplf{}), and \osqplf{} at~$\varepsilon=10^{-8}$ ($27.57$ against~$29.88$).

\paragraph{MIPLIB}
The MIPLIB 2017 collection has 1{,}065 instances.
We use the 992 that remain after excluding the 42 with indicator constraints, the 27 that the collection lists as unbounded or of unknown status, and 4 infeasible ones whose LP relaxation Gurobi does not solve to optimality.
On these 992 LP relaxations (\Cref{tab:MIPLIB-sgm}), \osqplf{} is the fastest solver and solves the most instances at both tolerances.
At~$\varepsilon=10^{-4}$ it solves 986 instances in~$2.07$ seconds.
It is~$1.5$ times faster than cuPDLPx, $2.2$ times faster than HPR-LP and the PDLP of cuOpt, and~$2.6$ times faster than the interior-point method of Gurobi.
At~$\varepsilon=10^{-8}$ it solves 971 instances in~$4.66$ seconds, $1.6$ times faster than cuPDLPx, $1.8$ times faster than HPR-LP, and~$2.1$ times faster than Gurobi.
It solves 13 more instances than HPR-LP and 22 more than cuPDLPx.
The advantage is largest on the 845 small instances, where \osqplf{} is~$2.3$ and~$2.1$ times faster than cuPDLPx.
Hedging accounts for factors of~$1.8$ and~$1.3$ of this (\Cref{sec:hedging-results}), and the single variant for the rest.
On the 24 large instances it is~$1.3$ times slower than HPR-LP at~$\varepsilon=10^{-4}$ and~$1.9$ times slower than cuPDLPx at~$\varepsilon=10^{-8}$.
The presolved instances show the same picture, with \osqplf{}~$1.5$ and~$1.7$ times faster than cuPDLPx.

\subsection{Batching parametric instances}
\label{sec:parametric-results}
We now evaluate the batching of instances that share~$P$ and~$A$ (\Cref{subsec:batch-mechanism}) on a four-year portfolio backtest.
On each of 960 trading days, we solve the mean--variance portfolio problem of~\Cref{app:portfolio} for 20 values of the risk-aversion parameter, which gives 19,200 QPs.
The 20 QPs of a day differ only in the vector~$q$, so they form one batch.
We run \osqplf{} in three ways.
\begin{itemize}
    \item \emph{Cold start.} It solves the 19,200 QPs one after the other.
    \item \emph{Warm start.} It also starts each QP from the solution of the previous day for the same risk-aversion parameter.
    \item \emph{Batched, warm start.} It also solves the 20 QPs of a day as one batch.
\end{itemize}
We compare with Gurobi, Clarabel, and the interior-point method of cuOpt, which solve the QPs one after the other.

\begin{figure}
    \centering
    \includegraphics[width=\linewidth]{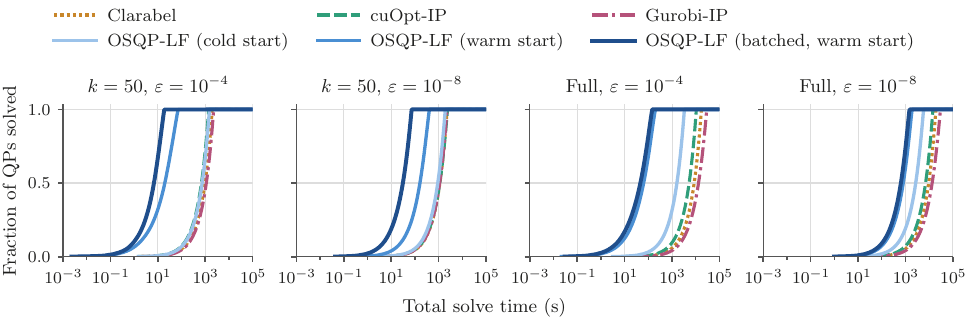}
    \caption{Fraction of the 19,200 QPs of the portfolio backtest with $n = 3{,}000$ stocks solved against the total solve time, with a risk model of $k = 50$ factors (left two panels) and with the full covariance matrix (right two panels), for tolerances $\varepsilon = 10^{-4}$ and $\varepsilon = 10^{-8}$.
    Each day has 20 QPs, one for each value of the risk-aversion parameter, and the batched run solves a full day as one batch.}
    \label{fig:portfolio-top3000-days}
\end{figure}

\Cref{fig:portfolio-top3000-days} shows the fraction of the QPs solved within a given total time for $n = 3{,}000$ stocks, with a risk model of $k = 50$ factors and with the full covariance matrix, and~\Cref{tab:portfolio-time-3000} reports the total solve times.
At~$\varepsilon=10^{-4}$ and with 50 factors, the cold-started solver takes about as long as the fastest interior-point solver, roughly 1{,}500 seconds.
Warm starting reduces this to 74 seconds, and batching to 19 seconds.
The batched solver is then~$73$ times faster than cuOpt and more than~$100$ times faster than Clarabel and Gurobi.
With the full covariance matrix, warm starting and batching make the solver~$22$ times faster than the cold-started run, and~$71$ to~$196$ times faster than the interior-point solvers.
Here~$\Sigma$ is a dense~$3{,}000\times3{,}000$ matrix, so the products with~$P$ dominate each iteration, and their cost grows with the number of columns (\Cref{fig:batch-cost}).
Batching still cuts the time of the warm-started run by 30\%, from 208 to 148 seconds.
On the logarithmic time axis of~\Cref{fig:portfolio-top3000-days}, this gain places the two curves close together in the right two panels.
At~$\varepsilon=10^{-8}$ the gains are smaller but still large: about~$30$ times with 50 factors, and~$10$ to~$20$ times with the full covariance matrix.
The results for the other numbers of stocks and factors are similar, and we report them in~\Cref{app:detailed-results}.

\subsection{Effect of batching and hedging}
\label{sec:hedging-results}
We now isolate the two ingredients of~\Cref{sec:gpu-batching}: the cost of running a batch, and the benefit of running different variants in it.

\begin{figure}[!ht]
    \centering
    \includegraphics[width=0.48\linewidth]{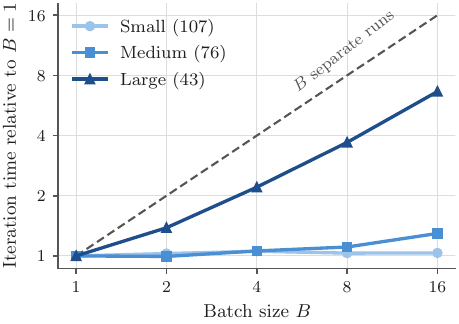}
    \caption{Cost of batching. Median time of one iteration at batch size~$B$ divided by the time at batch size~$1$, on the Maros--M\'esz\'aros, SuiteSparse, and QPLIB instances, by size group, with the number of instances in parentheses.
    The dashed line is the cost of running the~$B$ variants one after the other.}
    \label{fig:batch-cost}
\end{figure}

\paragraph{Cost of batching}
We measure the time of one iteration of~\osqplf{} with batch sizes~$B \in \{1, 2, 4, 8, 16\}$ on the 226 instances of the Maros--M\'esz\'aros, SuiteSparse, and QPLIB benchmark sets.
\Cref{fig:batch-cost} shows the median of the ratio between the time at batch size~$B$ and the time at batch size~$1$, for each size group.
On the small instances, an iteration with 16 variants costs about as much as an iteration with one.
On the medium instances it costs about $1.3$ times as much.
On the large instances the cost grows with the batch size, but an iteration with 16 variants still takes less than half the time of 16 separate iterations.
This is why the batch size rule of~\Cref{subsec:hedging-implementation} gives large instances a small batch.
The LP benchmark sets behave in the same way, with a median ratio of~$1.8$ at~$B=16$ on their small instances and~$8$ on the medium and large ones.
With the sparse products of cuSPARSE alone, the same ratio on Maros--M\'esz\'aros is~$3.4$ instead of~$1.1$, which is the reason for the custom kernels of~\Cref{app:algebra-engine}.

\begin{table}[!ht]
    \centering
    \footnotesize
    \setlength{\tabcolsep}{4.5pt}
    \caption{Winning variants for tolerance~$\varepsilon = 10^{-4}$ on the original instances: the number of instances, how many run more than one variant ($E > 1$), how many the first variant of the order wins (First), and how many distinct variants win at least once.}
    \label{tab:hedging-results}
    \begin{tabular}{lrrrr}
        \toprule
         & \multicolumn{3}{c}{Instances} & Variants\\
        \cmidrule(lr){2-4} \cmidrule(lr){5-5}
        Benchmark set & All & $E > 1$ & First & Winning\\
        \midrule
        MIPLIB & 992 & 780 & 304 & 16\\
        Mittelmann & 49 & 6 & 44 & 3\\
        Maros--M\'esz\'aros & 138 & 136 & 6 & 16\\
        SuiteSparse & 52 & 34 & 20 & 11\\
        QPLIB & 36 & 32 & 8 & 11\\
        \bottomrule
    \end{tabular}
\end{table}

\paragraph{Winning variants}
For each instance, we record the variant of the ensemble that satisfies the termination criteria first.
\Cref{tab:hedging-results} summarizes the results for tolerance~$\varepsilon = 10^{-4}$ on the original instances.
The first variant of the order is the fastest on 304 of the 992 MIPLIB instances, but it runs alone on 212 of them.
On the 780 instances that run more than one variant, it is the fastest fewer than 100 times.
On Maros--M\'esz\'aros it is the fastest on 6 of the 138 instances.
All 16 variants are the fastest on some instance of MIPLIB and of Maros--M\'esz\'aros.
The Mittelmann instances behave differently: 43 of the 49 are too large for a batch and run the first variant alone.

\begin{figure}[!ht]
    \centering
    \includegraphics[width=\linewidth]{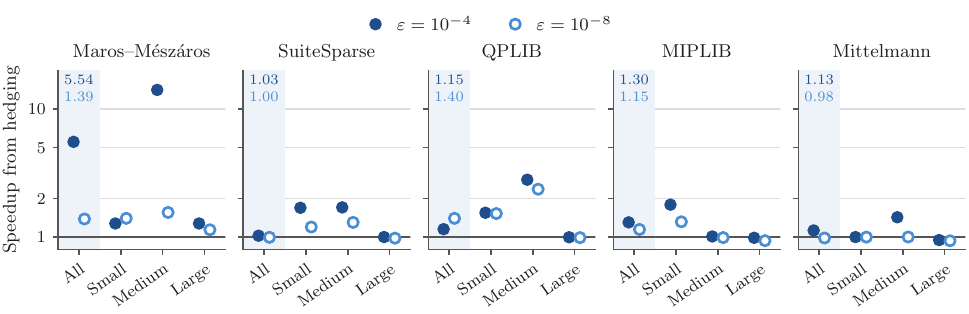}
    \caption{Hedging against a single variant on the original instances. Shifted geometric mean of the solve times of the single-variant run divided by that of the hedged run, for all instances and by size group, at both tolerances.
    Values above~$1$ favor hedging.
    The numbers in the shaded band are the values for all instances, at~$\varepsilon=10^{-4}$ (top) and~$\varepsilon=10^{-8}$ (bottom).}
    \label{fig:hedging-speedup}
\end{figure}

\paragraph{Hedging against a single variant}
We rerun every benchmark set on the original instances with the ensemble size fixed to~$E=1$, so that only the first variant of the order runs, and compare with the runs of~\Cref{sec:qp-benchmarks,sec:lp-benchmarks}.
\Cref{fig:hedging-speedup} shows the ratio between the shifted geometric means of the two runs.
Hedging never solves fewer instances than the single variant, and over the ten combinations of benchmark set and tolerance it solves 18 more, 11 of them on Maros--M\'esz\'aros.
The gain concentrates on the small and medium instances, where the batch is cheap.
There, hedging is between~$1.2$ and~$2.8$ times faster on the QP sets, and~$14$ times faster on the medium Maros--M\'esz\'aros instances at~$\varepsilon=10^{-4}$, where the single variant reaches the time limit on three instances.
On the large instances the batch size rule leaves little room for more than one variant, so the two runs take about the same time, except on Maros--M\'esz\'aros, where hedging still gains a factor of~$1.1$ to~$1.3$.
The small overhead that remains makes hedging at most~$6\%$ slower in those groups, and~$2\%$ slower on the full Mittelmann set at~$\varepsilon=10^{-8}$.

\section{Conclusion}
\label{sec:conclusion}
We presented a GPU-accelerated first-order solver for LPs and convex QPs that runs an ensemble of algorithm variants on the same problem as a single batch, and stops as soon as one of them converges.
This \emph{algorithm hedging} protects the solver against a poor choice of variant, since no single variant is best on every problem.
The solver builds on linear-system-free ADMM, in which one preconditioned projected gradient step replaces the linear system solve of ADMM, so that each iteration needs only matrix--vector products and simple projections.
The variants share every product with the problem matrices and differ in acceleration scheme, step-size rule, and relaxation parameter.
On a GPU, the ensemble therefore runs at a fraction of the cost of separate runs.
With restarted Halpern iteration, the method converges linearly on problems with polyhedral constraints, and its iterates certify primal and dual infeasibility.

Our open-source implementation, \osqplf{}, is a new GPU mode of the OSQP solver.
On five standard LP and QP benchmark sets, it is the fastest first-order solver in most cases, and it solves the most instances in nearly all of them.
It is already the fastest first-order solver on three of the five sets with a single variant, and hedging makes it up to five times faster on top of that, without ever solving fewer instances.
The same batching solves many instances of a parametric problem at once.
In a portfolio backtest, together with warm starting, it reduces the total solve time by one to two orders of magnitude compared with interior-point solvers.

\section*{Acknowledgments}
Liubomir Maxim Baicev, Jisun Park, and Bartolomeo Stellato are supported by the ONR Young Investigator Program Award N00014-25-1-2147.
Liubomir Maxim Baicev and Bartolomeo Stellato are also supported by the Alfred P. Sloan Foundation through a Sloan Research Fellowship (grant FG-2026-79479).
Bartolomeo Stellato is also supported by the NSF CAREER Award ECCS-2239771.
The experiments ran on a computing node that Princeton Research Computing, a consortium led by the Princeton Institute for Computational Science and Engineering and the Office of Information Technology's Research Computing, hosts and maintains at Princeton University.

\section*{Use of AI tools}
We used AI assistants (Claude, by Anthropic) to develop parts of the solver code, including the translation of early prototypes to CUDA, to edit the text, to check the proofs and formalize them in Lean 4, to write the scripts that generate the tables and figures from the benchmark results, and to verify the bibliography.
The authors reviewed all of this work and take full responsibility for the content of the paper.

\appendix
\crefalias{section}{appendix}
\crefalias{subsection}{subappendix}

\section{Derivation from FLiP-ADMM}
\label{app:flip}
\Cref{subsec:derivation} derives linear-system-free ADMM from a proximal ADMM.
This appendix derives it from the function-linearized proximal ADMM (FLiP-ADMM) of Ryu and Yin~\cite[Chapter~8]{ryuLargeScaleConvexOptimization2022} in two ways, and uses the second one to prove~\Cref{cor:lsf-averagedness}($b$).

We consider the more general problem
\begin{equation*}
    \begin{array}{ll}
        \textnormal{minimize} & f(x) \\
        \textnormal{subject to} & Ax = z, \quad x \in \cX, \quad z \in \cC,
    \end{array}
\end{equation*}
where~$x \in \reals^d$ and~$z \in \reals^m$ are the optimization variables and~$A$, $\cX$, and~$\cC$ are as in~\eqref{prob:qp}.
The objective splits into~$f = f_1 + f_2$, where~$f_1$ is closed, convex, and proper and~$f_2$ is convex and differentiable.
Problem~\eqref{prob:qp} is the case~$f(x) = (1/2) x^T P x + q^T x$.
The function-linearized augmented Lagrangian is
\begin{equation*}
    \lagrangian_k (x, z, y) =
        f_1(x)
        + f_2(x^k) + \nabla f_2 (x^k)^T (x - x^k)
        + y^T (Ax - z)
        + (\rho/2) \norm{Ax - z}^2,
\end{equation*}
which linearizes~$f_2$ at the iterate~$x^k$, with the multiplier~$y$ and the step size~$\rho$ of~\Cref{subsec:derivation}.
FLiP-ADMM is the iteration
\begin{equation}\label{algo:flip-admm-qp}
    \begin{array}{rcl}
        \xkk &=& \underset{x \in \cX}{\argmin} \left(
            \lagrangian_k (x, z^k, y^k)
            + (1/2) \norm[W]{x - x^k}^2
        \right) \\
        \zkk &=& \underset{z \in \cC}{\argmin} \, \lagrangian_k (x^{k+1}, z, y^k) 
        = \proj_\cC\!\left(A\xkk + \rho^{-1}\yk\right) \\
        \ykk &=& \yk + \rho (A \xkk - \zkk),
    \end{array}
    \tag{FLiP-ADMM}
\end{equation}
where~$W \succeq 0$ is a proximal matrix.
With~$f_2 = 0$ and the objective of~\eqref{prob:qp}, it is the semi-proximal ADMM~\eqref{algo:prox-admm}.
The next proposition chooses~$W$ so that the~$x$-update is explicit.

\begin{proposition}
    \label{prop:flip-explicit}
    Suppose that~$f_1$ is a convex quadratic function, and let~$R \succ 0$ be such that the proximal matrix~$W = R - (\nabla^2 f_1 + \rho A^T A)$ is positive semidefinite.
    Then the~$x$-update of~\eqref{algo:flip-admm-qp} is the projected gradient step
    \begin{equation*}
        \xkk = \proj_\cX^R \! \left( \xk - R^{-1} \nabla_x \lagrangian_k (\xk, \zk, \yk) \right),
    \end{equation*}
    where~$\proj_\cX^R$ is the projection onto~$\cX$ in the norm~$\|\cdot\|_R$ (\Cref{subsec:derivation}).
    If~$R$ is diagonal, then~$\proj_\cX^R = \proj_\cX$ is a componentwise clipping.
\end{proposition}

\begin{proof}
    Since~$f_1$ is quadratic, its Hessian~$\nabla^2 f_1$ is constant, and~$\nabla f_1 (x^{k+1}) = \nabla f_1 (x^k) + \nabla^2 f_1 (x^{k+1} - x^k)$.
    Substituting this identity and~$A x^{k+1} - z^k = A (x^{k+1} - x^k) + (A x^k - z^k)$ into the optimality condition of the~$x$-update gives
    \begin{equation*}
        0 \in
        \cN_\cX (x^{k+1})
        + \left( \nabla^2 f_1 + \rho A^T A + W \right) (x^{k+1} - x^k)
        + \nabla_x \lagrangian_k (x^k, z^k, y^k).
    \end{equation*}
    Since~$\nabla^2 f_1 + \rho A^T A + W = R$, this inclusion reads~$0 \in \cN_\cX (\xkk) + R (\xkk - \hat x^k)$ with~$\hat x^k = \xk - R^{-1} \nabla_x \lagrangian_k (\xk, \zk, \yk)$.
    It is the optimality condition for minimizing~$(1/2) \|x - \hat x^k\|_R^2$ over~$x \in \cX$, whose unique solution is~$\proj_\cX^R (\hat x^k)$ because~$\cX$ is closed and convex and~$R \succ 0$.
    For a diagonal~$R$, this problem separates over the entries of~$x$, and each entry of the solution is the clipping of the corresponding entry of~$\hat x^k$ to its interval.
\end{proof}

\begin{figure}
    \centering
    \begin{tikzpicture}[
    x=1cm, y=1cm, >=latex,
    method/.style={draw=black!45, rounded corners=2pt, align=center,
        minimum height=.75cm, inner xsep=5pt, inner ysep=3pt, font=\footnotesize},
    route/.style={method, minimum width=4.2cm},
    focus/.style={method, draw=black, line width=.8pt},
    edge/.style={->, draw=black!65, line width=.6pt},
    repr/.style={edge, dashed},
    lab/.style={font=\scriptsize, align=center, inner sep=2pt},
]
\node[method] (flip) at (0, 0) {FLiP-ADMM};
\node[route] (padmm) at (4.0, 1.05)
    {Proximal ADMM\\{\scriptsize proximal point on~$\opA$}};
\node[route] (cv) at (4.0, -1.05)
    {Condat--V\~u\\{\scriptsize forward--backward on~$\opA + \opH$}};
\node[focus] (lsf) at (8.9, 0) {Linear-system-free ADMM};
\node[method] (pdhg) at (13.1, 0) {PDHG};
\draw[edge] (flip) -- node[lab, above left] {route~$(i)$: $f_2 = 0$} (padmm);
\draw[repr] (flip) -- node[lab, below left] {route~$(ii)$: $f_1 = 0$;\\eliminate~$z$} (cv);
\draw[edge] (padmm) -- node[lab, above right] {$W = R - P - \rho A^T A$} (lsf);
\draw[edge] (cv) -- node[lab, below right] {$W = R - \rho A^T A$} (lsf);
\draw[edge] (lsf) -- node[lab, above] {LP} node[lab, below] {$P = 0$} (pdhg);
\end{tikzpicture}
     \caption{The two routes of~\Cref{prop:two-routes}.
    Solid arrows point from a method to a special case.
    The dashed arrow also eliminates~$z$, as in~\Cref{thm:lsf-osqp-condat-vu}($a$), since the Condat--V\~u method is stated in the variables~$(x,y)$.
    Route~$(i)$ keeps~$f$ in the proximal step, so~$\opTCV$ is a proximal point operator.
    Route~$(ii)$ linearizes~$f$, so~$\opTCV$ is a Condat--V\~u operator.
    Both routes give the same step~$\opF$.}
    \label{fig:two-routes}
\end{figure}
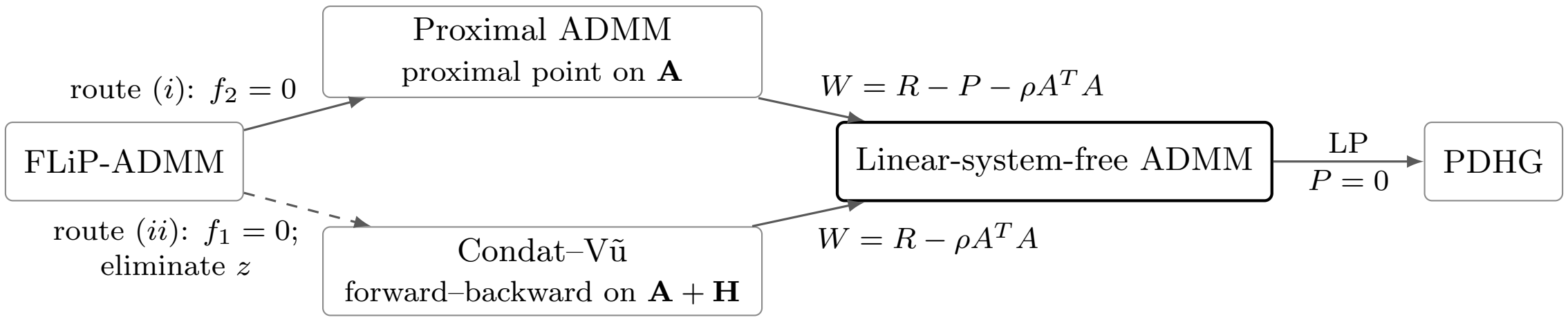

Two routes lead from FLiP-ADMM to the step~$\opF$ of linear-system-free ADMM (\Cref{fig:two-routes}).
Route~$(i)$ is the derivation of~\Cref{subsec:derivation}, in which~$\opF$ is a step of proximal ADMM.
In route~$(ii)$, the primal--dual form~\eqref{eq:lsf-primal-dual} of~$\opF$ is a step of the preconditioned Condat--V\~u method.
Each route gives a metric in which~$\opF$ is averaged.

\begin{proposition}
    \label{prop:two-routes}
    Suppose that~$f$ is the objective of~\eqref{prob:qp}, and let~$R \succ 0$ be such that both proximal matrices~$W$ below are positive semidefinite.
    Then the following two routes give the same~$x$-update of~\eqref{algo:flip-admm-qp}, which is the projected gradient step of~\Cref{prop:flip-explicit}.
    \begin{enumerate}
        \item[$(i)$]
            \emph{Majorize} (proximal ADMM): $f_1=f$ and~$f_2=0$, with~$W = R - (P + \rho A^T A)$.
        \item[$(ii)$]
            \emph{Linearize} (preconditioned Condat--V\~u): $f_1=0$ and~$f_2=f$, with~$W = R-\rho A^T A$.
    \end{enumerate}
\end{proposition}

\begin{proof}
    Routes~$(i)$ and~$(ii)$ are the cases~$\nabla^2 f_1 = P$ and~$\nabla^2 f_1 = 0$ of~\Cref{prop:flip-explicit}.
    Since~$f_1 + f_2 = f$ in both, the gradient~$\nabla_x \lagrangian_k (\xk, \zk, \yk)$ is the same.
\end{proof}

Route~$(ii)$ evaluates the gradient~$Px+q$ of the objective at the current iterate.
In the inclusion of~\Cref{thm:lsf-osqp-condat-vu}($a$), this moves the objective out of the operator~$\opA$ and of the metric~$\opM$.
Let
\begin{equation*}
    \opA_0(x,y)=\begin{bmatrix}
        \cN_\cX(x)+A^Ty \\
        \partial S_\cC(y)-Ax
    \end{bmatrix},
    \qquad
    \opH(x,y)=\begin{bmatrix}
        Px+q \\
        0
    \end{bmatrix},
    \qquad
    \opM_0=\begin{bmatrix}
        R & A^T \\
        A & \rho^{-1}\identity
    \end{bmatrix},
\end{equation*}
so that~$\opA=\opA_0+\opH$ and~$\opM=\opM_0-\diag(P,0)$.
For all~$u=(x,y)$ and~$u^+=(x^+,y^+)$, the terms~$Px^+$ cancel in~$\opM(u^+-u)+\opA u^+$, which gives
\begin{equation*}
    \opM(u^+-u)+\opA u^+=\opM_0(u^+-u)+\opA_0u^++\opH u.
\end{equation*}
Hence the inclusion~$(\opM+\opA)u^{k+1}\ni\opM u^k$ of~\Cref{thm:lsf-osqp-condat-vu}($a$) is the same as~$(\opM_0+\opA_0)u^{k+1}\ni(\opM_0-\opH)u^k$, and
\begin{equation}
    \label{eq:forward-backward-form}
    \opTCV=(\opM_0+\opA_0)^{-1}(\opM_0-\opH).
\end{equation}
This is the operator of the preconditioned Condat--V\~u method, a forward--backward step on~$\opA_0+\opH$ in the metric~$\opM_0$.
With~$W_0=R-\rho A^TA$ and~$\opW_0=\diag(W_0,\rho\identity)$, the block multiplication of~\Cref{thm:lsf-osqp-condat-vu}($c$) gives~$\Phi^T\opW_0\Phi=\opM_0$.
Hence~$\opM_0\succ0$ if and only if~$W_0\succ0$, and~$\Phi$ is also an isometry between the~$\opM_0$-norm and the~$\opW_0$-norm.

\begin{proof}[Proof of~\Cref{cor:lsf-averagedness}($b$)]
    By~\eqref{eq:forward-backward-form}, $\opTCV$ is a variable metric forward--backward operator with metric~$\opM_0\succ0$~\cite[\S2.8]{ryuLargeScaleConvexOptimization2022}.
    The backward step~$(\opM_0+\opA_0)^{-1}\opM_0$ is~$\opTCV$ for~$P=0$ and~$q=0$, so it is everywhere defined, single-valued, and firmly nonexpansive in the~$\opM_0$-norm by~\Cref{thm:lsf-osqp-condat-vu}($a$) and~\Cref{cor:lsf-averagedness}($a$).
    Since~$(\opM_0^{-1})_{11}=W_0^{-1}$, the gradient~$\opM_0^{-1}\opH$ of the convex quadratic~$(x,y)\mapsto(1/2)x^TPx+q^Tx$ in this metric has Lipschitz constant~$L$.
    For~$0<L<2$, this gradient is~$1/L$-cocoercive, so the forward--backward case ($B=0$) of~\cite[Theorem~28]{ryuLargeScaleConvexOptimization2022} with unit step size gives that~$\opTCV$ is~$\theta$-averaged in the~$\opM_0$-norm with~$\theta=2/(4-L)$.
    If~$L=0$, the forward step is a translation, so~$\opTCV$ is firmly nonexpansive.
    In both cases, the isometry~$\Phi$ and~$\opF=\Phi\circ\opTCV\circ\Phi^{-1}$ transfer the bound to~$\opF$ in the~$\opW_0$-norm, and the relaxation bounds follow from the definition of averagedness.
\end{proof}

\section{Termination criteria compared with OSQP}
\label{app:relation}
The termination criteria~\eqref{eq:termination} of~\Cref{subsec:termination-criteria} differ from those of OSQP~\cite{stellatoOSQPOperatorSplitting2020} in two respects.

First, they depend only on the point~$(x,y,w)$ and the problem data.
The primal residual of OSQP and its tolerance involve the iterate~$z$, whereas~$r_{\mathrm{prim}}$ and~$\varepsilon_\mathrm{prim}$ depend only on~$x$.

Second, they do not depend on how the input encodes the variable bounds.
Since OSQP has~$\cX=\reals^d$, the variable bounds enter as rows of the constraint matrix, either as a full identity matrix or as the rows of one for the bounded variables, and its criteria change with this choice.
The criteria~\eqref{eq:termination} take the same values whether a bound enters as a row of~$A$, with its multiplier in~$y$, or through~$\cX$, with its multiplier in~$w$.
Indeed, a bound row moves an entry of~$r_{\mathrm{prim}}$ between its two blocks and adds an entry of~$x$ to~$Ax$, which~$\varepsilon_\mathrm{prim}$ already covers through~$\|x\|_\infty$.
The other residuals and tolerances contain the multipliers only through the sums~$A^Ty+w$ and~$S_\cC(y)+S_\cX(w)$, which do not change when a multiplier moves between~$y$ and~$w$.

\section{KKT residuals at the projected point}
\label{app:projected-residual}
This appendix bounds the KKT residuals of~\eqref{prob:qp} at the point~$(\xtkk,z^k,y^k,w^{k+1})$ by the fixed-point residual.
\Cref{subsec:sharpness-free} uses the bound to transfer the linear rate of restarted Halpern iteration to the KKT residuals, and~\Cref{subsec:termination-criteria} evaluates the termination criteria at this point.
By~\Cref{subsec:kkt-residuals}, the point satisfies every inclusion of the KKT conditions~\eqref{eq:qp-kkt}, so only the two equations can be violated.

\begin{lemma}[KKT residuals at the projected point]
    \label{lem:projected-residual}
    The residuals at~$(\xtkk,z^k,y^k,w^{k+1})$ satisfy
    \begin{equation}
        \label{eq:kkt-residual-projected}
        \bar r^k
        =\begin{bmatrix}
            P\xtkk+q+A^Ty^k+w^{k+1}\\
            A\xtkk-z^k
        \end{bmatrix}
        =-\begin{bmatrix}
            R-P-\rho A^TA & \rho A^T\\
            0 & -\identity
        \end{bmatrix}
        (\opF-\identity)s^k.
    \end{equation}
    If~$W=R-P-\rho A^TA\succeq0$, then for every~$\alpha\in(0,2]$,
    \begin{equation}
        \label{eq:residual-transfer-bound}
        \norm{\bar r^k}
        \leq\kappa_\mathrm{res}\norm[\opW]{(\opF-\identity)s^k}
        =\frac{\kappa_\mathrm{res}}{\alpha}\norm[\opW]{(\opT-\identity)s^k},
        \qquad
        \kappa_\mathrm{res}=\lambdamax(\opM)^{1/2},
    \end{equation}
    where~$\opW$ and~$\opM$ are the matrices of~\Cref{subsec:nonexpansiveness}.
    If, in addition, $\cS^\star$ is nonempty, then~$\norm{\bar r^k}\le(2\kappa_\mathrm{res}/\alpha)\dist_{\opW}(s^k,\cS^\star)$.
\end{lemma}

\begin{proof}
    Set~$\delta s^k=(\opF-\identity)s^k=(\delta x^k,\delta v^k)$, so that~$\delta x^k=\xtkk-x^k$.
    By~\eqref{eq:v-update}, $\delta v^k=A\xtkk-z^k$, which is the second row of~\eqref{eq:kkt-residual-projected}.
    For the first row, $P\xtkk+q+A^Ty^k+w^{k+1}=r_d^k+P\delta x^k$, and~\eqref{eq:inverse-kkt-residual-map} gives~$r_d^k=-(R-\rho A^TA)\delta x^k-\rho A^T\delta v^k$.
    The matrix in~\eqref{eq:kkt-residual-projected} factors as~$J\opW$ with
    \begin{equation*}
        J=\begin{bmatrix}
            \identity & A^T\\
            0 & -\rho^{-1}\identity
        \end{bmatrix},
        \qquad
        J\opW J^T
        =\begin{bmatrix}
            W+\rho A^TA & -A^T\\
            -A & \rho^{-1}\identity
        \end{bmatrix}
        =D\opM D,
    \end{equation*}
    where~$D=\diag(\identity,-\identity)$ is orthogonal.
    Hence~$\norm{\bar r^k}=\norm{J\opW^{1/2}\opW^{1/2}\delta s^k}\le\norm{J\opW^{1/2}}\norm[\opW]{\delta s^k}$ with~$\norm{J\opW^{1/2}}^2=\lambdamax(J\opW J^T)=\lambdamax(\opM)$, which is the inequality in~\eqref{eq:residual-transfer-bound}.
    The equality follows from~$\opT-\identity=\alpha(\opF-\identity)$, which holds since~$\opT=(1-\alpha)\identity+\alpha\opF$.
    For every~$s^\star\in\cS^\star$, the nonexpansiveness of~$\opT$ in the~$\opW$-seminorm (\Cref{cor:lsf-averagedness}($a$)) gives~$\norm[\opW]{(\opT-\identity)s^k}\le\norm[\opW]{\opT s^k-\opT s^\star}+\norm[\opW]{s^\star-s^k}\le2\norm[\opW]{s^k-s^\star}$; taking the infimum over~$s^\star$ gives the last bound.
\end{proof}

\section{Sharpness constants and complexity bounds}
\label{app:complexity}
This appendix proves the claims on explicit rates and lower bounds of~\Cref{subsec:sharpness-free}.
It computes the best sharpness constant when the operator~$\opT$ is affine, compares the resulting complexity with the lower bounds for first-order methods, and shows that no restart rule improves the factor~$\kappa_P$ in the complexity of restarted Halpern iteration on unconstrained strongly convex QPs.

\subsection{Exact sharpness constants}
\label{app:exact-sharpness}
As in~\Cref{sec:sharpness}, we suppose that the strict majorization condition~\eqref{eq:free-strict-majorization} holds.
When~$\cS^\star\neq\emptyset$ and~$\opT\neq\identity$, the \emph{best sharpness constant}
\begin{equation*}
    \eta_{\opW}^\star
    =\inf_{s\notin\cS^\star}\frac{\norm[\opW]{\opT s-s}}{\dist_{\opW}(s,\cS^\star)}
\end{equation*}
is the largest~$\eta$ such that~$\norm[\opW]{\opT s-s}\ge\eta\,\dist_{\opW}(s,\cS^\star)$ for every~$s\in\reals^{d+m}$.
If~$\cX=\reals^d$ and~$\cC$ is an affine subspace, then both projections in~$\opT$ are affine, and so is~$\opT$.
\Cref{thm:affine-sharpness} computes~$\eta_{\opW}^\star$ in this case, and~\Cref{thm:model-sharpness} specializes the result to unconstrained strongly convex QPs and equality-constrained LPs.

\begin{theorem}[Exact sharpness of affine linear-system-free ADMM]
    \label{thm:affine-sharpness}
    Suppose that~$\cX=\reals^d$, that~$\cC$ is an affine subspace, and that~$\cS^\star\neq\emptyset$.
    Write~$\proj_\cC(v)=\Pi v+c$, where~$c=\proj_\cC(0)$ and~$\Pi$ is the orthogonal projector onto the subspace~$\cC-c$, and define
    \begin{equation*}
        S_\Pi
        =
        \begin{bmatrix}
            \identity-W^{1/2}R^{-1}W^{1/2}
            &
            \sqrt\rho\,W^{1/2}R^{-1}A^T(\identity-2\Pi)
            \\
            -\sqrt\rho\,AR^{-1}W^{1/2}
            &
            \rho AR^{-1}A^T(\identity-2\Pi)+\Pi
        \end{bmatrix}.
    \end{equation*}
    If~$S_\Pi=0$, which happens only when~$P=0$, $A=0$, and~$\Pi=0$, then~$\opT=\identity$.
    Otherwise,
    \begin{equation}
        \label{eq:eta-Q-exact}
        \eta_{\opW}^\star
        =\alpha\,\sigma_{\min}^+(S_\Pi),
    \end{equation}
    where~$\sigma_{\min}^+$ denotes the smallest positive singular value.
\end{theorem}

\begin{proof}
    Since~$\cX=\reals^d$, the vector~$w^{k+1}$ in~\eqref{eq:free-normal-vector} vanishes.
    For a state~$s=(x,v)$, let~$r(s)$ be the stacked residual~$(r_d,r_p)$ of~\Cref{lem:free-residual-map} evaluated at~$(x,z,y)$ with~$z=\proj_\cC(v)$ and~$y=\rho(v-z)$.
    Fix~$s^\star=(x^\star,v^\star)\in\cS^\star$, so that~$r(s^\star)=0$, and set~$\Delta s=(\Delta x,\Delta v)=s-s^\star$.
    The relations~$\Delta z=\Pi\Delta v$ and~$\Delta y=\rho(\identity-\Pi)\Delta v$ give
    \begin{equation*}
        r(s)
        =
        N\Delta s,
        \qquad
        N=
        \begin{bmatrix}
            P & \rho A^T(\identity-\Pi)\\
            A & -\Pi
        \end{bmatrix}.
    \end{equation*}
    By~\eqref{eq:free-residual-matrix},
    \begin{equation*}
        (\opT-\identity)s=-\alpha EN\Delta s.
    \end{equation*}
    Direct multiplication using~$W=R-(P+\rho A^TA)$ gives
    \begin{equation*}
        \opW^{1/2}EN\opW^{-1/2}=S_\Pi.
    \end{equation*}
    The upper-left block~$\identity-W^{1/2}R^{-1}W^{1/2}$ vanishes only if~$R=W$, that is, $P+\rho A^TA=0$, which forces~$P=0$ and~$A=0$; then~$S_\Pi=\diag(0,\Pi)$.
    Hence~$S_\Pi=0$ only if~$P=0$, $A=0$, and~$\Pi=0$, and then~$(\opT-\identity)s=0$ for every~$s$.
    Otherwise, set~$h=\opW^{1/2}\Delta s$.
    Then~$\norm[\opW]{\opT s-s}=\alpha\norm{S_\Pi h}$, and since~$U$ is invertible, $\opW^{1/2}(\cS^\star-s^\star)=\nullspace(S_\Pi)$, so~$\dist_{\opW}(s,\cS^\star)=\dist(h,\nullspace(S_\Pi))$.
    The bound~$\norm{S_\Pi h}\ge\sigma_{\min}^+(S_\Pi)\dist(h,\nullspace(S_\Pi))$ holds with equality along a right singular vector for~$\sigma_{\min}^+(S_\Pi)$, which proves~\eqref{eq:eta-Q-exact}.
\end{proof}

\begin{theorem}[Exact sharpness constants]
    \label{thm:model-sharpness}
    Suppose that~$\cX=\reals^d$ and that~$R=\tau^{-1}\identity$ satisfies~\eqref{eq:free-strict-majorization}.
    \begin{enumerate}
        \item[$(i)$]
            If~$P\succ0$ and~$\cC=\reals^m$, so that~\eqref{prob:qp} is an unconstrained strongly convex QP, then~$\eta_{\opW}^\star=\alpha\tau\lambdamin(P)$.
        \item[$(ii)$]
            If~$P=0$, $\cC=\{b\}$, and~$A\neq0$, so that~\eqref{prob:qp} is an equality-constrained LP, and if~$\cS^\star\neq\emptyset$, then~$\eta_{\opW}^\star=\alpha\sqrt{\rho\tau}\,\sigma_{\min}^+(A)$.
    \end{enumerate}
\end{theorem}

\begin{proof}[Proof of~\Cref{thm:model-sharpness}]
    In both cases~$\cC$ is an affine subspace.
    In case~$(i)$, the point~$x=-P^{-1}q$ with~$z=Ax$ and~$y=w=0$ satisfies the KKT conditions~\eqref{eq:qp-kkt}, so~$\cS^\star\neq\emptyset$ by~\Cref{subsec:kkt-residuals}.
    Hence~\Cref{thm:affine-sharpness} applies in both cases, and~$S_\Pi\ne0$ because its upper-left block is nonzero when~$P\ne0$ or~$A\ne0$.

    $(i)$
    The set~$\cC=\reals^m$ gives~$\Pi=\identity$.
    The next steps hold for every~$R$ that satisfies~\eqref{eq:free-strict-majorization}.
    Define
    \begin{equation*}
        B=\begin{bmatrix}W^{1/2}\\ \sqrt\rho\,A\end{bmatrix}.
    \end{equation*}
    Since~$B^TB=R-P$,
    \begin{equation*}
        S_{\identity}=\identity-BR^{-1}B^T,
        \qquad
        S_{\identity}^{-1}=\identity+BP^{-1}B^T.
    \end{equation*}
    Hence~$S_\identity$ is symmetric with~$0\prec S_\identity\preceq\identity$, so~$\sigma_{\min}^+(S_\identity)=\lambdamin(S_\identity)=1/\lambdamax(\identity+BP^{-1}B^T)$.
    The matrices~$BP^{-1}B^T$ and~$P^{-1/2}B^TBP^{-1/2}=P^{-1/2}RP^{-1/2}-\identity$ have the same nonzero eigenvalues, so~\eqref{eq:eta-Q-exact} gives
    \begin{equation}
        \label{eq:eta-unconstrained}
        \eta_{\opW}^\star
        =\frac{\alpha}{\lambdamax(P^{-1/2}RP^{-1/2})}.
    \end{equation}
    For~$R=\tau^{-1}\identity$, $\lambdamax(P^{-1/2}RP^{-1/2})=1/(\tau\lambdamin(P))$, which proves~$(i)$.

    $(ii)$
    The set~$\cC=\{b\}$ gives~$\Pi=0$.
    Let~$A=U\Sigma V^T$ and set~$\gamma_i=\rho\tau\sigma_i^2$ for each singular value~$\sigma_i>0$, so that~$\gamma_i<1$ by~$\tau^{-1}>\rho\norm{A}^2$.
    In~$\opW$-rescaled singular-vector coordinates, the linear part of~$\opT$ acts on the corresponding two-dimensional block as the matrix
    \begin{equation*}
        \begin{bmatrix}
            1-\alpha\gamma_i & -\alpha\sqrt{\gamma_i(1-\gamma_i)}\\
            \alpha\sqrt{\gamma_i(1-\gamma_i)} & 1-\alpha\gamma_i
        \end{bmatrix}.
    \end{equation*}
    This matrix is normal, and its eigenvalues~$\lambda$ satisfy~$|1-\lambda|=\alpha\sqrt{\gamma_i}$.
    Null-space directions are tangent to~$\cS^\star$.
    Formula~\eqref{eq:eta-Q-exact} therefore gives
    \begin{equation*}
        \eta_{\opW}^\star
        =\min_{\sigma_i>0}\alpha\sqrt{\gamma_i}
        =\alpha\sqrt{\rho\tau}\,\sigma_{\min}^+(A).
        \qedhere
    \end{equation*}
\end{proof}

In case~$(i)$ with a general preconditioner~$R$, \eqref{eq:eta-unconstrained} and~$P^{-1/2}RP^{-1/2}\preceq\lambdamax(R)P^{-1}$ give~$\eta_{\opW}^\star\ge\alpha\lambdamin(P)/\lambdamax(R)$.
The strict majorization condition places no upper bound on~$R$, so it alone does not bound~$\eta_{\opW}^\star$ below by a function of~$\kappa_P$ only.
If~$R\preceq c\lambdamax(P)\identity$ for some~$c\ge1$, then~$\eta_{\opW}^\star\ge\alpha/(c\kappa_P)$.

\subsection{Lower bounds for span-respecting methods}
\label{app:lower-bounds}
The lower bounds of Applegate et al.~\cite{applegateFasterFirstorderPrimaldual2023} and of Lu and Yang~\cite{luPracticalOptimalFirstOrder2025} hold for the methods that are span-respecting on the saddle function
\begin{equation*}
    \lagrangian(x,y)
    =(1/2)x^T P x+q^T x+y^T(Ax-b),
\end{equation*}
which is the Lagrangian of~\eqref{prob:qp} with~$\cX=\reals^d$ and~$\cC=\{b\}$.
We show that, on the instances that attain these bounds, restarted linear-system-free ADMM with the scalar preconditioner~$R=\tau^{-1}\identity$ is span-respecting in the following sense, up to a factor of two in the number of evaluations of~$\opT$.

\begin{definition}[Span-respecting primal--dual method, cf.~{\cite[Def.~3]{applegateFasterFirstorderPrimaldual2023}}]
    \label{def:span-respecting}
    A primal--dual method applied to~$\lagrangian$ is \emph{span-respecting} if its iterates~$(x^k, y^k)$ satisfy the following for all~$k \ge 1$:
    \begin{align*}
        x^k &\in x^0+\spann \{\nabla_x\lagrangian(x^i,y^j)\mid 0\le i,j\le k-1\}, \\
        y^k &\in y^0+\spann \{\nabla_y\lagrangian(x^i,y^j)\mid 0\le i \le k,\, 0\le j\le k-1\}.
    \end{align*}
\end{definition}
We call the step from~$k-1$ to~$k$ an \emph{oracle round}: it adds the gradients at the latest iterates to the two spans.

\begin{proposition}[Span-respecting property on the hard families]
    \label{prop:free-span}
    Let~$R=\tau^{-1}\identity$.
    \begin{itemize}
        \item[($a$)]
            On an equality-constrained LP ($P=0$, $\cC=\{b\}$), each evaluation of~$\opT$ can be represented by at most two oracle rounds in the sense of~\Cref{def:span-respecting}.
        \item[($b$)]
            On an unconstrained QP ($A=0$, $b=0$), each evaluation of~$\opT$ is one oracle round.
            Since~$\lagrangian$ does not depend on~$y$, \Cref{def:span-respecting} reduces to the class of first-order methods of Nesterov~\cite[Assumption~2.1.4]{nesterovLecturesConvexOptimization2018}.
    \end{itemize}
    Thus a lower bound on the number of oracle rounds, divided by two, bounds the number of evaluations of~$\opT$ that~\Cref{algo:restarted-halpern} needs on linear-system-free ADMM.
\end{proposition}

\begin{proof}
    On the bilinear family, $P=0$ and~$\cC=\{b\}$.
    A first oracle round queries~$Ax^k-b$ and forms~$\widehat y^k=y^k+\rho(Ax^k-b)$ while leaving the primal point unchanged.
    The second round uses~$q+A^T\widehat y^k$ to form~$\xtkk=x^k-\tau(q+A^T\widehat y^k)$ and the next primal iterate.
    In a Halpern update, the coefficient of~$\xtkk$ in the affine combination defining~$x^{k+1}$ is positive, so~$A\xtkk-b$ is an affine combination of the dual gradients at~$x^{k+1}$, $x^k$, and the anchor.
    At a restart, the new primal iterate is~$\xtkk$ itself.
    Thus two span-respecting rounds suffice.
    If~$A=0$ and~$b=0$, then~$\nabla_y\lagrangian=0$ and~$\nabla_x\lagrangian(x,y)=Px+q$ does not depend on~$y$, so~\Cref{def:span-respecting} requires only~$x^k\in x^0+\spann\{Px^i+q\mid0\le i\le k-1\}$, which is Nesterov's class, and the step reduces to the gradient step~$\tilde x^{k+1}=x^k-\tau(Px^k+q)$, one oracle round.
    Halpern updates and restarts are affine combinations and therefore preserve these spans.
\end{proof}

The following lower bounds therefore apply to linear-system-free ADMM, up to the factor~$2$ of~\Cref{prop:free-span}.
They grow linearly in~$\kappa_A$ on equality-constrained LPs and as~$\sqrt{\kappa_P}$ on unconstrained QPs.

\begin{theorem}[Condition-number lower bounds, {\cite[Corollary~1]{applegateFasterFirstorderPrimaldual2023}, \cite[Theorem~2.1.13]{nesterovLecturesConvexOptimization2018}}]
    \label{thm:lower_bound}
    Let~$\mathcal Z^\star$ be the set of saddle points of~$\lagrangian$, let~$\dist$ be the Euclidean distance, and let~$K\ge1$ and~$\varepsilon>0$.
    \begin{itemize}
        \item[($a$)]
            For every~$\kappa_A>1$ there is an equality-constrained LP, with~$P=0$, $\cC=\{b\}$, and~$\norm{A}/\sigma_{\min}^+(A)\le\kappa_A$, whose dimension depends on~$K$ and~$\kappa_A$, such that every span-respecting primal--dual method started at~$z^0=0$ that reaches~$\dist(z^k,\mathcal Z^\star)\le\varepsilon\dist(z^0,\mathcal Z^\star)$ within~$K$ rounds needs
            \begin{equation*}
                k\ge\frac{\kappa_A-1}{2}\log\frac{1}{\sqrt2\,\varepsilon}
                =\Omega\!\left(\kappa_A\log\frac{1}{\varepsilon}\right)
            \end{equation*}
            rounds.
        \item[($b$)]
            For every~$\kappa_P>1$ there is an unconstrained strongly convex QP, with~$A=0$, $b=0$, and~$\lambdamax(P)/\lambdamin(P)\le\kappa_P$, whose dimension depends on~$K$ and~$\kappa_P$, such that every such method needs
            \begin{equation*}
                k\ge\frac{\sqrt{\kappa_P}-1}{2}\log\frac{1}{\sqrt2\,\varepsilon}
                =\Omega\!\left(\sqrt{\kappa_P}\log\frac{1}{\varepsilon}\right)
            \end{equation*}
            rounds.
    \end{itemize}
    The $\Omega$ bounds hold uniformly when the corresponding condition number is at least~$2$ and~$0<\varepsilon\le1/2$.
\end{theorem}

Both parts are Nesterov's lower bound for the quadratic~$\tfrac12x^THx+h^Tx$ with a tridiagonal~$H$~\cite[Theorem~2.1.13]{nesterovLecturesConvexOptimization2018}, whose factor per round is~$(\sqrt\kappa-1)/(\sqrt\kappa+1)$ with~$\kappa=\lambdamax(H)/\lambdamin(H)$; we state it with this factor, truncated to a finite dimension, which costs the~$\sqrt2$ inside the logarithm.
Part~(a) is the bilinear embedding~$A=H^{1/2}$ of Applegate et al.~\cite[proof of Corollary~1]{applegateFasterFirstorderPrimaldual2023}, for which~$\sqrt\kappa=\kappa_A$.

On equality-constrained LPs, the sharpness constant of~\Cref{thm:model-sharpness}$(ii)$ and~\Cref{thm:restart} give an upper bound of the same order, up to constants and a logarithmic term.

\begin{corollary}[Optimality on equality-constrained LPs]
    \label{cor:lp-optimality}
    If~$\tau^{-1}=c\rho\norm{A}^2$ for fixed~$c>1$ and fixed algorithmic parameters~$\alpha$ and~$\beta$, then restarted linear-system-free ADMM requires at most
    \begin{equation*}
        O\!\left(
            \frac{\norm{A}}{\sigma_{\min}^+(A)}
            \log\frac{\norm{A}}{\sigma_{\min}^+(A)\varepsilon}
        \right)
    \end{equation*}
    evaluations of~$\opT$ to reach~$\dist_{\opW}(s^{n,0},\cS^\star)\le\varepsilon\dist_{\opW}(s^{0,0},\cS^\star)$ for any~$\varepsilon\in(0,1/e]$, with a constant that depends only on~$c$, $\alpha$, and~$\beta$.
    By~\Cref{prop:free-span}($a$), this matches~\Cref{thm:lower_bound}($a$), which counts oracle rounds, up to constants and a logarithmic term.
    Since~$\opW\succ0$, conversion to Euclidean distance changes only the constants inside the logarithm for each fixed problem instance.
\end{corollary}

Lu and Yang~\cite[Remark~1]{luRestartedHalpernPDHG2024} note the same near-optimality for restarted Halpern PDHG on LPs.

\subsection{Tightness of restarted Halpern iteration}
\label{app:tightness}
This subsection shows that the factor~$\kappa_P$ in the complexity of restarted Halpern iteration on unconstrained strongly convex QPs cannot be improved by the restart rule.
Throughout, \eqref{prob:qp} is an unconstrained QP with~$\mu>0$.
The preconditioner is~$R=\tau^{-1}\identity$ and satisfies~\eqref{eq:free-strict-majorization}, so that~$\tau^{-1}>L$, and the relaxation parameter~$\alpha\in(0,2]$ is fixed.

On an affine map, the Halpern iterate~$s^k$ is the average of the first~$k+1$ plain iterates, as Lu and Yang~\cite[Proposition~1]{luRestartedHalpernPDHG2024} observe for PDHG on unconstrained bilinear problems.

\begin{lemma}[Halpern iteration on affine maps; cf.~{\cite[Proposition~1]{luRestartedHalpernPDHG2024}}]
    \label{lem:halpern-ergodic}
    Let~$\opT$ be affine with linear part~$T$ and a fixed point~$s^\star$.
    Then the Halpern iterates~\eqref{eq:halpern} satisfy
    \begin{equation*}
        s^k-s^\star
        =p_k(T)\,(s^0-s^\star),
        \qquad
        p_k(t)=\frac{1}{k+1}\sum_{j=0}^{k}t^j.
    \end{equation*}
\end{lemma}

\begin{proof}
    Write~$\Delta s^k=s^k-s^\star$, so that~$\opT s^k-s^\star=T\Delta s^k$ and~\eqref{eq:halpern} reads~$\Delta s^{k+1}=\frac{k+1}{k+2}T\Delta s^k+\frac{1}{k+2}\Delta s^0$.
    The case~$k=0$ is immediate, and if~$\Delta s^k=\frac{1}{k+1}\sum_{j=0}^{k}T^j\Delta s^0$, then
    \begin{equation*}
        \Delta s^{k+1}
        =\frac{1}{k+2}\sum_{j=0}^{k}T^{j+1}\Delta s^0+\frac{1}{k+2}\Delta s^0
        =\frac{1}{k+2}\sum_{j=0}^{k+1}T^j\Delta s^0.
        \qedhere
    \end{equation*}
\end{proof}

Averaging cannot contract a slow eigenmode faster than the mode itself decays, and restarting only multiplies such averages.
The next theorem quantifies this for every restart schedule.

\begin{theorem}[Tightness for every restart schedule]
    \label{thm:restart-tightness}
    In the setting above, suppose that~$\delta=\max\{\alpha,1\}\,\tau\mu\le1/2$, which holds in particular whenever~$\kappa_P\ge4$, since~$\delta<2/\kappa_P$.
    Then there exists an initialization~$s^{0,0}$ such that every realization of~\Cref{algo:restarted-halpern}, under any restart rule and any choice of outer iteration lengths, satisfies
    \begin{equation}
        \label{eq:restart-tightness}
        \dist_{\opW}(s,\cS^\star)
        \ge e^{-\frac{5}{2}\delta N}\dist_{\opW}(s^{0,0},\cS^\star)
    \end{equation}
    for every iterate~$s$ produced after~$N$ evaluations of~$\opT$.
    Consequently, for~$0<\varepsilon\le1$, reaching~$\dist_{\opW}(s,\cS^\star)\le\varepsilon\,\dist_{\opW}(s^{0,0},\cS^\star)$ requires
    \begin{equation}
        \label{eq:restart-tightness-complexity}
        N
        \ge\frac{2}{5\max\{\alpha,1\}\,\tau\mu}\log\frac1\varepsilon
        \ge\frac{2}{5\max\{\alpha,1\}}\,\frac{L}{\mu}\log\frac1\varepsilon
        =\Omega\!\left(\kappa_P\log\frac1\varepsilon\right).
    \end{equation}
\end{theorem}

\begin{proof}
    On this family~$\proj_\cC$ is the identity, so~$\opT$ is affine, and by the proof of~\Cref{thm:affine-sharpness} its linear part in the rescaled coordinates~$\tilde s=\opW^{1/2}(s-s^\star)$ is~$\identity-\alpha S_\identity$.
    The matrix~$S_\identity$ is symmetric positive definite by the proof of~\Cref{thm:model-sharpness}$(i)$, so the fixed point~$s^\star$ is unique and~$\dist_{\opW}(s,\cS^\star)=\norm{\tilde s}$.
    Comparing~\eqref{eq:eta-Q-exact} with the constant~$\eta_{\opW}^\star=\alpha\tau\mu$ of~\Cref{thm:model-sharpness}$(i)$ gives~$\lambdamin(S_\identity)=\tau\mu$.
    Let~$u$ be a unit eigenvector of~$S_\identity$ for~$\lambdamin(S_\identity)$ and initialize with~$\tilde s^{0,0}=u$, so that the linear part of~$\opT$ acts on~$\spann\{u\}$ as multiplication by~$t=1-\alpha\tau\mu\in[1/2,1)$.
    The linear part of~$\opF$, which is~$\opT$ with~$\alpha=1$, acts on it as multiplication by~$1-\tau\mu$.
    Every step of~\Cref{algo:restarted-halpern}, including a restart from~$\opF s=s+\alpha^{-1}(\opT s-s)$, is an affine combination of~$\opT$-images, iterates, and the anchor of the current outer iteration, so all iterates remain in~$s^\star+\opW^{-1/2}\spann\{u\}$ and the error reduces to a scalar with~$\dist_{\opW}(s,\cS^\star)$ equal to its absolute value.
    In particular, any restart rule, adaptive or not, only selects the outer iteration lengths.

    Consider an outer iteration whose anchor has scalar error~$\xi$ and which restarts after~$k+1$ evaluations of~$\opT$.
    \Cref{lem:halpern-ergodic} gives the inner errors~$\xi\,p_j(t)$ for~$0\le j\le k$, and the next anchor~$s^{n+1,0}=\opF s^{n,k}$ has error~$\xi\,(1-\tau\mu)\,p_k(t)$ after~$k+1$ evaluations of~$\opT$.
    Let~$\delta'=\alpha\tau\mu\le\delta$.
    Since~$t^{k+1}=(1-\delta')^{k+1}\le e^{-(k+1)\delta'}$,
    \begin{equation*}
        p_k(t)
        =\frac{1-t^{k+1}}{(k+1)\delta'}
        \ge\frac{1-e^{-(k+1)\delta'}}{(k+1)\delta'}
        \ge e^{-(k+1)\delta'/2}
        \ge e^{-(k+1)\delta/2},
    \end{equation*}
    where the third step uses~$(1-e^{-a})/a\ge e^{-a/2}$ for~$a>0$, equivalent to~$2\sinh(a/2)\ge a$.
    Together with~$1-\tau\mu\ge e^{-2\tau\mu}\ge e^{-2\delta}$ for~$\tau\mu\le\delta\le1/2$, one outer iteration with~$k+1$ evaluations multiplies the error by at least
    \begin{equation*}
        (1-\tau\mu)\,p_k(t)
        \ge e^{-2\delta}\,e^{-(k+1)\delta/2}
        \ge e^{-\frac{5}{2}\delta(k+1)},
    \end{equation*}
    using~$k\ge0$ in the last inequality.
    Multiplying these factors over completed outer iterations, and bounding the factor of the current outer iteration by~$p_k(t)\ge e^{-(k+1)\delta/2}\ge e^{-\frac{5}{2}\delta k}$ for~$k\ge1$, proves~\eqref{eq:restart-tightness} for every iterate produced after~$N$ total evaluations.
    Solving~$e^{-\frac{5}{2}\delta N}\le\varepsilon$ for~$N$ and using~$\tau<1/L$ yields~\eqref{eq:restart-tightness-complexity}.
\end{proof}

The proof requires only that the initial error have a nonzero component along the eigenvector~$u$.
For a general~$s^{0,0}$, the same argument applies to this component, so~\eqref{eq:restart-tightness} holds with~$\dist_{\opW}(s^{0,0},\cS^\star)$ replaced by the magnitude of this component.

The plain iteration~$s^{k+1}=\opT s^k$ attains the same order on these problems.
In the rescaled coordinates~$\tilde s=\opW^{1/2}(s-s^\star)$,
its linear part~$\identity-\alpha S_\identity$ is symmetric with eigenvalues in~$[1-\alpha,\,1-\alpha\tau\mu]$, because~$\tau\mu\identity\preceq S_\identity\preceq\identity$ by the proofs of~\Cref{thm:model-sharpness}$(i)$ and~\Cref{thm:restart-tightness}.
For~$\alpha\le1$, it therefore satisfies~$\|s^{k+1}-s^\star\|_{\opW}\le(1-\alpha\tau\mu)\|s^k-s^\star\|_{\opW}$ and needs~$O\big((\alpha\tau\mu)^{-1}\log(1/\varepsilon)\big)$ evaluations of~$\opT$, which is~$O(\kappa_P\log(1/\varepsilon))$ for~$\tau^{-1}=\Theta(L)$.
Hence the Halpern update does not improve the worst-case order of complexity of the plain iteration on these problems.

\section{Proofs for infeasible problems}
\label{app:infeasibility}
This appendix proves~\Cref{thm:minimal-displacement} and~\Cref{prop:candidate-convergence} of~\Cref{subsec:infeasible}.
Throughout, the strict majorization condition~\eqref{eq:free-strict-majorization} holds, so that~$\opM\succ0$.
The proofs use the sets~\eqref{eq:gap-sets}, which we write as
\begin{equation*}
    \cG_\mathrm{p}=A\cX-\cC,
    \qquad
    \cG_\mathrm{d}=q+\range(P)+\barc\cX+A^T\barc\cC.
\end{equation*}
They also use that the barrier cone of a box is the polar of its recession cone, \ie, $\barc\cC=\{ y \mid y^Te\le0 \text{ for all } e\in\rec\cC \}$, and likewise for~$\cX$.
By~\Cref{subsec:infeasible},
\begin{equation}
    \label{eq:feasibility-gap-sets}
    \text{\eqref{prob:qp} is primal feasible}\iff0\in\cG_\mathrm{p},
    \qquad
    \text{\eqref{prob:qp} is dual feasible}\iff0\in\cG_\mathrm{d}.
\end{equation}

\begin{lemma}[Range of the displacement]
    \label{lem:range-displacement}
    The range of~$\opTCV-\identity$ is~$-\opM^{-1}\range(\opA)$.
\end{lemma}

\begin{proof}
    The operator~$\opM+\opA$ is strongly monotone, because~$\opA$ is the sum of a subdifferential and a skew-symmetric linear map and~$\opM\succ0$.
    Hence~$(\opM+\opA)^{-1}$ is single-valued, and~$u^+=\opTCV u$ if and only if~$\opM(u-u^+)\in\opA u^+$.
    This shows~$\opTCV u-u\in-\opM^{-1}\range(\opA)$.
    Conversely, let~$a\in\opA u'$ and set~$u=u'+\opM^{-1}a$.
    Then~$\opM u\in(\opM+\opA)u'$, so~$\opTCV u=u'$ and~$\opTCV u-u=-\opM^{-1}a$.
\end{proof}

\begin{lemma}[Range of the saddle operator]
    \label{lem:range-saddle}
    The sets~\eqref{eq:gap-sets} are nonempty, convex, and closed, and the range of~$\opA$ is~$\cG_\mathrm{d}\times(-\cG_\mathrm{p})$.
\end{lemma}

\begin{proof}
    Both sets are polyhedral, hence closed.
    Indeed, $A\cX$ is the linear image of a polyhedron, hence polyhedral~\cite[Theorem~19.3]{rockafellarConvexAnalysis1970}, so~$\cG_\mathrm{p}$ is the sum of two polyhedra, which is polyhedral~\cite[Corollary~19.3.2]{rockafellarConvexAnalysis1970}.
    The set~$\cG_\mathrm{d}$ is a translate of the sum of a subspace and two finitely generated cones, which is again finitely generated.

    Fix~$(r,c')\in\reals^d\times\reals^m$.
    The inclusion~$(r,c')\in\opA(x,y)$ holds if and only if there exists~$w\in\reals^d$ such that
    \begin{equation*}
        Px+(q-r)+A^Ty+w=0,
        \qquad
        w\in\cN_\cX(x),
        \qquad
        Ax+c'\in\cC,
        \qquad
        y\in\cN_\cC(Ax+c').
    \end{equation*}
    With~$z=Ax$, these are the KKT conditions of the shifted problem
    \begin{equation*}
        \begin{array}{ll}
            \textnormal{minimize} & (1/2)x^TPx+(q-r)^Tx\\
            \textnormal{subject to} & l-c'\le Ax\le u-c',\\
            & l_x\le x\le u_x,
        \end{array}
    \end{equation*}
    where~$x\in\reals^d$ is the optimization variable.
    Its constraint box~$\cC-c'$ has the same recession cone as~$\cC$.
    By~\eqref{eq:feasibility-gap-sets} applied to the shifted data, the shifted problem is primal feasible if and only if~$-c'\in\cG_\mathrm{p}$, and dual feasible if and only if~$r\in\cG_\mathrm{d}$.
    A KKT point implies primal and dual feasibility.
    Conversely, under primal and dual feasibility, weak duality bounds the objective from below on the nonempty polyhedral feasible set, so the Frank--Wolfe theorem~\cite{frankAlgorithmQuadraticProgramming1956} gives a minimizer, at which the KKT conditions hold because the constraints are polyhedral.
    Hence~$(r,c')\in\range(\opA)$ if and only if~$(r,c')\in\cG_\mathrm{d}\times(-\cG_\mathrm{p})$.
\end{proof}

By~\Cref{lem:range-displacement,lem:range-saddle},
\begin{equation}
    \label{eq:range-formula}
    \range(\opTCV-\identity)=-\opM^{-1}\left(\cG_\mathrm{d}\times(-\cG_\mathrm{p})\right),
\end{equation}
which is closed and convex.
Jiang, Moursi, and Vavasis~\cite{jiangRangeDisplacementOperator2026} characterize the range of the displacement operator of PDHG in a similar way.

\begin{proof}[Proof of~\Cref{thm:minimal-displacement}]
    The proof has three steps.

    \emph{Step 1: sign conditions.}
    The projection~$r_p^\star$ is characterized by~$r_p^\star\in\cG_\mathrm{p}$ and
    \begin{equation}
        \label{eq:proj-p}
        (r_p^\star)^T(r'-r_p^\star)\ge0
        \quad\text{for all } r'\in\cG_\mathrm{p}.
    \end{equation}
    Since~$\cG_\mathrm{p}-\rec\cC\subseteq\cG_\mathrm{p}$ and~$\cG_\mathrm{p}+A\rec\cX\subseteq\cG_\mathrm{p}$, taking~$r'=r_p^\star-e$ and~$r'=r_p^\star+Ae'$ in~\eqref{eq:proj-p} gives~$(r_p^\star)^Te\le0$ for~$e\in\rec\cC$ and~$(r_p^\star)^TAe'\ge0$ for~$e'\in\rec\cX$.
    Hence~$\delta y^\star\in\barc\cC$ and~$-A^T\delta y^\star\in\barc\cX$.
    Similarly, $r_d^\star\in\cG_\mathrm{d}$ and
    \begin{equation}
        \label{eq:proj-d}
        (R^{-1}r_d^\star)^T(r-r_d^\star)\ge0
        \quad\text{for all } r\in\cG_\mathrm{d}.
    \end{equation}
    The set~$\cG_\mathrm{d}$ is invariant under adding elements of~$\range(P)$, $\barc\cX$, and~$A^T\barc\cC$.
    Since~$\range(P)$ is a subspace, \eqref{eq:proj-d} gives~$P\delta x^\star=0$.
    It also gives~$(\delta x^\star)^Tw\le0$ for~$w\in\barc\cX$ and~$(A\delta x^\star)^Ty\le0$ for~$y\in\barc\cC$, so~$\delta x^\star\in\rec\cX$ and~$A\delta x^\star\in\rec\cC$ by the bipolar theorem.
    In particular, $(A\delta x^\star)^T\delta y^\star\le0$.

    \emph{Step 2: minimality of~$(\delta x^\star,\delta y^\star)$.}
    By~\eqref{eq:range-formula}, $(\delta x^\star,\delta y^\star)$ is an element of~$\range(\opTCV-\identity)$ if and only if~$-\opM(\delta x^\star,\delta y^\star)\in\cG_\mathrm{d}\times(-\cG_\mathrm{p})$.
    Using~$P\delta x^\star=0$ and~$R\delta x^\star=-r_d^\star$, the two blocks of~$-\opM(\delta x^\star,\delta y^\star)$ are
    \begin{equation*}
        r_d^\star-A^T\delta y^\star\in\cG_\mathrm{d}+\barc\cX\subseteq\cG_\mathrm{d},
        \qquad
        -(r_p^\star+A\delta x^\star)\in-(\cG_\mathrm{p}+A\rec\cX)\subseteq-\cG_\mathrm{p},
    \end{equation*}
    by Step~1.
    The element of minimal~$\opM$-norm of the closed convex set~\eqref{eq:range-formula} is its unique element~$h^\star$ with~$(h^\star)^T\opM(h-h^\star)\ge0$ for all elements~$h$~\cite[Theorem~3.16]{bauschkeConvexAnalysisMonotone2017}.
    Writing~$h=-\opM^{-1}a$ with~$a\in\cG_\mathrm{d}\times(-\cG_\mathrm{p})$, this condition reads
    \begin{equation}
        \label{eq:variational}
        (h^\star)^T(\opM h^\star+a)\le0
        \quad\text{for all } a\in\cG_\mathrm{d}\times(-\cG_\mathrm{p}).
    \end{equation}
    Let~$h^\star=(\delta x^\star,\delta y^\star)$ and~$a=(r,-r')$ with~$r\in\cG_\mathrm{d}$ and~$r'\in\cG_\mathrm{p}$.
    Since~$P\delta x^\star=0$, $(h^\star)^T\opM h^\star=\norm[R]{\delta x^\star}^2+2(A\delta x^\star)^T\delta y^\star+\rho^{-1}\norm{\delta y^\star}^2$.
    By~\eqref{eq:proj-d}, $(\delta x^\star)^Tr\le-\norm[R^{-1}]{r_d^\star}^2=-\norm[R]{\delta x^\star}^2$, and by~\eqref{eq:proj-p}, $-(\delta y^\star)^Tr'\le-\rho\norm{r_p^\star}^2=-\rho^{-1}\norm{\delta y^\star}^2$.
    Adding the three relations gives
    \begin{equation}
        \label{eq:variational-bound}
        (h^\star)^T(\opM h^\star+a)\le2(A\delta x^\star)^T\delta y^\star\le0,
    \end{equation}
    so~$h^\star$ satisfies~\eqref{eq:variational}, which proves~($d$).

    \emph{Step 3: identities and equivalences.}
    Choosing~$a=-\opM h^\star$ in~\eqref{eq:variational-bound}, which is allowed by Step~2, gives~$0\le2(A\delta x^\star)^T\delta y^\star\le0$.
    Together with~$P\delta x^\star=0$, this proves~($c$).
    For~\eqref{eq:primal-identity}, the definition of~$\cG_\mathrm{p}$ and~\eqref{eq:proj-p} give
    \begin{equation*}
        S_\cC(\delta y^\star)+S_\cX(-A^T\delta y^\star)
        =\sup_{r'\in\cG_\mathrm{p}}\left(-(\delta y^\star)^Tr'\right)
        =-\rho\inf_{r'\in\cG_\mathrm{p}}(r_p^\star)^Tr'
        =-\rho\norm{r_p^\star}^2
        =-\rho^{-1}\norm{\delta y^\star}^2.
    \end{equation*}
    For~\eqref{eq:dual-identity}, \eqref{eq:proj-d} gives~$\sup_{r\in\cG_\mathrm{d}}(\delta x^\star)^Tr=-\norm[R^{-1}]{r_d^\star}^2=-\norm[R]{\delta x^\star}^2$.
    By Step~1, this supremum equals~$q^T\delta x^\star$, because~$\delta x^\star$ is orthogonal to~$\range(P)$ and has a nonpositive inner product with the two cones in~$\cG_\mathrm{d}$, with equality at~$0$.
    For the last three statements of the theorem, \eqref{eq:closed-form} and~\eqref{eq:feasibility-gap-sets} give~$\delta y^\star=0$ if and only if the problem is primal feasible, and~$\delta x^\star=0$ if and only if it is dual feasible.
    If~$\delta y^\star\ne0$, then~\eqref{eq:primal-identity} gives~\eqref{eq:primal-certificate}, and if~$\delta x^\star\ne0$, then~($a$) gives~\eqref{eq:dual-certificate}.
    Finally, $(\delta x^\star,\delta y^\star)$ is the element of minimal~$\opM$-norm in~$\range(\opTCV-\identity)$ by~($d$), so~$(\delta x^\star,\delta y^\star)=0$ if and only if~$0\in\range(\opTCV-\identity)$, \ie, if and only if~$\opTCV$ has a fixed point.
    The fixed points of~$\opTCV$ are the zeros of~$\opA$.
    These are the pairs~$(x,y)$ of the KKT points~$(x,z,y,w)$ of~\eqref{prob:qp}, with~$z=Ax$ and~$w=-(Px+q+A^Ty)$, and~$\Phi$ maps them onto~$\cS^\star$ by~\Cref{lem:free-residual-map}.
\end{proof}

\begin{proof}[Proof of~\Cref{prop:candidate-convergence}]
    Since~$\opT-\identity=\alpha(\opF-\identity)=\alpha\,\Phi\circ(\opTCV-\identity)\circ\Phi^{-1}$ and~$\Phi$ is a linear isometry, the element of~$\range(\opT-\identity)$ with minimal~$\opW$-norm is~$\alpha\Phi(\delta x^\star,\delta y^\star)$.
    Park and Ryu~\cite{parkAcceleratedInfeasibilityDetection2023} call the element of~$\range(\identity-\opT)$ with minimal norm the infimal displacement vector~$v$ of~$\opT$, so here~$v=-\alpha\Phi(\delta x^\star,\delta y^\star)$, with the opposite orientation to the certificate candidates.
    By the same isometry,
    \begin{equation*}
        \norm[\opM]{(\delta x^k,\delta y^k)-(\delta x^\star,\delta y^\star)}
        =\alpha^{-1}\norm[\opW]{(\identity-\opT)s^k-v},
        \qquad
        \norm[\opM]{(\delta x^k,\delta y^k)}
        =\alpha^{-1}\norm[\opW]{(\identity-\opT)s^k}.
    \end{equation*}
    The set~$\cS_\delta=\{ s \mid (\opF-\identity)s=\Phi(\delta x^\star,\delta y^\star) \}$ is nonempty because~$\Phi(\delta x^\star,\delta y^\star)\in\range(\opF-\identity)$, and every~$s_\delta\in\cS_\delta$ satisfies~$s_\delta-\opT s_\delta=v$.
    We apply the results of~\cite{parkAcceleratedInfeasibilityDetection2023} with the inner product~$(a,b)\mapsto a^T\opW b$.
    In part~($a$), the bound holds for every~$s_\delta\in\cS_\delta$, so we take the infimum over~$s_\delta$.

    For~($a$), the operator~$\opN=2\opF-\identity$ is nonexpansive in~$\norm[\opW]{\cdot}$, because~$\opF$ is firmly nonexpansive by~\Cref{cor:lsf-averagedness}($a$).
    Since~$\opT=(1-\alpha/2)\identity+(\alpha/2)\opN$, the plain iteration is the Krasnosel'ski\u{\i}--Mann iteration of~\cite{parkAcceleratedInfeasibilityDetection2023} applied to~$\opN$ with the constant parameter~$\lambda_k=1-\alpha/2\in(0,1)$.
    For every~$s$, $s-\opN s=(2/\alpha)(s-\opT s)$, so the infimal displacement vector of~$\opN$ is~$(2/\alpha)v$, attained on~$\cS_\delta$.
    Hence~\cite[Theorem~4]{parkAcceleratedInfeasibilityDetection2023}, with~$\sum_{i=0}^k\lambda_{i+1}(1-\lambda_{i+1})=(k+1)\alpha(2-\alpha)/4$, gives, for~$k\ge1$,
    \begin{equation*}
        \left(\norm[\opM]{(\delta x^k,\delta y^k)}-\norm[\opM]{(\delta x^\star,\delta y^\star)}\right)^2
        \le\frac{\dist_{\opW}(s^0,\cS_\delta)^2}{\alpha(2-\alpha)(k+1)}.
    \end{equation*}
    The certificate candidates lie in~$\range(\opTCV-\identity)$, whose element of minimal~$\opM$-norm is~$(\delta x^\star,\delta y^\star)$, so the proof of~\cite[Lemma~2]{parkAcceleratedInfeasibilityDetection2023} gives
    \begin{equation*}
        \norm[\opM]{(\delta x^k,\delta y^k)-(\delta x^\star,\delta y^\star)}^2
        \le\norm[\opM]{(\delta x^k,\delta y^k)}^2-\norm[\opM]{(\delta x^\star,\delta y^\star)}^2,
    \end{equation*}
    whose right-hand side tends to zero by the previous bound, so the certificate candidates converge to~$(\delta x^\star,\delta y^\star)$.
    Convergence alone also follows from the Baillon--Bruck--Reich theorem~\cite{baillonAsymptoticBehaviorNonexpansive1978}.

    For~($b$), $\opT$ is nonexpansive in~$\norm[\opW]{\cdot}$ by~\Cref{cor:lsf-averagedness}($a$), and~\cite[Theorem~8]{parkAcceleratedInfeasibilityDetection2023} gives a point~$\bar s\in\cS_\delta$ such that, for~$k\ge1$,
    \begin{equation*}
        \norm[\opM]{(\delta x^k,\delta y^k)-(\delta x^\star,\delta y^\star)}^2
        \le\frac{\left(\sqrt{H_k+4}+1\right)^2}{\alpha^2(k+1)^2}\,
        \norm[\opW]{s^0-\bar s}^2,
        \qquad
        H_k=\sum_{i=1}^k\frac{1}{i}\le1+\log k,
    \end{equation*}
    whose right-hand side is~$O(\log k/k^2)$ and hence tends to zero.
\end{proof}

\section{Algebra engine}
\label{app:algebra-engine}
This appendix describes the algebra engine of~\Cref{sec:implementation}, which evaluates the sparse products of~\Cref{algo:batched-h-osqp} for all~$B$ columns in parallel.

\paragraph{Two products per iteration}
As written, \Cref{algo:batched-h-osqp} needs the products~$PX$, $AX$, $A\bar X$, and one product with~$A^T$.
The engine needs two: the stacked product~$\left[\begin{smallmatrix}A\\ P\end{smallmatrix}\right]X$, which gives~$AX$ and~$PX$ in one pass over the matrix, and the product~$A^TY$.
To avoid the product~$A\bar X$, the engine stores~$\Psi=\brho\otimes(V-AX)$ in place of~$V$.
Each column of~$(X,\Psi)$ is then the pair~$\Phi^{-1}(x,v)$ of~\Cref{thm:lsf-osqp-condat-vu}($c$), and the primal--dual form~\eqref{eq:lsf-primal-dual} gives
\begin{equation*}
    \bar X=\proj_{\cX}\big(X-R^{-1}\otimes\left[PX+Q+2A^TY-A^T\Psi\right]\big),
    \qquad
    \bar\Psi=Y.
\end{equation*}
The relaxation and the acceleration scheme are affine in the state, so the engine updates the product~$A^T\Psi$ with the same affine combination, without another sparse product.
For example, for a Halpern column with counter~$k$ and anchor~$\Psi^0$, the new product is
\begin{equation*}
    \frac{1}{k+2}A^T\Psi^0+\frac{k+1}{k+2}\left[(1-\alpha)A^T\Psi+\alpha A^TY\right].
\end{equation*}

\paragraph{Two implementations}
For every product, the engine selects between an implementation using the cuSPARSE library~\cite{nvidiaCuSPARSELibrary2019} and custom kernels written for the batched iteration.
Both evaluate the same formulas, so the choice changes the speed and the rounding of an iteration, but not the iteration itself.
The custom kernels tailor their behavior to the sparsity structure detected in the matrices during setup and perform their computations in a deterministic way, providing bitwise identical iterates from run to run, which the cuSPARSE algorithms do not guarantee.

\paragraph{Rows of different lengths}
The number of nonzeros per row varies widely, so the custom kernels split the rows of each matrix into three classes and assign threads differently in each.
A \emph{light} row, with at most~$32$ nonzeros, is processed by one thread per column.
These~$B$ threads are adjacent in a warp, the group of~$32$ threads that the GPU executes in lockstep, so each nonzero is loaded once for all of them.
A \emph{mid} row, with at most~$4096$ nonzeros, is processed by one warp for a group of columns.
Its threads read every~$32$nd nonzero, apply it to all columns of the group, and add their partial sums along a fixed binary tree.
A \emph{heavy} row is split into chunks of about~$8192/B$ nonzeros, and each chunk is processed by one thread block for all~$B$ columns.
The partial sums of the chunks are added in a fixed order.
For small batches, or when the long rows are too few to occupy the GPU, the selection below can change the number of threads per row and the boundaries between the classes.

\paragraph{Selection at setup}
The engine is configured once per problem, at setup, after the batch size has been chosen.
For each product, it first times the cuSPARSE algorithms against each other, then tunes the custom kernels over the options in~\Cref{tab:engine-choices}, and finally keeps the faster of the two, with ties resolved in favor of cuSPARSE.
Each option is timed over a fixed number of repetitions after warm-up runs.
The two products of an iteration are timed in alternation, because a product timed alone keeps its data in the L2 cache and appears faster than it is inside the iteration.
The incumbent option is only changed if the new option is faster by more than~$3\%$, which keeps timing noise from changing a choice.

\begin{table}
    \centering
    \small
    \caption{Options that the engine compares for each product at setup.
    The last column marks the choices that leave the iterates bitwise unchanged.
    The other choices change the order of the floating-point operations.
    The light/mid boundary~$128$ is tried only for~$B\le2$.}
    \label{tab:engine-choices}
    \begin{tabular}{llll}
        \toprule
        Choice & Options & Applies to & Bitwise \\
        \midrule
        Implementation & custom, cuSPARSE & every product & \xmark \\
        cuSPARSE algorithm & CSR algorithms of cuSPARSE & cuSPARSE products & \xmark \\
        Columns per warp & $1,2,4,8,16$ & mid rows & \cmark \\
        Threads per block & $128$, $256$ & custom kernels & \cmark \\
        Threads per light row & $1,4,8,16$ & $B\le16$ & \xmark \\
        Light/mid boundary & $8$, $32$, $128$ nonzeros & $B\le16$ & \xmark \\
        Mid/heavy boundary & $256$, $1024$, $4096$ nonzeros & underfilled GPU & \xmark \\
        Chunk length & $8192/B$, sized to fill the GPU & heavy rows & \xmark \\
        Heavy-row kernel & after, alongside the other rows & first product & \cmark \\
        Index loads & explicit, implicit & consecutive columns & \cmark \\
        Single-chunk rows & partial-sum slot, direct update & heavy rows & \cmark \\
        \bottomrule
    \end{tabular}
\end{table}

\section{Portfolio problem formulations}
\label{app:portfolio}
This appendix gives the portfolio problem of the backtest of~\Cref{sec:parametric-results}, the generation of its data, and the estimation of the risk model.

\paragraph{Problem formulation}
On each trading day, we solve the portfolio optimization problem
\begin{equation*}
    \begin{array}{ll}
        \textnormal{maximize} & \mu^T x - \gamma x^T \Sigma x \\
        \textnormal{subject to} & \ones^T x = 1, \quad x \geq 0,
    \end{array}
\end{equation*}
where~$x \in \reals^n$ is the optimization variable, which represents the portfolio.
The problem data are the vector of expected returns~$\mu \in \reals^n$, the risk-aversion parameter~$\gamma > 0$, and the covariance matrix~$\Sigma \in \symm_+^n$ of the risk model.
The~$20$ values of~$\gamma$ are logarithmically spaced over the interval~$[10^{-2}, 10^{2}]$.
The risk model is the sum of a matrix of rank~$k < n$ and a diagonal matrix,
\begin{equation*}
    \Sigma = F F^T + D,
\end{equation*}
where~$F \in \reals^{n \times k}$ is the factor loading matrix and the diagonal matrix~$D \in \reals^{n \times n}$ describes the asset-specific risk.
With this model, the problem is equivalent to
\begin{equation*}
    \begin{array}{ll}
        \textnormal{minimize} & x^T D x + y^T y - \gamma^{-1} \mu^T x \\
        \textnormal{subject to} & y = F^T x, \quad \ones^T x = 1, \quad x \geq 0,
    \end{array}
\end{equation*}
where~$x \in \reals^n$ and~$y \in \reals^k$ are the optimization variables.
The parameter~$\gamma$ enters only the linear term of the objective.
We also consider the \emph{full covariance} model, in which~$\Sigma$ is the sample covariance matrix of the returns, with no factors.

\paragraph{Data generation}
We take from Yahoo Finance the daily returns of the 3,000 stocks listed on the NYSE and the NASDAQ with the largest market capitalization, over the three years from January~1, 2023, to January~1, 2026.
A stock that was not listed before 2023 is replaced by the stock with the next largest market capitalization that was.
Hence every trading day of the three years has the returns of exactly 3,000 stocks.
The set with 1,500 stocks is built in the same way from the largest stocks.
To generate four years of synthetic returns, we fit a two-tailed generalized Pareto distribution to the returns of every stock by the peaks-over-threshold method.
The largest and the smallest~$15$ percent of the returns of a stock form its two tails, whose parameters we estimate by maximum likelihood, and the empirical quantiles describe the returns between the two thresholds.
The stocks are coupled through a Gaussian copula with the sample correlation of the historical returns.
The synthetic days are drawn independently of one another.

\paragraph{Risk model estimation}
We estimate~$\mu$, $F$, and~$D$ from the returns of the previous~$h$ trading days.
The vector~$\mu$ is the mean of these returns.
Let~$\lambda_1 \geq \cdots \geq \lambda_n$ be the eigenvalues of their sample covariance matrix~$\hat\Sigma$, with orthonormal eigenvectors~$v_1, \ldots, v_n$.
We keep the~$k$ leading principal components and set
\begin{equation*}
    F = \begin{bmatrix} \sqrt{\lambda_1}\, v_1 & \cdots & \sqrt{\lambda_k}\, v_k \end{bmatrix},
    \qquad
    D_{ii} = \hat\Sigma_{ii} - (FF^T)_{ii}, \quad i = 1, \ldots, n.
\end{equation*}

\paragraph{Rolling window}
The window has~$h = 3k$ days, and~$h_0$ denotes the number of historical trading days.
For the full covariance model, the window has~$h = h_0$ days.
The backtest reads the last~$\min\{h, h_0\}$ historical days followed by the synthetic days, and builds problem~$i$ from days~$i+1, \ldots, i+h$ of this sequence.
If~$h \leq h_0$, the first window is historical.
Otherwise, it holds every historical day and~$h - h_0$ synthetic days.
Every later problem adds one synthetic day and drops its oldest day.
We re-estimate~$\mu$ every day, and~$F$ and~$D$ every 60 trading days, \ie, every three months.
Each new estimate is blended entrywise with the previous one as~$\mu \gets (1-\varphi)\mu^\mathrm{new} + \varphi\mu$, and likewise for~$F$ and~$D$, with~$\varphi = 0.7$.

\Cref{app:detailed-results} reports the total solve times for all the numbers of stocks and risk models that we consider.

\section{Perturbed instances for the variant order}
\label{app:variant-order}
This appendix describes the perturbed copies of the MIPLIB, Mittelmann, Maros--M\'esz\'aros, and SuiteSparse instances from which we compute the variant order of~\Cref{subsec:variant-order}.

\paragraph{Perturbation}
Let~$\gamma>0$ be the noise level and~$\phi\in(0,1]$ the fraction of perturbed entries.
In each of~$A$, $q$, $l$, $u$, $l_x$, and~$u_x$, we select a random fraction~$\phi$ of the finite nonzero entries and replace each selected entry~$a$ by
\begin{equation*}
    a'=a(1+\gamma\chi),
    \qquad \chi\sim\cN(0,1),
\end{equation*}
with independent samples~$\chi$.
To keep the objective convex, we perturb~$P$ as
\begin{equation*}
    P'=DPD,
    \qquad
    D=\diag(e^{\gamma\chi_1},\ldots,e^{\gamma\chi_d}),
\end{equation*}
where~$\chi_j\sim\cN(0,1)$ for a random fraction~$\phi$ of the indices~$j$ and~$\chi_j=0$ otherwise.
A variable whose two bounds are equal stays fixed, and if the perturbed bounds of a variable cross, we set both to their midpoint.
We denote the perturbed data by~$A'$, $q'$, $\cC'=[l',u']$, and~$\cX'=[l_x',u_x']$.

\paragraph{Repair}
The perturbed instance may be infeasible or unbounded, so we modify~$l'$, $u'$, and~$q'$ until it has a solution.
Let~$(x^\star,y^\star,w^\star)$ be a primal--dual solution of the original instance, which we compute with Gurobi~\cite{gurobioptimizationllcGurobiOptimizerReference2026} with tolerance~$10^{-6}$ and a time limit of 30 seconds.
We apply the following steps in order, for every row~$i$:
\begin{align*}
    x^\circ&=\proj_{\cX'}(x^\star), \\
    l_i'&\gets\min\{l_i',(A'x^\circ)_i\},
    \qquad
    u_i'\gets\max\{u_i',(A'x^\circ)_i\}, \\
    y^\circ&=\proj_{\barc\cC'}(y^\star),
    \qquad
    w^\circ=\proj_{\barc\cX'}\big(-q'-P'x^\circ-A'^Ty^\circ\big), \\
    q'&\gets-(P'x^\circ+A'^Ty^\circ+w^\circ).
\end{align*}
For an equality row we set~$l_i'=u_i'=(A'x^\circ)_i$ in the second step.
The first two steps make~$x^\circ$ feasible for the perturbed instance.
The last two make~$(x^\circ,y^\circ,w^\circ)$ feasible for its dual, where~$\barc$ denotes the barrier cone of~\Cref{sec:lsf-osqp}.
The last step changes~$q_j'$ only if~$x_j$ lacks a finite lower or upper bound.
By weak duality the objective is then bounded below on a nonempty feasible set, so the perturbed instance has an optimal solution~\cite{frankAlgorithmQuadraticProgramming1956}.
We use~$\phi=0.7$ and~$\gamma=0.1$, and we generate~$3$ perturbed copies of every instance.

\section{Detailed benchmark results}
\label{app:detailed-results}
This appendix reports the results of every solver on the benchmark sets of~\Cref{sec:qp-benchmarks,sec:lp-benchmarks} and on the portfolio backtest of~\Cref{sec:parametric-results}.

\paragraph{QP and LP benchmarks}
\Cref{tab:maros-sgm,tab:suite-sgm,tab:qplib-sgm,tab:Mittelmann-LP-sgm,tab:MIPLIB-sgm} report the results by size group.
For each set, the first subtable has the original instances and the second the instances obtained after Gurobi's presolve routine.
In each table, \# is the number of solved instances and SGM is the shifted geometric mean of the solve times in seconds.
The size of each group appears in parentheses.
Highlighted entries mark the lowest time in each column, separately for the interior-point and for the first-order methods.
Gurobi runs on the CPU and every other solver on the GPU.
The result of Clarabel on the presolved QPLIB set (\Cref{tab:qplib-sgm-sidebyside-presolved}) varies between runs: it fails on one instance in about two of five runs, at either tolerance.

\begin{table}[!htp]
    \centering
    \caption{Maros--M\'esz\'aros: solved instances (\#) and shifted geometric mean of the solve times in seconds (SGM), for tolerances $\varepsilon \in \{10^{-4}, 10^{-8}\}$, by size group.}
    \label{tab:maros-sgm}

    \begin{subtable}{\linewidth}
    \centering
    \caption{Original instances.}
    \label{tab:maros-sgm-sidebyside}
    \footnotesize
    \setlength{\tabcolsep}{4pt}
    \begin{tabular}{lllrrrrrrrr}
        \toprule
        \multirow{2}{*}{$\varepsilon$} & \multirow{2}{*}{Method} & \multirow{2}{*}{Solver} & \multicolumn{2}{c}{Small (83)} & \multicolumn{2}{c}{Medium (47)} & \multicolumn{2}{c}{Large (8)} & \multicolumn{2}{c}{Total (138)}\\
        \cmidrule(lr){4-5}\cmidrule(lr){6-7}\cmidrule(lr){8-9}\cmidrule(lr){10-11}
         & & & \# & SGM & \# & SGM & \# & SGM & \# & SGM\\
        \midrule
        \multirow{6}{*}{$10^{-4}$} & \multirow{3}{*}{Interior point} & Clarabel & 83 & 0.08 & 46 & 1.50 & 6 & 36.00 & 135 & 1.51\\
         &  & cuOpt-IP & 83 & \best{0.07} & 46 & 1.51 & 8 & \best{0.68} & 137 & \best{0.58}\\
         &  & Gurobi-IP & 82 & 0.76 & 47 & \best{0.46} & 8 & 2.96 & 137 & 0.77\\
        \cmidrule(lr){2-11}
         & \multirow{3}{*}{First order} & HPR-QP & 79 & 4.44 & 43 & 19.95 & 8 & 6.97 & 130 & 8.69\\
         &  & \osqplf{} & 83 & \best{0.42} & 47 & \best{1.00} & 8 & \best{3.14} & 138 & \best{0.76}\\
         &  & SCS & 81 & 7.50 & 40 & 19.93 & 6 & 48.75 & 127 & 12.54\\
        \midrule
        \multirow{6}{*}{$10^{-8}$} & \multirow{3}{*}{Interior point} & Clarabel & 80 & 2.47 & 38 & 21.36 & 6 & 36.27 & 124 & 8.42\\
         &  & cuOpt-IP & 83 & \best{0.09} & 46 & \best{1.57} & 8 & \best{0.74} & 137 & \best{0.61}\\
         &  & Gurobi-IP & 81 & 1.55 & 43 & 7.09 & 8 & 2.99 & 132 & 3.29\\
        \cmidrule(lr){2-11}
         & \multirow{3}{*}{First order} & HPR-QP & 78 & 5.76 & 40 & 30.98 & 8 & 20.71 & 126 & 12.69\\
         &  & \osqplf{} & 83 & \best{0.99} & 47 & \best{16.99} & 8 & \best{10.73} & 138 & \best{5.48}\\
         &  & SCS & 78 & 10.38 & 35 & 84.22 & 6 & 92.45 & 119 & 27.70\\
        \bottomrule
    \end{tabular}
    \end{subtable}

    \vspace{1.5em}
    \begin{subtable}{\linewidth}
    \centering
    \caption{Presolved instances.}
    \label{tab:maros-sgm-sidebyside-presolved}
    \footnotesize
    \setlength{\tabcolsep}{4pt}
    \begin{tabular}{lllrrrrrrrr}
        \toprule
        \multirow{2}{*}{$\varepsilon$} & \multirow{2}{*}{Method} & \multirow{2}{*}{Solver} & \multicolumn{2}{c}{Small (83)} & \multicolumn{2}{c}{Medium (47)} & \multicolumn{2}{c}{Large (8)} & \multicolumn{2}{c}{Total (138)}\\
        \cmidrule(lr){4-5}\cmidrule(lr){6-7}\cmidrule(lr){8-9}\cmidrule(lr){10-11}
         & & & \# & SGM & \# & SGM & \# & SGM & \# & SGM\\
        \midrule
        \multirow{6}{*}{$10^{-4}$} & \multirow{3}{*}{Interior point} & Clarabel & 83 & 0.08 & 46 & 1.51 & 6 & 35.85 & 135 & 1.51\\
         &  & cuOpt-IP & 83 & \best{0.06} & 45 & 3.00 & 8 & \best{0.57} & 136 & 1.01\\
         &  & Gurobi-IP & 82 & 0.75 & 47 & \best{0.48} & 8 & 2.52 & 137 & \best{0.75}\\
        \cmidrule(lr){2-11}
         & \multirow{3}{*}{First order} & HPR-QP & 80 & 3.46 & 43 & 19.50 & 8 & \best{5.57} & 131 & 7.73\\
         &  & \osqplf{} & 83 & \best{0.31} & 47 & \best{1.15} & 8 & 5.59 & 138 & \best{0.84}\\
         &  & SCS & 79 & 6.55 & 42 & 10.85 & 6 & 86.77 & 127 & 9.83\\
        \midrule
        \multirow{6}{*}{$10^{-8}$} & \multirow{3}{*}{Interior point} & Clarabel & 81 & 1.62 & 39 & 17.65 & 6 & 36.21 & 126 & 6.91\\
         &  & cuOpt-IP & 83 & \best{0.08} & 45 & \best{3.06} & 8 & \best{0.64} & 136 & \best{1.04}\\
         &  & Gurobi-IP & 82 & 0.75 & 45 & 3.35 & 7 & 15.81 & 134 & 2.18\\
        \cmidrule(lr){2-11}
         & \multirow{3}{*}{First order} & HPR-QP & 79 & 4.31 & 40 & 29.83 & 8 & 20.65 & 127 & 11.19\\
         &  & \osqplf{} & 83 & \best{0.57} & 47 & \best{16.46} & 8 & \best{12.68} & 138 & \best{5.10}\\
         &  & SCS & 78 & 8.77 & 35 & 51.75 & 5 & 155.60 & 118 & 21.94\\
        \bottomrule
    \end{tabular}
    \end{subtable}
\end{table}

\begin{table}[!htp]
    \centering
    \caption{SuiteSparse: solved instances (\#) and shifted geometric mean of the solve times in seconds (SGM), for tolerances $\varepsilon \in \{10^{-4}, 10^{-8}\}$, by size group.}
    \label{tab:suite-sgm}

    \begin{subtable}{\linewidth}
    \centering
    \caption{Original instances.}
    \label{tab:suite-sgm-sidebyside}
    \footnotesize
    \setlength{\tabcolsep}{4pt}
    \begin{tabular}{lllrrrrrrrr}
        \toprule
        \multirow{2}{*}{$\varepsilon$} & \multirow{2}{*}{Method} & \multirow{2}{*}{Solver} & \multicolumn{2}{c}{Small (22)} & \multicolumn{2}{c}{Medium (8)} & \multicolumn{2}{c}{Large (22)} & \multicolumn{2}{c}{Total (52)}\\
        \cmidrule(lr){4-5}\cmidrule(lr){6-7}\cmidrule(lr){8-9}\cmidrule(lr){10-11}
         & & & \# & SGM & \# & SGM & \# & SGM & \# & SGM\\
        \midrule
        \multirow{6}{*}{$10^{-4}$} & \multirow{3}{*}{Interior point} & Clarabel & 22 & 0.05 & 8 & 0.10 & 22 & \best{10.37} & 52 & \best{3.56}\\
         &  & cuOpt-IP & 21 & 3.18 & 6 & 34.16 & 14 & 147.67 & 41 & 35.36\\
         &  & Gurobi-IP & 22 & \best{0.02} & 8 & \best{0.09} & 19 & 63.41 & 49 & 13.30\\
        \cmidrule(lr){2-11}
         & \multirow{3}{*}{First order} & HPR-QP & 22 & 0.13 & 8 & 0.26 & 21 & 11.44 & 51 & 3.94\\
         &  & \osqplf{} & 22 & \best{0.09} & 8 & \best{0.07} & 22 & \best{3.96} & 52 & \best{1.57}\\
         &  & SCS & 22 & 2.72 & 8 & 11.48 & 20 & 49.49 & 50 & 16.48\\
        \midrule
        \multirow{6}{*}{$10^{-8}$} & \multirow{3}{*}{Interior point} & Clarabel & 22 & 0.06 & 8 & 0.11 & 20 & \best{21.19} & 50 & \best{6.25}\\
         &  & cuOpt-IP & 21 & 3.24 & 6 & 34.29 & 13 & 184.49 & 40 & 39.70\\
         &  & Gurobi-IP & 22 & \best{0.03} & 8 & \best{0.10} & 13 & 313.91 & 43 & 33.67\\
        \cmidrule(lr){2-11}
         & \multirow{3}{*}{First order} & HPR-QP & 22 & 0.82 & 8 & 1.05 & 20 & 25.28 & 50 & 7.90\\
         &  & \osqplf{} & 22 & \best{0.31} & 8 & \best{0.28} & 20 & \best{15.88} & 50 & \best{5.21}\\
         &  & SCS & 22 & 2.23 & 8 & 9.57 & 17 & 116.85 & 47 & 25.36\\
        \bottomrule
    \end{tabular}
    \end{subtable}

    \vspace{1.5em}
    \begin{subtable}{\linewidth}
    \centering
    \caption{Presolved instances.}
    \label{tab:suite-sgm-sidebyside-presolved}
    \footnotesize
    \setlength{\tabcolsep}{4pt}
    \begin{tabular}{lllrrrrrrrr}
        \toprule
        \multirow{2}{*}{$\varepsilon$} & \multirow{2}{*}{Method} & \multirow{2}{*}{Solver} & \multicolumn{2}{c}{Small (22)} & \multicolumn{2}{c}{Medium (8)} & \multicolumn{2}{c}{Large (19)} & \multicolumn{2}{c}{Total (49)}\\
        \cmidrule(lr){4-5}\cmidrule(lr){6-7}\cmidrule(lr){8-9}\cmidrule(lr){10-11}
         & & & \# & SGM & \# & SGM & \# & SGM & \# & SGM\\
        \midrule
        \multirow{6}{*}{$10^{-4}$} & \multirow{3}{*}{Interior point} & Clarabel & 22 & 0.05 & 8 & 0.11 & 19 & \best{9.27} & 49 & \best{2.95}\\
         &  & cuOpt-IP & 22 & 0.07 & 8 & 0.12 & 12 & 141.57 & 42 & 18.84\\
         &  & Gurobi-IP & 22 & \best{0.02} & 8 & \best{0.09} & 17 & 56.02 & 47 & 10.84\\
        \cmidrule(lr){2-11}
         & \multirow{3}{*}{First order} & HPR-QP & 22 & 0.19 & 8 & 0.29 & 18 & 11.73 & 48 & 3.69\\
         &  & \osqplf{} & 22 & \best{0.07} & 8 & \best{0.08} & 19 & \best{5.86} & 49 & \best{2.01}\\
         &  & SCS & 22 & 2.17 & 8 & 13.81 & 18 & 26.93 & 48 & 10.88\\
        \midrule
        \multirow{6}{*}{$10^{-8}$} & \multirow{3}{*}{Interior point} & Clarabel & 22 & 0.07 & 7 & 11.08 & 14 & \best{65.17} & 43 & \best{14.77}\\
         &  & cuOpt-IP & 21 & 3.20 & 6 & 34.10 & 10 & 273.01 & 37 & 42.77\\
         &  & Gurobi-IP & 22 & \best{0.02} & 7 & \best{11.03} & 12 & 235.94 & 41 & 29.12\\
        \cmidrule(lr){2-11}
         & \multirow{3}{*}{First order} & HPR-QP & 22 & 0.64 & 7 & 15.92 & 14 & 59.22 & 43 & 15.43\\
         &  & \osqplf{} & 22 & \best{0.26} & 8 & \best{6.65} & 16 & \best{41.56} & 46 & \best{10.77}\\
         &  & SCS & 21 & 3.26 & 7 & 25.22 & 16 & 113.87 & 44 & 26.99\\
        \bottomrule
    \end{tabular}
    \end{subtable}
\end{table}

\begin{table}[!htp]
    \centering
    \caption{QPLIB: solved instances (\#) and shifted geometric mean of the solve times in seconds (SGM), for tolerances $\varepsilon \in \{10^{-4}, 10^{-8}\}$, by size group.}
    \label{tab:qplib-sgm}

    \begin{subtable}{\linewidth}
    \centering
    \caption{Original instances.}
    \label{tab:qplib-sgm-sidebyside}
    \footnotesize
    \setlength{\tabcolsep}{4pt}
    \begin{tabular}{lllrrrrrrrr}
        \toprule
        \multirow{2}{*}{$\varepsilon$} & \multirow{2}{*}{Method} & \multirow{2}{*}{Solver} & \multicolumn{2}{c}{Small (2)} & \multicolumn{2}{c}{Medium (21)} & \multicolumn{2}{c}{Large (13)} & \multicolumn{2}{c}{Total (36)}\\
        \cmidrule(lr){4-5}\cmidrule(lr){6-7}\cmidrule(lr){8-9}\cmidrule(lr){10-11}
         & & & \# & SGM & \# & SGM & \# & SGM & \# & SGM\\
        \midrule
        \multirow{6}{*}{$10^{-4}$} & \multirow{3}{*}{Interior point} & Clarabel & 2 & 0.07 & 19 & 7.86 & 11 & 18.67 & 32 & 10.52\\
         &  & cuOpt-IP & 2 & 0.08 & 19 & 8.00 & 11 & 23.61 & 32 & 11.84\\
         &  & Gurobi-IP & 2 & \best{0.05} & 20 & \best{4.06} & 13 & \best{8.82} & 35 & \best{5.33}\\
        \cmidrule(lr){2-11}
         & \multirow{3}{*}{First order} & HPR-QP & 2 & \best{0.07} & 21 & 1.62 & 13 & \best{8.01} & 36 & \best{3.51}\\
         &  & \osqplf{} & 2 & \best{0.07} & 21 & \best{0.52} & 12 & 17.60 & 35 & 4.86\\
         &  & SCS & 2 & 0.10 & 19 & 18.69 & 12 & 26.57 & 33 & 19.56\\
        \midrule
        \multirow{6}{*}{$10^{-8}$} & \multirow{3}{*}{Interior point} & Clarabel & 2 & 0.08 & 19 & \best{7.90} & 9 & 58.98 & 30 & \best{18.22}\\
         &  & cuOpt-IP & 2 & 0.09 & 17 & 21.56 & 7 & 165.18 & 26 & 45.00\\
         &  & Gurobi-IP & 2 & \best{0.05} & 18 & 13.98 & 10 & \best{37.65} & 30 & 19.28\\
        \cmidrule(lr){2-11}
         & \multirow{3}{*}{First order} & HPR-QP & 2 & 0.14 & 20 & 6.26 & 11 & 64.54 & 33 & 17.45\\
         &  & \osqplf{} & 2 & \best{0.07} & 21 & \best{4.64} & 11 & \best{47.05} & 34 & \best{13.43}\\
         &  & SCS & 2 & 0.12 & 19 & 22.60 & 10 & 129.04 & 31 & 41.58\\
        \bottomrule
    \end{tabular}
    \end{subtable}

    \vspace{1.5em}
    \begin{subtable}{\linewidth}
    \centering
    \caption{Presolved instances.}
    \label{tab:qplib-sgm-sidebyside-presolved}
    \footnotesize
    \setlength{\tabcolsep}{4pt}
    \begin{tabular}{lllrrrrrrrr}
        \toprule
        \multirow{2}{*}{$\varepsilon$} & \multirow{2}{*}{Method} & \multirow{2}{*}{Solver} & \multicolumn{2}{c}{Small (2)} & \multicolumn{2}{c}{Medium (21)} & \multicolumn{2}{c}{Large (13)} & \multicolumn{2}{c}{Total (36)}\\
        \cmidrule(lr){4-5}\cmidrule(lr){6-7}\cmidrule(lr){8-9}\cmidrule(lr){10-11}
         & & & \# & SGM & \# & SGM & \# & SGM & \# & SGM\\
        \midrule
        \multirow{6}{*}{$10^{-4}$} & \multirow{3}{*}{Interior point} & Clarabel & 2 & 0.07 & 19 & 7.83 & 10 & 33.23 & 31 & 13.78\\
         &  & cuOpt-IP & 2 & 0.09 & 18 & 13.59 & 11 & 23.20 & 31 & 15.46\\
         &  & Gurobi-IP & 2 & \best{0.06} & 20 & \best{3.83} & 13 & \best{8.24} & 35 & \best{5.02}\\
        \cmidrule(lr){2-11}
         & \multirow{3}{*}{First order} & HPR-QP & 2 & 0.11 & 21 & 2.31 & 12 & 31.75 & 35 & 8.92\\
         &  & \osqplf{} & 2 & \best{0.06} & 21 & \best{0.55} & 12 & \best{26.35} & 35 & \best{6.45}\\
         &  & SCS & 2 & \best{0.06} & 20 & 7.09 & 11 & 52.86 & 33 & 16.56\\
        \midrule
        \multirow{6}{*}{$10^{-8}$} & \multirow{3}{*}{Interior point} & Clarabel & 2 & 0.08 & 19 & \best{7.87} & 11 & \best{18.43} & 32 & \best{10.47}\\
         &  & cuOpt-IP & 2 & 0.10 & 18 & 13.66 & 7 & 151.96 & 27 & 35.20\\
         &  & Gurobi-IP & 2 & \best{0.05} & 18 & 13.86 & 9 & 58.34 & 29 & 23.25\\
        \cmidrule(lr){2-11}
         & \multirow{3}{*}{First order} & HPR-QP & 2 & 0.12 & 20 & 5.86 & 10 & 76.74 & 32 & 18.57\\
         &  & \osqplf{} & 2 & \best{0.07} & 21 & \best{4.76} & 11 & \best{47.02} & 34 & \best{13.54}\\
         &  & SCS & 2 & 0.08 & 19 & 13.23 & 10 & 117.15 & 31 & 30.97\\
        \bottomrule
    \end{tabular}
    \end{subtable}
\end{table}

\begin{table}[!htp]
    \centering
    \caption{Mittelmann: solved instances (\#) and shifted geometric mean of the solve times in seconds (SGM), for tolerances $\varepsilon \in \{10^{-4}, 10^{-8}\}$, by size group.}
    \label{tab:Mittelmann-LP-sgm}

    \begin{subtable}{\linewidth}
    \centering
    \caption{Original instances.}
    \label{tab:Mittelmann-LP-sgm-sidebyside}
    \footnotesize
    \setlength{\tabcolsep}{4pt}
    \begin{tabular}{lllrrrrrrrr}
        \toprule
        \multirow{2}{*}{$\varepsilon$} & \multirow{2}{*}{Method} & \multirow{2}{*}{Solver} & \multicolumn{2}{c}{Small (11)} & \multicolumn{2}{c}{Medium (27)} & \multicolumn{2}{c}{Large (11)} & \multicolumn{2}{c}{Total (49)}\\
        \cmidrule(lr){4-5}\cmidrule(lr){6-7}\cmidrule(lr){8-9}\cmidrule(lr){10-11}
         & & & \# & SGM & \# & SGM & \# & SGM & \# & SGM\\
        \midrule
        \multirow{7}{*}{$10^{-4}$} & \multirow{2}{*}{Interior point} & Clarabel & 9 & 46.50 & 22 & 88.10 & 8 & 1001.06 & 39 & 136.33\\
         &  & Gurobi-IP & 11 & \best{6.95} & 24 & \best{53.85} & 7 & \best{745.02} & 42 & \best{72.55}\\
        \cmidrule(lr){2-11}
         & \multirow{5}{*}{First order} & cuOpt-PDLP & 11 & 5.54 & 25 & 26.82 & 11 & 52.44 & 47 & 24.16\\
         &  & cuPDLPx & 11 & \best{4.65} & 27 & 10.10 & 10 & 89.35 & 48 & \best{16.80}\\
         &  & Gurobi-PDHG & 10 & 25.51 & 24 & 61.39 & 10 & 185.45 & 44 & 66.52\\
         &  & HPR-LP & 11 & 6.42 & 26 & 18.18 & 11 & \best{41.60} & 48 & 18.59\\
         &  & \osqplf{} & 10 & 17.14 & 27 & \best{6.83} & 11 & 83.87 & 48 & 17.55\\
        \midrule
        \multirow{7}{*}{$10^{-8}$} & \multirow{2}{*}{Interior point} & Clarabel & 8 & 100.10 & 18 & 272.47 & 5 & 1924.07 & 31 & 342.11\\
         &  & Gurobi-IP & 10 & \best{25.26} & 21 & \best{116.42} & 7 & \best{990.50} & 38 & \best{141.01}\\
        \cmidrule(lr){2-11}
         & \multirow{5}{*}{First order} & cuOpt-PDLP & 10 & 40.40 & 24 & 86.21 & 10 & 233.75 & 44 & 92.52\\
         &  & cuPDLPx & 10 & 29.06 & 24 & 48.56 & 9 & 256.37 & 43 & 65.12\\
         &  & Gurobi-PDHG & 7 & 240.35 & 15 & 610.38 & 4 & 3382.66 & 26 & 731.02\\
         &  & HPR-LP & 10 & \best{28.17} & 24 & 52.90 & 10 & \best{130.00} & 44 & \best{57.29}\\
         &  & \osqplf{} & 10 & 31.91 & 25 & \best{35.08} & 10 & 286.96 & 45 & 57.72\\
        \bottomrule
    \end{tabular}
    \end{subtable}

    \vspace{0.5em}
    \begin{subtable}{\linewidth}
    \centering
    \caption{Presolved instances.}
    \label{tab:Mittelmann-LP-sgm-sidebyside-presolved}
    \footnotesize
    \setlength{\tabcolsep}{4pt}
    \begin{tabular}{lllrrrrrrrr}
        \toprule
        \multirow{2}{*}{$\varepsilon$} & \multirow{2}{*}{Method} & \multirow{2}{*}{Solver} & \multicolumn{2}{c}{Small (11)} & \multicolumn{2}{c}{Medium (27)} & \multicolumn{2}{c}{Large (11)} & \multicolumn{2}{c}{Total (49)}\\
        \cmidrule(lr){4-5}\cmidrule(lr){6-7}\cmidrule(lr){8-9}\cmidrule(lr){10-11}
         & & & \# & SGM & \# & SGM & \# & SGM & \# & SGM\\
        \midrule
        \multirow{7}{*}{$10^{-4}$} & \multirow{2}{*}{Interior point} & Clarabel & 9 & 40.94 & 24 & 36.23 & 7 & 968.73 & 40 & 83.75\\
         &  & Gurobi-IP & 11 & \best{6.85} & 25 & \best{26.28} & 10 & \best{115.93} & 46 & \best{30.39}\\
        \cmidrule(lr){2-11}
         & \multirow{5}{*}{First order} & cuOpt-PDLP & 11 & 3.87 & 27 & 10.50 & 11 & 12.11 & 49 & 9.10\\
         &  & cuPDLPx & 11 & \best{1.24} & 27 & 6.72 & 10 & 35.47 & 48 & 9.15\\
         &  & Gurobi-PDHG & 10 & 19.33 & 26 & 29.57 & 10 & 77.25 & 46 & 34.19\\
         &  & HPR-LP & 11 & 1.27 & 27 & 8.84 & 11 & \best{9.76} & 49 & \best{6.97}\\
         &  & \osqplf{} & 11 & 1.92 & 27 & \best{5.63} & 11 & 24.55 & 49 & 7.57\\
        \midrule
        \multirow{7}{*}{$10^{-8}$} & \multirow{2}{*}{Interior point} & Clarabel & 9 & 42.25 & 20 & 124.27 & 6 & 1351.94 & 35 & 172.75\\
         &  & Gurobi-IP & 11 & \best{7.15} & 25 & \best{27.57} & 10 & \best{140.61} & 46 & \best{33.03}\\
        \cmidrule(lr){2-11}
         & \multirow{5}{*}{First order} & cuOpt-PDLP & 11 & 17.40 & 26 & 33.74 & 10 & 182.05 & 47 & 44.89\\
         &  & cuPDLPx & 10 & 21.71 & 24 & 42.35 & 10 & 114.28 & 44 & 46.80\\
         &  & Gurobi-PDHG & 6 & 621.05 & 16 & 418.96 & 3 & 6297.55 & 25 & 845.32\\
         &  & HPR-LP & 11 & 11.46 & 27 & 28.02 & 11 & \best{73.31} & 49 & 29.88\\
         &  & \osqplf{} & 11 & \best{6.49} & 26 & \best{21.35} & 11 & 123.40 & 48 & \best{27.57}\\
        \bottomrule
    \end{tabular}
    \end{subtable}
\end{table}

\begin{table}[!htp]
    \centering
    \caption{MIPLIB: solved instances (\#) and shifted geometric mean of the solve times in seconds (SGM), for tolerances $\varepsilon \in \{10^{-4}, 10^{-8}\}$, by size group.}
    \label{tab:MIPLIB-sgm}

    \begin{subtable}{\linewidth}
    \centering
    \caption{Original instances.}
    \label{tab:MIPLIB-sgm-sidebyside}
    \footnotesize
    \setlength{\tabcolsep}{4pt}
    \begin{tabular}{lllrrrrrrrr}
        \toprule
        \multirow{2}{*}{$\varepsilon$} & \multirow{2}{*}{Method} & \multirow{2}{*}{Solver} & \multicolumn{2}{c}{Small (845)} & \multicolumn{2}{c}{Medium (123)} & \multicolumn{2}{c}{Large (24)} & \multicolumn{2}{c}{Total (992)}\\
        \cmidrule(lr){4-5}\cmidrule(lr){6-7}\cmidrule(lr){8-9}\cmidrule(lr){10-11}
         & & & \# & SGM & \# & SGM & \# & SGM & \# & SGM\\
        \midrule
        \multirow{7}{*}{$10^{-4}$} & \multirow{2}{*}{Interior point} & Clarabel & 803 & 4.45 & 106 & 38.88 & 21 & 369.88 & 930 & 8.19\\
         &  & Gurobi-IP & 833 & \best{2.14} & 115 & \best{33.52} & 19 & \best{354.15} & 967 & \best{5.44}\\
        \cmidrule(lr){2-11}
         & \multirow{5}{*}{First order} & cuOpt-PDLP & 835 & 2.99 & 119 & 9.04 & 21 & 234.77 & 975 & 4.62\\
         &  & cuPDLPx & 836 & 1.94 & 121 & 6.21 & 22 & 127.51 & 979 & 3.15\\
         &  & Gurobi-PDHG & 821 & 4.34 & 118 & 16.00 & 21 & 400.14 & 960 & 6.74\\
         &  & HPR-LP & 831 & 3.34 & 122 & 7.40 & 22 & \best{107.15} & 975 & 4.54\\
         &  & \osqplf{} & 842 & \best{0.84} & 122 & \best{5.55} & 22 & 137.25 & 986 & \best{2.07}\\
        \midrule
        \multirow{7}{*}{$10^{-8}$} & \multirow{2}{*}{Interior point} & Clarabel & 755 & 10.12 & 86 & 98.89 & 13 & 1385.15 & 854 & 17.48\\
         &  & Gurobi-IP & 801 & \best{5.13} & 108 & \best{50.23} & 17 & \best{528.68} & 926 & \best{9.58}\\
        \cmidrule(lr){2-11}
         & \multirow{5}{*}{First order} & cuOpt-PDLP & 812 & 6.90 & 113 & 23.98 & 16 & 1000.69 & 941 & 10.34\\
         &  & cuPDLPx & 815 & 4.69 & 114 & 20.21 & 20 & \best{287.79} & 949 & 7.28\\
         &  & Gurobi-PDHG & 692 & 26.88 & 80 & 139.39 & 9 & 3771.43 & 781 & 39.06\\
         &  & HPR-LP & 820 & 5.91 & 117 & 21.02 & 21 & 307.70 & 958 & 8.58\\
         &  & \osqplf{} & 834 & \best{2.24} & 118 & \best{14.89} & 19 & 546.69 & 971 & \best{4.66}\\
        \bottomrule
    \end{tabular}
    \end{subtable}

    \vspace{0.5em}
    \begin{subtable}{\linewidth}
    \centering
    \caption{Presolved instances.}
    \label{tab:MIPLIB-sgm-sidebyside-presolved}
    \footnotesize
    \setlength{\tabcolsep}{4pt}
    \begin{tabular}{lllrrrrrrrr}
        \toprule
        \multirow{2}{*}{$\varepsilon$} & \multirow{2}{*}{Method} & \multirow{2}{*}{Solver} & \multicolumn{2}{c}{Small (827)} & \multicolumn{2}{c}{Medium (122)} & \multicolumn{2}{c}{Large (24)} & \multicolumn{2}{c}{Total (973)}\\
        \cmidrule(lr){4-5}\cmidrule(lr){6-7}\cmidrule(lr){8-9}\cmidrule(lr){10-11}
         & & & \# & SGM & \# & SGM & \# & SGM & \# & SGM\\
        \midrule
        \multirow{7}{*}{$10^{-4}$} & \multirow{2}{*}{Interior point} & Clarabel & 792 & 3.57 & 105 & 32.82 & 20 & 115.47 & 917 & 6.56\\
         &  & Gurobi-IP & 819 & \best{1.48} & 118 & \best{19.55} & 22 & \best{65.37} & 959 & \best{3.53}\\
        \cmidrule(lr){2-11}
         & \multirow{5}{*}{First order} & cuOpt-PDLP & 822 & 1.48 & 121 & 5.45 & 23 & 33.78 & 966 & 2.32\\
         &  & cuPDLPx & 820 & 1.40 & 122 & \best{3.06} & 23 & 29.68 & 965 & 1.96\\
         &  & Gurobi-PDHG & 813 & 2.41 & 117 & 13.00 & 23 & 62.92 & 953 & 4.01\\
         &  & HPR-LP & 821 & 1.26 & 120 & 5.11 & 22 & \best{26.72} & 963 & 2.03\\
         &  & \osqplf{} & 826 & \best{0.47} & 121 & 4.99 & 23 & 32.30 & 970 & \best{1.34}\\
        \midrule
        \multirow{7}{*}{$10^{-8}$} & \multirow{2}{*}{Interior point} & Clarabel & 758 & 7.34 & 96 & 52.88 & 20 & \best{121.91} & 874 & 11.43\\
         &  & Gurobi-IP & 794 & \best{3.67} & 113 & \best{28.60} & 18 & 165.11 & 925 & \best{6.58}\\
        \cmidrule(lr){2-11}
         & \multirow{5}{*}{First order} & cuOpt-PDLP & 814 & 3.74 & 119 & 17.45 & 23 & 101.91 & 956 & 5.78\\
         &  & cuPDLPx & 806 & 3.18 & 114 & 16.10 & 23 & \best{56.19} & 943 & 4.94\\
         &  & Gurobi-PDHG & 639 & 33.78 & 73 & 188.27 & 14 & 476.77 & 726 & 46.15\\
         &  & HPR-LP & 814 & 3.51 & 119 & 15.95 & 23 & 76.45 & 956 & 5.35\\
         &  & \osqplf{} & 823 & \best{1.24} & 116 & \best{12.25} & 23 & 94.42 & 962 & \best{2.94}\\
        \bottomrule
    \end{tabular}
    \end{subtable}
\end{table}

\paragraph{Portfolio backtest}
\Cref{tab:portfolio-time} reports the total solve times on the portfolio backtest with~$1{,}500$ and~$3{,}000$ stocks for every risk model that we consider.
Gurobi runs on the CPU, and all other solvers on the GPU.
\osqplf{} runs in the three ways of~\Cref{sec:parametric-results}.

\begin{table}[!ht]
    \centering
    \caption{Portfolio backtest: total solve time in seconds of the 19,200 QPs, for tolerances $\varepsilon \in \{10^{-4}, 10^{-8}\}$, with risk models of $k \in \{50, 75, 100\}$ factors and with the full covariance matrix.}
    \label{tab:portfolio-time}

    \begin{subtable}{\linewidth}
    \centering
    \caption{$n = 1{,}500$ stocks.}
    \label{tab:portfolio-time-1500}
    \footnotesize
    \setlength{\tabcolsep}{4pt}
    \begin{tabular}{lllrrrr}
        \toprule
         & & & \multicolumn{4}{c}{Total solve time (s)}\\
        \cmidrule(lr){4-7}
        $\varepsilon$ & Method & Solver & $k = 50$ & $k = 75$ & $k = 100$ & Full\\
        \midrule
        \multirow{6}{*}{$10^{-4}$} & \multirow{3}{*}{Interior point} & Clarabel & 1459.73 & 1528.90 & \best{1550.73} & 7596.16\\
         &  & cuOpt-IP & \best{1231.18} & \best{1419.15} & 1590.16 & \best{3246.60}\\
         &  & Gurobi-IP & 4947.55 & 1627.46 & 2017.00 & 10249.89\\
        \cmidrule(lr){2-7}
         & \multirow{3}{*}{\osqplf{}} & Cold start & 1011.91 & 678.54 & 729.25 & 1169.62\\
         &  & Warm start & 47.03 & 40.25 & 36.86 & 57.39\\
         &  & Batched, warm start & \best{13.38} & \best{14.57} & \best{14.53} & \best{29.91}\\
        \midrule
        \multirow{6}{*}{$10^{-8}$} & \multirow{3}{*}{Interior point} & Clarabel & 1872.53 & 1922.81 & \best{1929.15} & 9353.79\\
         &  & cuOpt-IP & \best{1451.11} & \best{1778.38} & 2192.15 & \best{3945.48}\\
         &  & Gurobi-IP & 4983.16 & 5292.46 & 2180.29 & 11705.68\\
        \cmidrule(lr){2-7}
         & \multirow{3}{*}{\osqplf{}} & Cold start & 1293.36 & 995.81 & 1069.54 & 1779.12\\
         &  & Warm start & 280.44 & 274.30 & 283.01 & 508.23\\
         &  & Batched, warm start & \best{38.88} & \best{58.49} & \best{58.68} & \best{285.42}\\
        \bottomrule
    \end{tabular}
    \end{subtable}

    \vspace{1.5em}
    \begin{subtable}{\linewidth}
    \centering
    \caption{$n = 3{,}000$ stocks.}
    \label{tab:portfolio-time-3000}
    \footnotesize
    \setlength{\tabcolsep}{4pt}
    \begin{tabular}{lllrrrr}
        \toprule
         & & & \multicolumn{4}{c}{Total solve time (s)}\\
        \cmidrule(lr){4-7}
        $\varepsilon$ & Method & Solver & $k = 50$ & $k = 75$ & $k = 100$ & Full\\
        \midrule
        \multirow{6}{*}{$10^{-4}$} & \multirow{3}{*}{Interior point} & Clarabel & 1978.29 & 2172.51 & 2320.56 & 17327.96\\
         &  & cuOpt-IP & \best{1401.61} & \best{1746.99} & \best{2191.26} & \best{10571.94}\\
         &  & Gurobi-IP & 2269.02 & 2946.57 & 3713.61 & 29073.13\\
        \cmidrule(lr){2-7}
         & \multirow{3}{*}{\osqplf{}} & Cold start & 1512.57 & 923.63 & 1022.48 & 3317.90\\
         &  & Warm start & 73.84 & 61.35 & 71.74 & 207.88\\
         &  & Batched, warm start & \best{19.13} & \best{21.27} & \best{24.47} & \best{148.20}\\
        \midrule
        \multirow{6}{*}{$10^{-8}$} & \multirow{3}{*}{Interior point} & Clarabel & 2286.68 & \best{2490.01} & \best{2658.26} & 19290.29\\
         &  & cuOpt-IP & 2373.03 & 2618.06 & 3725.33 & \best{14407.90}\\
         &  & Gurobi-IP & \best{2258.76} & 2980.32 & 3770.84 & 29876.27\\
        \cmidrule(lr){2-7}
         & \multirow{3}{*}{\osqplf{}} & Cold start & 1914.21 & 1328.03 & 1546.09 & 5857.66\\
         &  & Warm start & 409.19 & 387.12 & 474.43 & 1933.31\\
         &  & Batched, warm start & \best{73.11} & \best{86.65} & \best{95.39} & \best{1516.41}\\
        \bottomrule
    \end{tabular}
    \end{subtable}
\end{table}

\end{document}